\documentclass[12pt,a4paper]{amsart}
\usepackage{amsmath,amsfonts,amsthm,amssymb,paralist,pstricks,pst-plot,pst-node}
\usepackage[margin=3cm,twoside]{geometry}

\newcommand{\genmat}{\lambda}

\newcommand{\vnorm}[1]{\left\|  #1 \right\|}

\newcommand{\Cn}{\mathbb{C}^n}
\newcommand{\Cd}{\mathbb{C}^d}
\newcommand{\Cm}{\mathbb{C}^m}
\newcommand{\C}{\mathbb{C}}
\newcommand{\CN}{\mathbb{C}^N}

\newcommand{\Rd}{\mathbb{R}^d}

\newcommand{\R}{\mathbb{R}}

\newcommand{\N}{\mathbb{N}}

\newcommand{\onenorm}[1]{{\left\| #1 \right\|}_{1}}
\newcommand{\deltanorm}[1]{{\left\| #1 \right\|}_{\Delta}}

\newcommand{\poly}{\mathcal{P}}

\newcommand{\mcl}{\mathcal{C}}

\newcommand{\autM}{{\Aut (M,0)}}
\newcommand{\autMp}{{\Aut (M,p)}}
\newcommand{\holmaps}{\mathcal{H}(\CN,0)}
\newcommand{\holmapsn}{\mathcal{H}(\Cn,0)}

\newcommand{\njetsp}[2]{J_{#1}^{#2} }
\newcommand{\jetm}[2]{ j_{#1}^{#2} }
\newcommand{\glnc}{\mathsf{GL_n}(\C)}
\newcommand{\glc}{\mathsf{GL_{(m+1)n}}(\C)}

\newlength{\extendaxesby}

\DeclareMathOperator{\imag}{Im} \DeclareMathOperator{\real}{Re}

\DeclareMathOperator{\proj}{pr}

\DeclareMathOperator{\Aut}{Aut}

\newtheorem{thm}{Theorem}
\newtheorem{lem}[thm]{Lemma}
\newtheorem{prop}[thm]{Proposition}
\newtheorem{cor}[thm]{Corollary}

\theoremstyle{definition}
\newtheorem{defin}[thm]{Definition}
\newtheorem{exa}[thm]{Example}
\newtheorem{rem}[thm]{Remark}

\newcommand{\biholmaps}{\mathcal{B}^n}
\newcommand{\omeganorm}[1]{{\left\| #1 \right\|}_{\Omega}}

\newcommand{\M}{\mathcal{M}}
\newcommand{\D}{\mathcal{D}}

\newcommand{\PP}{\mathcal{P}}

\begin{document}
\title[Parametrization of local CR Automorphisms ]{Parametrization of local CR Automorphisms
by finite Jets and Applications}
\author{Bernhard Lamel}
\address{Universit\"at Wien\\ Fakult\"at f\"ur Mathematik\\
Nordbergstrasse 15\\ A-1090 Wien\\ Austria}
\email{lamelb@member.ams.org}%
\thanks{The first author was supported by the FWF, Projekt P17111.}
\author{Nordine Mir}
\address{Universit\'e de Rouen\\ Laboratoire de Math\'ematiques 
Rapha\"el Salem, UMR 6085 CNRS\\
Avenue de l'Universit\'e, B.P. 12\\ 76801 Saint Etienne du Rouvray\\ France}
\email{Nordine.Mir@univ-rouen.fr}
\begin{abstract}
For any real-analytic hypersurface $M\subset \CN$, which does not contain any
complex-analytic subvariety of positive dimension, we show that
for every point $p\in M$ the local real-analytic CR automorphisms
of $M$ fixing $p$ can be parametrized real-analytically by their
$\ell_p$ jets at $p$. As a direct application, we derive a Lie group
structure for the topological group $\autMp$. Furthermore, we also
show that the order $\ell_p$ of the jet space in which the
group $\autMp$ embeds can be chosen to depend
upper-semicontinuously on $p$. As a first consequence, it follows that
that given any compact real-analytic hypersurface $M$ in $\CN$,
there exists an integer $k$ depending only on $M$ such for every
point $p\in M$ germs at $p$ of CR diffeomorphisms mapping $M$ into
another real-analytic hypersurface in $\CN$ are uniquely
determined by their $k$-jet at that point. Another consequence is
the following boundary version of H.\ Cartan's uniqueness theorem:
given any bounded domain $\Omega$ with smooth real-analytic
boundary, there exists an integer $k$ depending only on $\partial
\Omega$ such that if $H\colon \Omega\to \Omega$ is a proper
holomorphic mapping extending smoothly up to $\partial \Omega$
near some point $p\in
\partial \Omega$ with the same $k$-jet at $p$ with that of the
identity mapping, then necessarily $H={\rm Id}$.

Our parametrization theorem also holds for the stability group of any essentially finite minimal
real-analytic CR manifold of arbitrary codimension. One of the new main tools developed in the paper, which may be  of independent interest,  is a parametrization theorem for invertible solutions of a certain kind of singular analytic equations, which roughly speaking consists of inverting certain families of parametrized maps with singularities. 
\end{abstract}

\numberwithin{equation}{section}
\numberwithin{thm}{section}

\maketitle

\tableofcontents

\section{Introduction}\label{s:intro}
The goal of this work is to study  the structure of the group of
all local real-analytic CR automorphisms of a germ of a real-analytic 
hypersurface $(M,p)$ in $\CN$,  or equivalently, of all
local biholomorphisms $H\colon (\CN,p)\to (\CN,p)$ mapping the
germ $(M,p)$ into itself;
we also deal with the more general situation of a real-analytic CR
manifold of arbitrary codimension (see \S \ref{s:highcodim}). For
every fixed point $p\in M$, we refer to this group of local
biholomorphisms as the {\em stability group} of $(M,p)$ and denote
it by $\autMp$. Endowed with the topology of uniform convergence
on compact neighbourhoods of $p$, it  becomes a topological group.
This natural (inductive limit)
topology can be described in the following way: 
a sequence of germs of biholomorphisms $(H_j)\subset  \autMp$ converges to another such germ 
$H\in \autMp$ if there exists a  neighbourhood $U$ of $p$ in $\CN$ to which all mappings $H_j$ extend holomorphically
and such that $H_j \to H$ uniformly on $U$. 

This group arises in a variety of different circumstances; one of
the most important ones where it naturally appears is the {\em
biholomorphic equivalence problem} which consists in deciding whether
two given germs of real-analytic hypersurfaces are
biholomorphically equivalent. The biholomorphic equivalence
problem has been studied since the beginning of the last century
when Poincar\'e \cite{Po1} observed, probably for the first time,
that two real hypersurfaces in $\CN$ for $N\geq 2$ are not necessarily 
 locally equivalent via a holomorphic transformation (though they are via  a ${\mathcal C}^\infty$ transformation). Today, 
the biholomorphic equivalence
problem for strongly pseudoconvex (or Levi-nondegenerate)
hypersurfaces is pretty well understood, based on the 
works of E.\ Cartan \cite{Ca1, Ca2} in dimension two,  and the
results of Tanaka \cite{Ta1} and Chern-Moser \cite{CM} in arbitrary
dimension.

%One of the approaches to the equivalence problem is the
%construction of a normal form for a given family of hypersurfaces (this is e.g. carried
%out in \cite{CM} for Levi-nondegenerate hypersurfaces). 

%Let us sketch how this 
%approach works: One finds a way to attach to each
%hypersurface in the given family a certain special choice of coordinates
%and  a special form of a defining function which is then called the
%normal form. This choice is unique up to the action of a certain group, namely the 
%automorphism group of a certain ,,flat model hypersurface'' (in the 
%Levi-nondegenerate case, the models are the hyperquadrics). The coefficients of 
%the power series expansion of the normal form then depend on 
%the hypersurface being normalized and
%the special choice of normalization and solve the 
%biholomorphic equivalence problem in the following sense: Two hypersurfaces 
%are equivalent if and only if for they can be normalized to the same 
%normal form. 

%In principle, such a normal form solves the biholomorphic equivalence problem (there are
%actually infinitely many biholomorphic invariants
%attached to a real-analytic hypersurface), provided that the automorphism group
%of the model hypersurface is finite-dimensional; if this is the case, 
%for a generic hypersurface finitely many of the coefficients of the 
%normal form will actually determine a special choice of normalization. 

%In this
%paper, one of our central results is that the stability group of
%any piece of a real-analytic boundary of a bounded domain in $\CN$
%is always a Lie group (Theorem \ref{t:main3hsf}).

Thanks to the works initiated in the 70's, a number of remarkable properties
of the stability group of strongly pseudoconvex real-analytic
hypersurfaces have been discovered (see e.g.\ the surveys by
Vituskhin \cite{Vi}, Huang \cite{H4} and Baouendi, Ebenfelt and
Rothschild \cite{BERbull} for complete discussions on this matter). These properties depend on whether the
germ $(M,p)$ can be biholomorphically  mapped onto a piece of the
unit sphere; if there exists such a biholomorphism, one says that
$(M,p)$ is {\em spherical}. For instance, if $(M,p)$ is not
spherical the group $\autMp$ is compact (see Vitushkin~\cite{Vi1}); every element of such a
group is uniquely determined by its derivative at $p$ and also extends
holomorphically to a fixed neighbourhood (independent of the
automorphism) of $p$ in $\CN$. On the other hand, if $(M,p)$ is
spherical the group $\autMp$ is not compact and every local CR
automorphism of $(M,p)$ is uniquely determined by its $2$-jet (but
not $1$-jet) at $p$. Beyond their interest in their own right,
these local results also have direct applications to global
biholomorphic mappings of strongly pseudoconvex bounded domains
with smooth real-analytic boundaries by the classical reflection
principle \cite{Fe1, Le, P1}.

There are two main approaches to derive results about the stability
group of a strongly pseudoconvex (or Levi-nondegenerate)
real-analytic hypersurface in $\CN$. One consists in using the
powerful Chern-Moser theory (see for example the results of
Burns-Shnider \cite{BS} and the survey by Vitushkin \cite{Vi}).
The other one, initiated by Webster \cite{W4}, uses the invariant
family of Segre varieties attached to any real-analytic
hypersurface in a complex manifold.

While it seems difficult to extend the first approach to
understand the structure of the stability group of Levi-degenerate
hypersurfaces (see e.g.\ Ebenfelt \cite{E5}), the second one can
be carried over to that setting, as first observed by Baouendi,
Ebenfelt and Rothschild \cite{BER3}. In fact, it was shown in
\cite{BER3} that many interesting properties of the stability
group of a germ of a real-analytic hypersurface $(M,p)$ may be
obtained after showing that its local CR automorphisms are
parametrized (in a suitable sense) by their jets at $p$ of some
finite order. Such a program was successfully carried out in
\cite{BER3, Z2, BER4} for some classes of real-analytic
Levi-degenerate hypersurfaces and even CR manifolds of arbitrary
codimension. However, the degeneracies allowed for the
hypersurface $M$ in these results are quite restrictive and in
particular fail to hold at every point of real-analytic
hypersurfaces of $\CN$ containing no complex-analytic subvariety
of positive dimension. Such a class of real hypersurfaces is 
of particular interest since, by a well-known result of
Diederich-Forn\ae ss \cite{DF2}, boundaries of bounded domains
with smooth real-analytic boundary are of this type.

In dimension $N=2$, the stability group of such real hypersurfaces
has recently been studied by Ebenfelt, Zaitsev and the first
author in \cite{ELZ1} where it was shown that their local CR
automorphisms are analytically parametrized by their $2$-jets. In
this paper, among other results, we show that such a parametrization property by $k$-jets 
holds at {\em every} point of arbitrary real-analytic
hypersurfaces containing no complex analytic subvariety of positive
dimension and for some $k$ depending on the point (see Theorem \ref{t:main1hsf}).

 In what follows, for every integer $k$ and every pair of points
$p,p'\in \CN$, we denote by $G^{k}_{p,p'}(\CN)$ the group of all
$k$-jets at $p$ of local biholomorphisms $H\colon (\CN,p)\to
(\CN,p')$. Given coordinates $Z=(Z_1,\ldots,Z_N)$ in $\CN$ we will
usually identify the $k$-jet $j_p^kH$ of a local biholomorphism
$H$ at $p$ with
$$j_p^kH=\left(\frac{\partial^{|\alpha|}H}{\partial Z^{\alpha}}(p)\right)_{|\alpha|\leq k}.$$
If, in addition, $M,M'$ are real-analytic hypersurfaces of $\C^N$ or more generally real-analytic generic submanifolds of the same dimension with $p\in M$ and $p'\in M'$, we denote by ${\mathcal F}(M,p;M',p')$ the set of all germs at $p$ of local biholomorphisms $H\colon (\CN,p)\to
(\CN,p')$ sending $M$ into $M'$.

The purpose of this paper is to provide general conditions on a real-analytic generic submanifold $M\subset \C^N$
so that for every $p\in M$ and every $M'$ and $p'$ as above, the set of mappings ${\mathcal F}(M,p;M',p')$ is analytically parametrized by a finite jet at $p$. Our results are best illustrated in the hypersurface case by the following theorem which was the main motivation of the present work.

\begin{thm}
  \label{t:main1hsf} Let  $M$ be a real-analytic hypersurface in  $\CN$ containing no complex analytic subvariety of
  positive dimension. Then for every real-analytic hypersurface $M'\subset \CN$, every $p'\in M'$ and every $p\in M$, there exist an integer $\ell_p$ $($depending only on $M$ and $p)$,  an open subset $\Omega \subset \CN
  \times G_{p,p'}^{\ell_p}(\CN)$ and a real-analytic map
  $\Psi (Z,\Lambda) \colon \Omega \to \CN$, 
  holomorphic in the first factor, such that the following holds: 
  for any $H\in {\mathcal F}(M,p;M',p')$, 
  the point
  $(p,j_p^{\ell_p}H)$ belongs to $\Omega$ and the following identity holds:
  \[H(Z) = \Psi(Z,\jetm{p}{\ell_p} H)\text{ for all } Z\in \CN \text{ near } p.\]
  Furthermore, $\ell_p$ can be chosen so that
  it depends upper-semicontinuously on $p$.
\end{thm}

Our first application of Theorem \ref{t:main1hsf} concerns the
group $\autMp$. By using standard arguments from
\cite{BER4}, we obtain a Lie group structure on $\autMp$
compatible with its topology.
%\begin{cor}\label{t:main4hsf}
%   Let  $M$ be a real-analytic hypersurface in $\CN$ containing no complex analytic subvariety of
%  positive dimension through some point $p\in M$. If $\autMp$ is compact then there exists a fixed neighbourhood of
%  $p$ in $\CN$ to which all elements of $\autMp$ extend holomorphically.
%\end{cor}
%For strongly pseudoconvex real-analytic hypersurfaces, Corollary
%\ref{t:main4hsf} follows from Vitushkin's extension theorem for
%non-spherical hypersurfaces \cite{Vi}. Note that the well-known
%example of the unit sphere shows that Theorem \ref{t:main4hsf}
%fails to hold in the non-compact case.

\begin{thm}
  \label{t:main3hsf}  Let  $M$ be a real-analytic hypersurface in $\CN$ containing no complex analytic subvariety of
  positive dimension. Then for every $p\in M$, there exists an integer $\ell_p$, depending upper-semicontinuously on $p$,
  such that the jet mapping
  \begin{equation}\label{e:injectivity}
   j_p^{\ell_p} \colon \autMp \to G_{p,p}^{\ell_p}(\CN)
   \end{equation}
  is a continuous group homomorphism that is a homeomorphism onto a real-algebraic Lie subgroup of $G_{p,p}^{\ell_p}(\CN)$.
\end{thm}
We should mention that the real-algebraicity of the Lie group
$\autMp$ in Theorem \ref{t:main3hsf} does not follow directly from
Theorem \ref{t:main1hsf} but from the more precise version given
by Theorem \ref{t:main1} in \S \ref{s:highcodim} below. Another
noteworthy consequence of Theorem \ref{t:main3hsf} is the finite
jet determination of local CR automorphisms of $(M,p)$ by their
$\ell_p$-jet at $p$. As stated, this finite jet determination
property would appear as already known since the existence of 
some integer $k$ for which the finite determination holds follows from the
work of Baouendi, Ebenfelt and Rothschild \cite{BER5}. However the
proof given in \cite{BER5} does not give any information on the
dependence of the jet order required with respect to a given point
of the hypersurface. The existence of an integer $\ell_p$
depending {\em upper-semicontinuously} on $p$ for which the
mapping given by \eqref{e:injectivity} is  injective is
completely new and crucial in order to get the following.

\begin{thm}
  \label{t:main2hsf} Let $M$ be a compact real-analytic hypersurface in $\CN$.
  Then there is an integer $k$, depending only on $M$,
  such that for every $p\in M$ and for every real-analytic hypersurface
  $M^\prime\subset \CN$, smooth local CR diffeomorphisms
  mapping a neighbourhood of $p$ in $M$ into $M^\prime$ are uniquely determined by their $k$-jets at $p$.
\end{thm}

Theorem \ref{t:main2hsf} follows from the upper-semicontinuity of
the map $p\mapsto \ell_p$ in Theorem \ref{t:main3hsf}, the fact
that compact real-analytic hypersurfaces do not contain any
holomorphic curves \cite{DF2} and from the regularity result for
smooth CR diffeomorphisms proved in \cite{BJT}.  For proper
holomorphic mappings of bounded domains with smooth real-analytic
boundaries, we also obtain the following rigidity result as a
consequence of Theorem \ref{t:main2hsf}.

\begin{cor}
  \label{c:cormain3}  Let $\Omega\subset\CN$ be a bounded domain with smooth real-analytic boundary. Then there
  exists an integer $k$, depending only on the boundary $\partial \Omega$, such that
  if $H\colon \Omega \to \Omega$ is a proper holomorphic mapping extending
  smoothly up to $\partial \Omega$ near some point $p\in \partial \Omega$
 which satisfies $H(z) = z+o ( \left|z - p \right|^k)$ when $z\in \Omega$ and $z\to p$, then $H$ is the identity mapping.
\end{cor}

Corollary \ref{c:cormain3} can be viewed as a boundary version of
H.\ Cartan's uniqueness theorem \cite{HCa} and seems to be new
even in the case of weakly pseudoconvex domains (see e.g.\
\cite{BK, H5, H6}). Note that the smooth extension up to the
boundary assumption is known to hold automatically near every
point $p\in
\partial \Omega$ in many cases e.g.\ if $\Omega$ is weakly
pseudoconvex (see \cite{BERbook} for more on that matter). For
bounded strongly pseudoconvex domains with ${\mathcal C}^\infty$
boundary, more precise results are known from the work of
Burns-Krantz \cite{BK} (see also Huang \cite{H5, H6}).

\section{Results for CR manifolds of higher codimension and main tools}\label{s:highcodim}
We state here our results for real-analytic generic submanifolds
of higher codimension from which all theorems mentioned in the
introduction will be derived. 
Recall that a CR submanifold $M$ of $\CN$ is
called {\em generic} if $T_pM+J(T_pM)=T_p\CN$ where $J$ is the
complex structure map of $\CN$ and $T_pM$ (resp.\ $T_p\CN$)
denotes the tangent space of $M$ (resp.\ of $\CN$) at $p$ (see
e.g.\ \cite{boggess, BERbook}).

In order to state our results, we need to impose two nondegeneracy
conditions on a given real-analytic generic submanifold $M\subset
\CN$. We shall assume that $M$ is {\em essentially finite} in the
sense of \cite{BJT}  and {\em minimal} in the sense of \cite{TU2}
at each of its points (see \S \ref{ss:nondegeneracy} for precise
definitions). Such conditions are very natural and hold in
general situations for compact as well as non-compact
real-analytic CR submanifolds in complex space. For instance, if
$M$ does not contain any complex analytic subvariety of positive
dimension then $M$ is necessarily essentially finite and also
automatically minimal when $M$ is furthermore of (real)
codimension one (see e.g.\ \cite{BERbook}).

Our first result in this section is the following.

\begin{thm}
  \label{t:main1} Let  $M$ be a real-analytic generic submanifold of $\CN$ that is essentially finite and minimal
  at each of its points. Then for every real-analytic generic submanifold $M'\subset \CN$, every point $p'\in M'$ and every $p\in M$, there exists an integer $\ell_p$ $($depending only on $M$ and on $p)$, an open subset
  $\Omega \subset \CN \times G_{p,p'}^{\ell_p}(\CN)$ and a real-analytic map
  $\Psi (Z,\Lambda) \colon \Omega\to \CN$ holomorphic in the first factor, such that the following hold :
  \begin{enumerate}
  \item[{\rm (i)}] for any $H\in {\mathcal F}(M,p;M',p')$ the point
  $(p,j_p^{\ell_p}H)$ belongs to $\Omega$ and the following identity holds:
  \[H(Z) = \Psi(Z,\jetm{p}{\ell_p} H)\text{ for all } Z\in \CN \text{ near } p;\]
  \item[{\rm (ii)}] the map $M\ni p\mapsto \ell_p\in \N$ is
  upper-semicontinuous;
  \item[{\rm (iii)}] the map $\Psi$ has the following formal Taylor
  expansion
  \begin{equation}
  \Psi (Z,\Lambda)=\sum_{\alpha \in \N^N}\frac{P_\alpha (\Lambda,\overline{\Lambda})}
  {(D(\Lambda^1))^{s_\alpha}(\overline{D({\Lambda^1})})^{r_\alpha}}\, (Z-p)^{\alpha},
  \end{equation}
  where for every $\alpha \in \N^N$, $s_\alpha$ and $r_\alpha$ are  nonnegative
  integers, $P_\alpha$ and $D$ are polynomials in their arguments, 
  $\Lambda^1$ denotes the linear part of the jet $\Lambda$ and $D(j_0^1H)\not =0$ for every 
  $H\in {\mathcal F}(M,p;M',p')$.
  \end{enumerate}
\end{thm}

To our knowledge, Theorem \ref{t:main1} is already new even in the
case of real hypersurfaces in $\CN$. (For $N=2$, Theorem
\ref{t:main1} (i) with $\ell_p=2$ is one of the main results of
\cite{ELZ1}.) On the other hand, for CR manifolds of codimension
bigger than one, only special cases of Theorem \ref{t:main1}
(including the Levi-nondegenerate case) were known as consequences
of the works in \cite{Z2, BER4}.

%, we first have the following direct
%consequence of Theorem \ref{t:main1}.

%\begin{cor}\label{t:main4}
%   Let  $M$ be a real-analytic generic submanifold of $\CN$.
%   Assume that $M$ is minimal and essentially finite at some point $p\in M$.
%   If $\autMp$ is compact then there exists a fixed neighbourhood of
%  $p$ in $\CN$ to which all elements of $\autMp$ extend holomorphically.
%\end{cor}

As in \S \ref{s:intro}, we also derive from Theorem \ref{t:main1} a Lie group structure on
the stability group of any essentially finite minimal
real-analytic generic submanifold of $\CN$.

\begin{thm}
  \label{t:main3} Let  $M$ be a real-analytic generic submanifold of $\CN$ that is essentially finite and minimal
  at each of its points. Then for every $p\in M$ there exists an integer $\ell_p$, depending upper-semicontinuously on $p$,
  such that the jet mapping
  \begin{equation}\label{e:hurryup}
   j_p^{\ell_p} \colon \autMp \to G_{p,p}^{\ell_p}(\CN)
   \end{equation}
  is a continuous group homomorphism that is a homeomorphism onto a real-algebraic Lie subgroup of $G_{p,p}^{\ell_p}(\CN)$.
\end{thm}
Even in the case of real hypersurfaces $M\subset \CN$ which do not
contain any complex analytic subvariety of positive dimension with
$N\geq 3$, the fact that $\autMp$ is a Lie group was an open
problem. Note also that the fact that the image of $\autMp$ under
the jet mapping $j_p^{\ell_p}$ is a {\em real-algebraic} subgroup
of $G_{p,p}^{\ell_p}(\CN)$ follows from the rational dependence of the
parametrization $\Psi$ with respect to the jet $\Lambda$ in
Theorem \ref{t:main1} (iii) (see \S \ref{s:conclusion}), which is
new even in the case $N=2$. For submanifolds of higher codimension
only special cases of Theorem \ref{t:main3} were known from the
works \cite{Z2, BER4}. We should point out that even in these
special cases the jet order $\ell_p$ given by our Theorem
\ref{t:main3} can be computed explicitly (from
integer-valued biholomorphic invariants of $(M,p)$)
and is always smaller
than the jet order required in the papers \cite{Z2, BER4}. For
details on that matter we refer the reader to \S
\ref{ss:maincrap}--\ref{ss:propertiesmcl}. We would also like to
recall that the statement concerning the injectivity of the
mapping given by \eqref{e:hurryup} is also new by itself since a
choice of an upper-semicontinuous bound on the jet order required
to get injectivity does not follow from the work \cite{BER5}. The
following finite jet determination result for compact
real-analytic CR manifolds can be viewed as a higher codimensional
version of Theorem \ref{t:main2hsf} and is a direct consequence of
the possibility of an upper-semicontinuous choice of $\ell_p$.

\begin{thm}
  \label{t:main2} Let $M$ be a compact real-analytic CR submanifold of $\CN$ minimal at each of its points.
  Then there is an integer $k$, depending only on $M$, such that
  for every $p\in M$ and for every real-analytic CR submanifold
  $M^\prime\subset \CN$ with the same CR dimension as that of $M$, smooth local CR diffeomorphisms
  mapping a neighbourhood of $p$ in $M$ into $M^\prime$ are uniquely determined by their $k$-jets at $p$.
\end{thm}

Let us explain now in some greater detail
the main new difficulties faced
in this paper compared to previous related work. The first main
ingredient of the proof of Theorem \ref{t:main1} relies on  the
full machinery of the Segre sets technique introduced by Baouendi,
Ebenfelt and Rothschild \cite{BERbook, BER4}. In order to use such
a tool one needs to show that every local biholomorphism fixing a
germ of a submanifold $(M,p)$ (satisfying the conditions of the
theorem) as well as all its jets, when restricted to any Segre set
attached to $p$, is analytically parametrized (in a suitable
manner) by its jet at $p$ of some finite precise order. 
However, in order to realize such a program, one immediately runs into the problem of inverting certain families of parametrized maps with singularities in multidimensional complex space, a phenomenon that does not appear in all previous related work. Here the existence of such singularities comes directly from the high Levi-degeneracy allowed
for the submanifold $M$ in Theorem \ref{t:main1}. To overcome the difficulty, we will be led to 
%solve an independent
%problem which consists of  that solutions of
 consider a precise type of linear and non-linear singular analytic systems
of equations near the origin in $\C^n$.  The most difficult type of singular systems we end up with
is of the form $A(w(z)) = b(z)$ where $A, b\colon(\Cn,0) \to (\Cn,0)$ are germs of holomorphic maps with $A=A(w)$ of generic rank $n$  and one wants to solve such a system for a biholomorphism
$w\colon(\Cn,0)\to(\Cn,0)$; the solution(s) should, for our
application, preserve the dependence of $b$ on arbitrary deformations, that is,  any continuous (resp.\ smooth, resp.\ holomorphic) perturbation of the term $b$ should lead to a solution $w=w(z)$ that behaves in the same manner. 
Unfortunately, there is, to our knowledge, no general theory to deal with this 
situation, so we have to develop our own machinery to parametrize
the invertible solutions of such singular systems.

Since our solution to the above problem may be of independent interest, we introduce some notation which allows
us to present it
in a particularly simple manner. It is useful here to work in the 
context of holomorphic maps of infinite dimensional complex vector spaces.
Let $\holmapsn$ be the space of germs of holomorphic maps
$H:(\Cn,0)\to \Cn$, $n\geq 1$, endowed with the topology of
uniform convergence on compact neighbourhoods of the origin; here again, a sequence of germs of holomorphic maps $(H_j)$ converges to another such germ 
$H$ if there exists a  neighbourhood $W$ of the origin in $\CN$ to which all $H_j$ extend holomorphically
and such that $H_j \to H$ uniformly on $W$. With
such a topology, the space $\holmapsn$ carries the structure of an
infinite dimensional topological complex vector space being an
inductive limit of Fr\'echet spaces (see \S \ref{ss:end?} for more details). Let $\holmapsn^0$ be the
subspace of elements of $\holmapsn$ preserving the origin and
$\biholmaps$ be the open subset of $\holmapsn^0$ consisting of
the germs of biholomorphic maps (preserving the origin). 

Given
any $A\in \holmapsn$, composition with $A$ naturally yields a
(holomorphic) map $A_{*}: \biholmaps \to \holmapsn $ given by $A_{*}(u):=A\circ u$ for $u\in \biholmaps$. Our above mentioned problem calls 
for the construction of 
a holomorphic ``left inverse'' to $A_{*}$, that
is, a holomorphic map $\Psi \colon \holmapsn\to \biholmaps$
satisfying $\Psi( A_{*} u) = u$. 
Here holomorphy of a map $f\colon X\supset \Omega\to Y$ between locally
convex topological complex vector spaces $X,Y$ with $\Omega$ open in $X$ is understood in the
sense of \cite[Definition 3.16]{Dineenbook}: $f$ is holomorphic if
it is continuous and its composition with any continuous linear
functional on $Y$ is holomorphic along finite dimensional affine
subspaces of $X$ intersected with $\Omega$. 

Since the mapping $A_{*}$ is
obviously not injective in general, there is of course no hope to
find such a left inverse. However, it was recently shown by the
first author in \cite{L4} that the holomorphic map $A_{*}\times
j_0^1\colon \biholmaps \ni u \mapsto (A_{*}(u),u'(0))\in \holmapsn
\times \glnc$ is injective {\em when $A$ is of generic rank} $n$. 
(As usual, $\glnc$ denotes the group of invertible
$n\times n$ matrices and $\jetm{0}{1} u = u'(0)$ 
the differential or $1$-jet of $u$ at the origin.) It thus
makes sense to try to find a left inverse to that latter map. 
Our solution provides such an inverse that we also call parametrization.

\begin{thm}
  \label{t:parcnfull} Let $A:(\Cn,0)\to \Cn$ be a germ of a holomorphic map of
  generic rank $n$. Then there exists a holomorphic map $\Psi \colon \holmapsn\times\glnc\to
  \biholmaps$ such that $\Psi$ satisfies $\Psi(A_{*}u, u'(0)) =
  u$ for all $u\in \biholmaps$; furthermore, $\Psi$ 
  can be chosen so that it satisfies $j_0^1 \left(\Psi(b,\lambda)\right) = 
  \lambda$ for every $b\in \holmapsn$ and $\lambda \in \glnc$.
\end{thm}

This result provides us the necessary information
about the structure of the invertible solutions to the equation
$A(w) = b(z)$ and allows us to proceed with the construction 
of the parametrization of the automorphism group of CR manifolds. 
The fact that we are even able to parametrize the
solutions of the above type of equations precisely by their
$1$-jets is one of the main reasons why we succeed in
getting an upper-semicontinuous dependence of the integer $\ell_p$ on
$p$ in Theorem \ref{t:main1}. Another  result in the 
spirit of Theorem~\ref{t:parcnfull}, for a certain kind of linear singular 
system that we also
need is Proposition \ref{t:lineqn2} which
provides the
parametrization of this type of equations. 
In both cases, we need a very careful construction 
of such a parametrization (which in the 
case of Theorem~\ref{t:parcnfull} is 
provided by the more precise Theorem~\ref{t:parcn1}) 
in order to prove statement (iii) of Theorem
\ref{t:main1}.

The paper is organized as follows. We start by providing the
necessary background on parametrization of non-linear singular systems
in \S\ref{s:nonlinearcase}; the proof of the precise 
results in that section 
and the proof of Theorem~\ref{t:parcnfull} are
given in 
\S\ref{ss:justwrite}--\S\ref{s:holomorphicity} (following the 
development of some necessary material on homogeneous
polynomials in
\S\ref{ss:tiredofwriting}); we also include an 
additional application of Theorem~\ref{t:parcnfull} in
\S\ref{s:onemore}. In \S \ref{s:linearcase}, we treat the necessary case of a different (and simpler) 
type of singular systems, namely
 linear ones, following the approach in the non-linear case.

After that, we proceed with the construction of the 
parametrization of the automorphism group of CR manifolds. 
In \S \ref{s:class} we introduce a class of real-analytic generic
submanifolds which is more general than the class of essentially
finite ones. After discussing a few properties of this class of
submanifolds and giving several examples, we state the most
general parametrization theorem of this paper (Theorem
\ref{t:further1}) which holds for the stability group of any
submanifold of the above mentioned class. The proof of Theorem
\ref{t:further1} is furnished in  \S \ref{s:jets}--\S
\ref{s:jetpar} and is completed in \S \ref{s:conclusion}, where we
also prove all the remaining theorems stated in \S \ref{s:intro}--\S \ref{s:highcodim}.

\begin{rem}[Jet space notation] Jet spaces of various kinds will be
  encountered frequently in this paper. For a complete discussion,
  we refer the interested reader to e.g. \cite{GG}
  or \cite{BCG}. Here we will use the 
  following notation: for positive integers $n,k,l$, 
  we denote by
$J_0^l(\C^n,\C^k)$ the jet space at the origin of order $l$ of
holomorphic mappings from $\C^n$ to $\C^k$.
For a germ of a holomorphic map $h\colon (\Cn,0)\to \C^k$, we
denote by $\jetm{0}{l}h$ the $l$-jet of $h$ at $0$. If $k=n$, we
simply write $J_0^l(\Cn)$ for $J_0^l(\Cn,\Cn)$. The 
jet space of holomorphic mappings $(\Cn,0) \to 
(\C^k,0)$ will be denoted by 
$J_{0,0}^l (\Cn,\C^k)$, or simply $J_{0,0}^l (\Cn)$
if $n=k$. We have already introduced above $G_{0,0}^l (\Cn)$ as
our notation for the jet group of order $l$
of biholomorphic mappings
$(\Cn,0)\to (\Cn,0)$. 
\end{rem}

\medskip

{\bf Acknowledgements:} The authors are grateful to the referee for a number of remarks improving the readability and the exposition of the paper.

%\part{Parametrization of solutions of singular analytic equations}

\section{Parametrizing invertible solutions of non-linear singular 
analytic equations}\label{s:nonlinearcase}

We will derive Theorem
\ref{t:parcnfull} in \S \ref{s:holomorphicity} from the 
more precise result given in Theorem
\ref{t:parcn1} below. For every
$\C^k$-valued (homogeneous polynomial) map $f\colon \C^n\to \C^k$,
and for any bounded set $F \subset \Cn$, we define
$$\vnorm{f}_{F}:= \max_{1\leq j\leq k} \sup_{z\in F}\, |f^j(z)|,\quad
f=(f^1,\ldots,f^k).$$

\begin{thm}
  \label{t:parcn1} Let $A:(\Cn,0)\to \Cn$ be a germ of a holomorphic map
 of generic rank $n$. Then there
  exists an integer $\ell_0$ and for every integer $\ell>0$ polynomial mappings
  $p_\ell:\Cn\times\glnc\times\njetsp{0}{\ell_0 + \ell}(\C^n)\to\Cn$,
  homogeneous of degree $\ell$ in their first variable, and an integer $\kappa (\ell)$
  such that for every germ of a holomorphic map
  $b \colon (\Cn,0)\to \Cn$,  if $u:(\Cn,0)\to(\Cn,0)$ is a germ of
  a biholomorphism satisfying $A(u(z)) = b(z)$, then necessarily
  \begin{equation}
    u_\ell (z) = \frac{1}{(\det u^\prime (0))^{\kappa (\ell)}} \, p_{\ell } ( z, u^\prime(0), \jetm{0}{\ell_0 + \ell} b),
    \label{e:parcn1}
  \end{equation}
  where $u(z) = \sum_{\ell \geq 1} u_\ell (z)$ is the decomposition of
  $u$ into homogeneous $\C^n$-valued polynomials terms.
  Furthermore, the polynomial maps $p_\ell$ and the integers $\kappa (\ell)$ can be chosen with the following two properties:
  \begin{enumerate}
  \item[{\rm (i)}] for every germ of a holomorphic map $b\colon (\Cn,0)\to
  \Cn$ and for every $\genmat \in \glnc$
  \begin{equation}\label{e:7}
  p_1(z,\genmat,j_0^{\ell_0+1}b)=\genmat z, \text{ and }
  \kappa(1)=0;
  \end{equation}
  \item[{\rm (ii)}] there exists a polydisc $\Delta\subset \C^n$ centered at the origin
   such that for every $C>0$ and every compact subset $L\subset\glnc$ there
  exists $K>0$ such that for every $\genmat\in L$ and every germ of a holomorphic map $b\colon (\Cn,0)\to \Cn$ written into homogeneous terms
    $b(z) = \sum_{k\geq 0} b_k (z)$ satisfying  $\deltanorm{b_k} \leq C^k$ for all $k\in \N$, we
    have for all $\ell \geq 1$
  \begin{equation}\deltanorm{ \frac{p_{\ell} (\cdot, \genmat, \jetm{0}{\ell_0 + \ell} b) }{(\det \genmat)^{\kappa (\ell)}}} \leq K^\ell.
    \label{e:parcnest}
  \end{equation}
  \end{enumerate}
\end{thm}

As explained in \S \ref{s:holomorphicity}, to every choice of polynomial mappings $(p_\ell)$ and integers $(\kappa (\ell))$ satisfying the conclusion of Theorem~\ref{t:parcn1}, the map 
\begin{equation}\label{e:formalstuff}
\Psi (b,\lambda):=\sum_{\ell \geq 1} \frac{p_{\ell} (\cdot,
\genmat, \jetm{0}{\ell_0 + \ell} b) }{(\det \genmat)^{\kappa
(\ell)}},\quad b\in
\holmapsn,\quad \lambda \in \glnc.
\end{equation}
will be a holomorphic parametrization satisfying the conclusions of Theorem \ref{t:parcnfull}. 

Let us add that it follows from the proof of Theorem
\ref{t:parcn1} that we can choose $\kappa(\ell) = \ell_0 \ell$.
However, the exact value of $\kappa$ is not important for our
purposes. If in Theorem \ref{t:parcn1} we just consider formal
maps $A$ instead of holomorphic ones the above theorem remains
true in the formal category without the estimates in (ii). In
fact, we have a formal counterpart of the above theorem as
follows. Consider  a formal power series mapping $A$ as in Theorem~\ref{t:parcn1} (not necessarily convergent) and 
the sets $\holmapsn_f$ and $\biholmaps_f$ of
$\Cn$-valued formal holomorphic maps and formal biholomorphic maps
preserving the origin respectively. For $A$ as above, let $(p_\ell)$ and $(\kappa (\ell))$ satisfy the conclusion of Theorem~\ref{t:parcn1} (i). Then for every $\lambda \in
\glnc$ and every $b\in\holmaps_f$, formula 
\eqref{e:formalstuff} defines a formal power series mapping that belongs to $\biholmaps_f$
and satisfies $\Psi(A_{*} u , u^\prime(0)) = u$ for
$u\in\biholmaps_f$. If we consider $b$'s which are elements of
$\mathcal{R}[ [ z ] ]$ for some algebra $\mathcal{R}$ over $\C$,
then we get that $\Psi(b,\genmat)\in \mathcal{R} \left( \genmat
\right) \left[ \left[ z \right] \right]$. This means that the
obtained parametrization $\Psi$ also preserves the dependence on
parameters in a formal way.

As mentioned before, the obtained parametrization in Theorem
\ref{t:parcnfull} or Theorem \ref{t:parcn1} allows us to consider
arbitrary perturbations on the right hand side of the singular
system. 
For the purposes of this paper, the dependence 
on holomorphic parameters given in Corollary~\ref{c:corparcn1} below is crucial: 
the precise result can be derived from Theorem \ref{t:parcn1} and Lemma
\ref{l:normsconv}, but we note that Corollary~\ref{c:corparcn1} (i) also follows directly from Theorem~\ref{t:parcnfull}.

\begin{cor}
  \label{c:corparcn1}   Let $A:(\Cn,0)\to \Cn$ be a germ of a holomorphic map
 of generic rank $n$, $X$ a complex manifold, and
  $b=b(z,\omega)$ be a $\Cn$-valued holomorphic map defined
  on an open neighbourhood of $\{0\}\times {X}\subset \Cn\times {X}$.
   Then there exists a holomorphic map $\Gamma=\Gamma(z,\genmat,\omega)\colon \Cn\times\glnc\times {X}\to \Cn$,
defined on an open neighbourhood $\Omega$ of $\{0\}\times
\glnc\times {X}$,
  such that if $u:(\Cn,0)\to(\Cn,0)$ is a germ of a biholomorphism satisfying
  $A(u(z)) = b(z,\omega_{0})$ for some $\omega_{0}\in {X}$, then necessarily
  $u(z) = \Gamma(z,u^\prime(0),\omega_{0})$. Furthermore, the map $\Gamma$
  can be chosen with the following properties :
  \begin{enumerate}
 \item[{\rm (i)}] for every $(\genmat,\omega)\in \glnc \times X$, $\Gamma(0,\genmat,\omega) = 0$ and $\Gamma_z (0,\genmat,\omega)=
 \genmat$; \label{e:Fprops}
  \item [{\rm (ii)}] if we write the Taylor expansions
  $$\Gamma(z,\genmat,\omega) = \sum_{\alpha\in \N^n} \Gamma_\alpha (\genmat,\omega) z^{\alpha},\quad b(z,\omega)=\sum_{\alpha \in \N^n}b_{\alpha}(\omega)z^\alpha,$$
then there exists an integer
  $\ell_0$ and for every $\alpha \in \N^n$ nonnegative integers $k_\alpha$ and polynomial mappings $p_\alpha \colon \glnc \times J_0^{\ell_0+|\alpha|}(\Cn)\to \C^n$
   such that
  \begin{equation}
    \Gamma_\alpha (\genmat,\omega) = \frac{p_\alpha (\genmat,(b_\beta(\omega))_{|\beta|\leq \ell_0+|\alpha|})}{(\det
    \genmat)^{k_\alpha}}.
    \label{e:polynomialpersist}
  \end{equation}
  \end{enumerate}
\end{cor}

%Another consequence of Theorem \ref{t:parcn1} is the $1$-jet determination of invertible
%solutions of the corresponding complex-analytic systems already mentioned above.

An outline of the proof of Theorem \ref{t:parcn1} is as follows:
First (in Lemma \ref{l:simplify1}), we show that we may
assume that the lowest order homogeneous terms of $A$ are
generically independent. Such a reduction then allows us to
linearize our initial system by considering homogeneous expansions
of the power series involved in the system. By using elementary
linear algebra on vector spaces of homogeneous polynomials, one
may then conclude the formal part of the proof. However, for the
last and more tedious part of the proof, namely the convergence
part (given by Theorem \ref{t:parcn1} (ii)), a more careful
analysis of the formal part is required in order to get the needed
estimates, which are obtained by using a majorant method.

Let us summarize the organization of the material in the following
sections again:
In \S\ref{ss:tiredofwriting}, we introduce the notation and
collect all the necessary facts about homogeneous polynomials
needed for the proof of Theorem \ref{t:parcn1}; a
particular normalization of the equations is discussed 
in \S\ref{ss:reduce} and the proof of Theorem~\ref{t:parcn1}
is completed
in \S\ref{ss:conclpfthm31}. In \S \ref{s:holomorphicity}, we show
how to derive Theorem~\ref{t:parcnfull} from Theorem
\ref{t:parcn1}. \S \ref{s:onemore} offers another application of
the above results which can be seen as the continuous version of
Corollary \ref{c:corparcn1}. In \S \ref{s:linearcase}, we
illustrate how the method of proof of Theorem \ref{t:parcn1} can
be used to give a precise result, which we also need, 
for the parametrization of solutions of other types of 
singular systems (which are simpler, because they are linear). 
This section actually serves as a guide to our
proof of Theorem~\ref{t:parcn1}, the steps being followed in
precisely the same manner, but the technicalities involved are
much easier to deal with. As a consequence, the reader may like to
start with that part before going through the details of the proof
of Theorem \ref{t:parcn1}.

\section{Proof of Theorem \ref{t:parcn1}}\label{ss:justwrite}
\subsection{Preliminary results about homogeneous
polynomials}\label{ss:tiredofwriting} Let $\poly_{n,d}$ be the
vector space of homogeneous polynomials of degree $d\geq 0$ in $z
= (z_1,\dots,z_n)\in \Cn$, $n\geq 1$. If $n$ is fixed, we are
dropping the $n$ and simply write $\poly_d$. Given any positive
integer $m$, note that the space $\poly_d^m$ can be identified
with the space of all homogeneous $\C^m$-valued polynomials maps
of degree $d$. Given any bounded open set $\Omega \subset \Cn$
containing the origin, we are going to consider the following
norms on the vector space $\poly_d$:
\begin{align*}
%\label{e:normsdefined}
 \onenorm{f} &:= \sum_{|\alpha|=k} \left| f_\alpha \right| , \text{ where } f (z)= \sum_{\left| \alpha \right| = d} f_\alpha z^\alpha,\\ 
  \omeganorm{f} &:= {\rm sup}_{z\in\Omega} \left| f(z) \right|.
\end{align*}
Given $(k_1,\ldots,k_m)\in \N^m$, we consider the following norms
on the cartesian product $\poly_{k_1}\times \ldots \times
\poly_{k_m}$:
\begin{align*}
%\begin{split}
 \onenorm{f} &:= \max_{j=1,\dots,m} \onenorm{f^j}, \quad f = (f^1,\dots,f^m), f^j \in \poly_{k_j},\\
 \omeganorm{f} &:= \max_{j=1,\dots,m} \omeganorm{f^j}, \quad f = (f^1,\dots,f^m), f^j \in \poly_{k_j}.
 % \label{e:normsdefined2}\end{split}
\end{align*}
Let us first note the following fact:
\begin{lem}
  \label{l:multmap} Let $M=M(z)$ be a $m\times m$ matrix with
  holomorphic coefficients near the origin in $\Cn$, $z=(z^1,\dots,z^n)$. Assume that  $\det M \not \equiv 0$ $($in any  neighbourhood of $0)$ and
  that for every $j=1,\ldots,m$, the $j$-th row of $M$ is formed by homogeneous polynomials of
degree $k_j\geq 0$. Then there exists a bounded open set  $\Omega
\subset \C^n$ containing the origin and a constant $C>0$, only
depending on $M$,  such that for every $d\in \N$ and for all
$k_1,\dots,k_m$ there exists a linear operator $B:
\poly_{k_1+d}\times \ldots \times \poly_{k_m+d}\to \poly_d^m$
which is inverse to matrix multiplication by $M$ of norm at most
$C$, that is $B$ satisfies
  \begin{equation}\label{e:constant}
    \vnorm{B s}_\Omega \leq C \vnorm{ s }_\Omega,\text{ for all } s\in  \poly_{k_1+d}\times \ldots \times \poly_{k_m+d},
    \text{ and } B(M \cdot p) = p \text{ for all }
  p\in\poly_d^m.
  \end{equation}
  After a linear change of coordinates $($normalizing $\det M$ in the $z^n$-direction$)$,
  we may actually assume that $\Omega$ is a polydisc.
\end{lem}
\begin{proof}
  Let $\widehat M(z)$ denote the transpose of the comatrix of $M(z)$ i.e.\ the matrix satisfying
  \[ M(z)\cdot \widehat M(z) = \widehat M(z)\cdot M(z) = \left( \det M(z) \right) I_m,\]
  where $I_m$ is the identity matrix. In particular, for any $d\in \N$ and for every $p\in\poly_d^m$ we have
  \[ (\det M(z)) p(z) = \widehat M(z)\cdot M(z)\cdot p(z).\]
  Since $\det M \not \equiv 0$, the Weierstrass division theorem (see e.g.\ \cite{Hor}) implies that
  we can find a linear change of coordinates (making the function $\det M$ $z^n$-regular of order $r\geq 1$) and, in the new coordinates, a polydisc $\Delta \subset \Cn$
  centered at the origin  and a constant $\tilde{C}$ (depending only on $M$) such that
  for any holomorphic function $f$ holomorphic in a neighborhood of $\overline{\Delta}$, 
$$f=q\, \det M+h,$$
where $h$ is a unique polynomial in $z^n$ of degree $<k$, and $q$ (the quotient) is a unique holomorphic function in $\Delta$ that satisfies $\deltanorm{q}\leq \tilde C \deltanorm{f}$.
  We define the map $B$ by associating to $s\in  \poly_{k_1+d}\times \ldots \times \poly_{k_m+d}$ the quotient given by
  Weierstrass division of $\widehat M(z)  \cdot s(z)$ by $\det M$. It is easily checked that this is a homogeneous polynomial map of
  degree $d$, and it satisfies
  $\deltanorm{ B s } \leq \tilde{C} \deltanorm{ \widehat M  s } \leq {C}
  \deltanorm{s}$, for some constant $C>0$. This proves \eqref{e:constant} on the open set $\Omega$ which
corresponds to the polydisc in the original coordinates.
\end{proof}

In what follows we will make extensive use of Faa di Bruno's
formula, which will be only needed for formal power series (so
that we actually only use the multinomial formula; we will still
refer to it as Faa di Bruno's formula). We give a statement in
this case (see e.g.\ \cite{KP}):
\begin{prop}
  \label{p:faadebruno} Assume that $f(t)= \sum_r a_r t^r$, $g(x) = \sum_q b_q x^q$, and
  $h(t)= \sum_s c_s t^s$
  are three complex-valued formal power series in one variable with $c_0=0$ and $f (t) = g (h (t))$.
  Then for every positive integer $r$, one has
  \begin{equation}
        a_r = \sum \frac{q!\, b_q}{k_1 ! \dots k_r!} c_1^{k_1} \dots c_r^{k_r} ,
    \label{e:faadebruno}
  \end{equation}
  where $q = \sum_{1\leq j\leq r} k_j$ and the sum goes over all nonnegative integers $k_1,\ldots,k_r$ such that
  \[ k_1 + 2 k_2 + \dots + r k_r = r.\]
\end{prop}

We now turn to some estimates for homogeneous parts of
compositions of certain types of formal power series mappings. We
use the following notation for decomposing a formal power series
mapping $u=(u^1,\dots,u^n):(\CN,0)\to(\Cn,0)$, $N\geq 1$, in
Taylor series or into polynomial homogeneous terms :
\begin{equation}\label{e:decomphomog}
 u(\zeta ) = \sum_{\alpha \in
\N^N} u_{\alpha} \zeta^\alpha = \sum_{j=1}^\infty u_{j} (\zeta), \quad
u_j (t \zeta) = t^j u_j(\zeta)\ {\rm for}\ {\rm all}\ (t,\zeta)\in
\C \times \C^N.
\end{equation}
In other words, for every $j\in \N^*$, $u_j \in (\poly_{N,j})^n$. From time to time, we also use the notation
$u|_j$ to denote $u_j$ (especially when the expression of the power series mapping $u$ is long). 
For each component $u^k$ of $u$, we use an analogous notation for
the decomposition of $u^k$ into homogeneous terms.  Note that for
every polydisc $\Delta \subset \Cn$ centered at the origin, we
have $\deltanorm{ u^k_j } \leq \deltanorm{ u_j }$ for every
$k=1,\ldots,n$.

Given a homogeneous polynomial $P\in\poly_{n,d}$ (of degree $d\geq
1$) and a formal power series mapping $u\colon (\CN,0)\to
(\Cn,0)$, we want to estimate, for any $j\in \N^*$, the
homogeneous part of order $d + j - 1$ of the composition $P\circ
u$ which, in agreement with the former notation, is denoted by
$(P\circ u)_{d+j-1}$.

\begin{lem}\label{l:homcompest}
  Let $P\in\poly_{n,d}$ and assume that $u:(\CN,0)\to(\Cn,0)$ is a formal power series mapping whose decomposition
  into homogeneous terms is given by \eqref{e:decomphomog}, $d\geq 1$.
  Then for every  polydisc $\Delta$ centered at the origin in $\CN$ and every $j\geq 1$, we have
\begin{equation}
  \deltanorm{(P\circ u)_{d+j-1}} \leq \onenorm{P}\sum \frac{d!}{k_1!\dots k_j!} \deltanorm{u_1}^{k_1} \dots \deltanorm{u_j}^{k_j},
  \label{e:homcompest1}
\end{equation}
where the sum goes over all $k_1,\dots,k_j$ such that $k_1 + 2 k_2
+ \dots + j k_j = d + j -1$ and $ k_1 + \dots + k_j = d$, and if
$j\geq 2$, we also have
\begin{equation}
  \deltanorm{(P\circ u)_{d+j-1} - (P^\prime \circ u_1) \cdot u_j } \leq \onenorm{P}\sum \frac{d!}{k_1!\dots k_{j-1}!} \deltanorm{u_1}^{k_1}
  \dots \deltanorm{u_{j-1}}^{k_{j-1}},
  \label{e:homcompest2}
\end{equation}
where the sum goes over all $k_1,\dots,k_{j-1}$, such that $k_1 +
2 k_2 + \dots + (j-1) k_{j-1} = d + j -1$ and $ k_1 + \dots +
k_{j-1} = d$. Here $\left((P^\prime \circ u_1) \cdot u_j\right) $
is the homogeneous polynomial of degree $d+j-1$ given by $\sum_{
\nu=1}^{n} \left(\displaystyle \frac{\partial P}{\partial w_{\nu}}\circ u_1\right) 
u_j^{\nu}$.
\end{lem}
\begin{proof}
  We set $v(t): = \sum_{j\in \N^*} \deltanorm{u_j} t^j\in \C
  [[t]]$. Writing $P(w):=\sum_{|\alpha|=d}P_\alpha w^\alpha$, we
  have
  \begin{align*}
    (P\circ u)_{d+j-1} (\zeta) &= \sum_\alpha P_\alpha  \left[ (u^1(\zeta))^{\alpha_1}\dots (u^n(\zeta))^{\alpha_n}\right]|_{d+j-1} \\
    &= \sum_{\alpha} P_\alpha \sum_{\sum \delta^{r,s}= d+j-1} u^1_{\delta^{1,1_{\,}}}(\zeta) \dots u^1_{\delta^{1,\alpha_1}}(\zeta) u^2_{\delta^{2,1_{\,}}}(\zeta)\dots u^n_{\delta^{n,\alpha_n}}(\zeta)
   \end{align*}
    Hence, taking the norms we get
    \begin{align*}
     \deltanorm{ (P\circ u)_{d+j-1} }&=\max_{\zeta \in \Delta}
    \left|(P\circ u)_{d+j-1} (\zeta)\right|\\
   &\hspace{-1cm}\leq  \sum_{\alpha} |P_\alpha| \sum_{\sum \delta^{r,s}= d+j-1}{\deltanorm{u^1_{\delta^{1,1}}} \dots \deltanorm{u^1_{\delta^{1,\alpha_1}}}
       \deltanorm{u^2_{\delta^{2,1}}}\dots
       \deltanorm{u^n_{\delta^{n,\alpha_n}}}}\\
        &\hspace{-1cm}\leq  \sum_{\alpha} |P_\alpha| \sum_{\sum \delta^{r,s}= d+j-1}{\deltanorm{u_{\delta^{1,1}}} \dots \deltanorm{u_{\delta^{1,\alpha_1}}}
       \deltanorm{u_{\delta^{2,1}}}\dots
       \deltanorm{u_{\delta^{n,\alpha_n}}}}
       \end{align*}
But since the length of every multiindex $\alpha \in \N^N$ is
exactly $d$, we get
\begin{align*}
         \deltanorm{ (P\circ u)_{d+j-1} }&
         \leq \sum_{\alpha} |P_\alpha| \sum_{k_1 + \dots + k_d = d+j-1} \deltanorm{u_{k_1}}\dots\deltanorm{u_{k_d}}\\
         &\leq \onenorm{P} \sum_{k_1 + \dots + k_d = d+j-1} \deltanorm{u_{k_1}}\dots\deltanorm{u_{k_d}} \\
    &\leq \onenorm{P} \frac{1}{(d+j-1)!}\left[\frac{{\rm d}^{d+j-1}}{{\rm dt}^{d+j-1}}(v(t))^d \right]\Biggr|_{t=0} \\
    &\leq \onenorm{P} \sum \frac{d!}{k_1! \dots k_j!} \deltanorm{u_1}^{k_1} \dots \deltanorm{u_j}^{k_j},
  \end{align*}
  the last sum goes over all $k_1,\dots,k_j$ such that $k_1 + 2 k_2 + \dots + j k_j = d + j -1$ and $ k_1 + \dots + k_j = d$, by Proposition~\ref{p:faadebruno}.
  So we have proved \eqref{e:homcompest1}.

The proof of the second estimate  \eqref{e:homcompest2} is similar
to what was done above for \eqref{e:homcompest1}. Indeed note that
one has
\begin{multline} (P\circ u)_{d+j-1} (\zeta) 
  -(P^\prime (u_1(\zeta))) \cdot u_j(\zeta)=\\
     \sum_{\alpha} P_\alpha 
     \sum_{\substack{\sum \delta^{r,s}= d+j-1 \\ 1\leq \delta^{r,s} \leq j-1}} 
     u^1_{\delta^{1,1_{\,}}}(\zeta) 
     \dots u^1_{\delta^{1,\alpha_1}}(\zeta) u^2_{\delta^{2,1_{\,}}}(\zeta)
     \dots
     u^n_{\delta^{n,\alpha_n}}(\zeta).
\end{multline}
Hence following the previous estimates, we get
\begin{align*}
         \deltanorm{ (P\circ u)_{d+j-1}-(P^\prime \circ u_1) \cdot u_j }
         &\leq \onenorm{P} \sum_{\substack{k_1 + \dots + k_d = d+j-1\\ 
	 1\leq k_\nu< j}} \deltanorm{u_{k_1}}\dots\deltanorm{u_{k_d}} \\
    &\leq \onenorm{P} \sum \frac{d!}{m_1! \dots m_{j-1}!} \deltanorm{u_1}^{m_1} \dots \deltanorm{u_{j-1}}^{m_{j-1}},
  \end{align*}
  the last sum going over all nonnegative $m_1,\dots,m_{j-1}$ such that $m_1 + 2 m_2 + \dots + (j-1) m_{j-1} = d + j -1$ and $m_1 + \dots + m_{j-1} = d$.
The proof of Lemma \ref{l:homcompest} is complete.
\end{proof}

We will also need the following lemma, which, although a bit
technical, is crucial for our treatment.
\begin{lem}
  \label{l:fdbcomp} Let $s\in\N^*$, $a>0$, $R>0$, $C>0$, $D>0$ and $\gamma_1>0$. Define $\gamma_0 = 0$, and for $l \geq 2$,
  \begin{equation}
    D \gamma_l := \sum \frac{(s+d)!}{k_1! \dots k_{l-1}!} R a^d \gamma_1^{k_1} 
    \dots \gamma_{l-1}^{k_{l-1}} \, + \, C^{l+s-1};
    \label{e:fdbcomp}
  \end{equation}
  the sum is over all nonnegative integers $k_1,\dots,k_{l-1}$
  satisfying $\sum_{1\leq \nu \leq l-1}k_\nu\geq s$ and
  $k_1 + 2 k_2 + \dots + (l-1) k_{l-1} = s + l - 1$, and we
  set  $s + d := k_1 + \dots + k_{l-1}$. Then the $($formal$)$ power series $\sum_{l\geq 0} \gamma_l t^l$ converges in a neighbourhood of $0$.
\end{lem}
\begin{proof}
  We are going to show that Lemma~\ref{l:fdbcomp} is a consequence of Proposition~\ref{p:faadebruno}. For this, we consider the following analytic
  equation where the unknown is a complex-valued formal power series $\varphi
  (t)$ in one variable :
\begin{equation}\label{e:bla}
(Rs\gamma_1^{s-1} + D) t^{s-1} \varphi (t) - (R(s-1)\gamma_1^{s} +
\gamma_1 D) t^s - \sum_{l>s} C^{l} t^l = \frac{R \varphi(t)^s }{ 1
- a \varphi(t)} ,
\end{equation}
We claim that there exists a unique formal power series solution
$\varphi (t)$ satisfying \eqref{e:bla} with $\varphi (0)=0$ and
$\varphi' (0)=\gamma_1$. For this, set $\varphi (t):=\sum_{l\geq
0}\gamma_lt^l$ and note that to prove the claim it suffices to
show that \eqref{e:bla} holds for a uniquely determined choice of
the coefficients $\gamma_l$, $l\geq 2$. This can be checked by
comparing homogeneous terms of the same order on both sides of
\eqref{e:bla}. We first note that homogeneous terms of order $s$
on both sides of \eqref{e:bla} are identical. For $m\geq 2$, the
homogeneous term of order $m+s-1$ in the left hand side of
\eqref{e:bla} is given by
$$\left((Rs\gamma_1^{s-1}+D)\gamma_m+C^{m+s-1}\right)t^{m+s-1}$$ while on the right hand side
such a term is of the form
\begin{equation}\label{e:left}
\left(Rs\gamma_1^{s-1}\gamma_m+Q_m(\gamma_1,\ldots,\gamma_{m-1})\right)t^{m+s-1},
\end{equation}
for some polynomial $Q_m$. Hence \eqref{e:bla} is satisfied by
$\varphi (t)$ if and only if for every $m\geq 2$,
$D\gamma_m=Q_m(\gamma_1,\ldots,\gamma_{m-1})-C^{m+s}$ which shows
the existence of a unique formal power series solution. Next, setting $\psi (t):=\varphi (t)/t$ and 
$$F(x,t):=(Rs\gamma_1^{s-1} + D)  x - (R(s-1)\gamma_1^{s} +
\gamma_1 D)  - \sum_{l>s} C^{l} t^{l-s} - \frac{R x^s }{ 1
- a t x},$$
it follows from  \eqref{e:bla} that $\psi$ satisfies the equation $F(\psi (t),t)=0$. Since $F$ is convergent and
$\displaystyle \frac{\partial F}{\partial x} (\gamma_1,0)=D\not =0$, from the implicit function theorem we get the convergence of $\psi$ and hence that of $\varphi$. (We take the opportunity here to thank an anonymous referee for
providing us a simple proof of the convergence of $\varphi$.) To complete
the proof of the lemma, we are left to check that the coefficients
$\gamma_l$, $l\geq 2$, satisfy \eqref{e:fdbcomp}. For this, we
want a more precise formula for the coefficient given in
\eqref{e:left}. To the right hand side of \eqref{e:bla}, we apply
Faa di Bruno's formula (Proposition~\ref{p:faadebruno}) with
  \[ h(t) = \varphi(t), \quad g(x) = \frac{R x^s}{1 - a x},\]
  to obtain that
  \begin{equation}\label{e:new}
   (Rs\gamma_1^{s-1} + D) \gamma_m - C^{m+s-1} = \sum_{k_1,\dots,k_m}  \frac{(k + d)!}{k_1!\dots k_m!} \, R a^d \, \gamma_1^{k_1}\dots\gamma_m^{k_m}
   \end{equation}
  where the sum goes over all nonnegative integers
  $k_1,\ldots,k_m$ such that $k_1+\ldots+k_m\geq s$, $\sum_{1\leq j\leq m} j k_j = s+ m
  -1$ and where  $s + d := \sum_{1\leq j\leq m} k_j$. (Indeed, one easily checks that since $g\circ h$ vanishes to order $k$,
  the $k_{m+j}$ in \eqref{e:faadebruno} for
  $1\leq j \leq s-1$ vanish, so we don't have to write them.)
  Moreover, in view of \eqref{e:left},
  the term on the right hand side of \eqref{e:new} containing $\gamma_m$ corresponds to $k_1 = s -1$, $k_2=0,\ldots, k_{m-1}=0$,
  and $k_m = 1$, so that it
  is equal to $R s \gamma_1^{s-1} \gamma_m$. Hence, \eqref{e:new} gives
  \eqref{e:fdbcomp}. The proof of Lemma \ref{l:fdbcomp} is
  complete.
\end{proof}

\subsection{Reducing the initial system to a simpler
one}\label{ss:reduce}
%Before proceeding to the proof of Theorem \ref{t:parcn1}, we first
%need to bring the holomorphic map $A\colon (\C^n,0)\to \C^n$ into
%some convenient canonical normal form.
Let \[ A(w) = \left( A^1 (w), \dots, A^n (w)
\right) \] be a $\Cn$-valued holomorphic map defined a
neighbourhood of the origin in $\Cn$. Without loss of generality,
we may assume that $A(0)=0$. We decompose each $A^j$ into
homogeneous polynomials
\begin{equation}\label{e:Aj}
A^j (w) = \sum_{k> k_j} A_k^j (w), \quad A^j_k (tw) = t^k A_j^k
(w),\ {\rm for}\ {\rm all}\ (t,w)\in \C\times C^n.
\end{equation}
 We will write
$A_\nu = \left( A^1_{k_1+1+\nu},\dots,A^n_{k_n+1 + \nu} \right)$
for all $\nu \in \N$. The first problem we have to face before proceeding to the proof of Theorem~\ref{t:parcn1} 
is due to the fact that even if $A$ is of generic full rank, its lowest order homogeneous terms need not be, in general, 
of full rank too. This difficulty will be overcome by using the following trick contained in
\cite{L4}. For the reader's convenience, we recall the proof.

\begin{lem} Let $A$ be as above and assume that $A$ is generically of full rank. Then there exists a
polynomial map $P: (\Cn,0)\to (\Cn,0)$ with the
  property that the holomorphic map $B\colon (\Cn,0)\to (\Cn,0)$ given by $B(w) := P ( A(w) )$
  has its lowest order homogeneous terms $B_0 (w)$ of generic full rank.
  \label{l:simplify1}
\end{lem}
\begin{proof}
  In what follows, for any local holomorphic map $h:(\Cn,0)\to (\Cn,0)$, we denote
  by $h'(w)$ its Jacobian matrix at the point $w$.

  Let $e$ be the order of vanishing of $\det A^\prime (w)$. We define $D = e - \sum_{1\leq j\leq n} k_j$. Then  clearly $A$ satisfies
  $\det A_0^\prime \neq 0$ if and only if $D = 0$. If $D=0$ we are, 
  therefore, done and hence we may assume
  that $\det A_0^\prime = 0$. In this case, the polynomials $A_{k_j + 1}^j$ for $j=1,\dots,n$
  are algebraically dependent over $\C$ (see e.\ g.\  \cite{HP}), that is, there exists a nonzero polynomial $p\in \C [ x_1 ,\dots x_n ]$ such that
  $p(A_{k_1+1}^1,\dots,A_{k_n+1}^n) = 0$. Since each $A_{k_j+1}^j$ is homogeneous of degree $k_j+1$,
  we may choose $p$ to be weighted homogeneous (where $x_j$ has weight $k_j+1$) of lowest possible degree
  $f$.  This means that $p(x_1,\ldots,x_n)=\sum_{\alpha\in \N^n}c_\alpha x^\alpha$
  with $\alpha_1(k_1+1)+\ldots+\alpha_n (k_n+1)=f$ and therefore for all
  $t\in \C$,
  $p(t^{k_1+1}x_1,\ldots,t^{k_n+1}x_n)=t^fp(x_1,\ldots,x_n)$.  Without loss
  of generality, we may also assume that $p_{x_1}$ is nonzero. We now
  consider the map $\tilde A \colon (\Cn,0)\to (\Cn,0)$ given by $\tilde A^1 = p(A)$, $\tilde A^j = A^j$, $j>1$ and
  claim that for this map, the associated number $\tilde D$ is strictly smaller than that associated to $A$, namely $D$.
  Indeed, first note that we have
  $$\det \tilde A^\prime(w)=p_{x_1}(A^1(w),\ldots,A^n(w))\cdot \det A^\prime(w).$$
Moreover, by our choice of $p$, the order of vanishing of
$p_{x_1}(A^1,\ldots,A^n)$ is exactly $f-(k_1+1)$ and, therefore, the
order of vanishing of $\det \tilde A^\prime$ is $e + f - (k_1 +
1)$. For $j=1,\ldots,n$, denote by $\tilde k_j+1$ the order of
vanishing of $\tilde A^j$. Then by our choice of $p$, one has
$\tilde k_1+1>f$ while $\tilde k_j = k_j$ for $j >1$. Thus we
obtain $\tilde D = e + f - (k_1 +1) - \tilde k_1 - \sum_{1< j\leq
n} k_j <  e - \sum_{1\leq j\leq n} k_j = D$, proving the claim.
  This in turn means that after a finite number of these 
  procedures, we get at a map $B$ satisfying
   the property claimed in the Lemma.
\end{proof}

\subsection{Completion of the proof}\label{ss:conclpfthm31}
By Lemma~\ref{l:simplify1}, there exists a polynomial map $P\colon
\C^n\to \C^n$ (vanishing at the origin) such that $B(w)=P(A(w))$
has lowest order homogeneous terms $B_0(w)$ of generic full rank
$n$. Now note that if we have proved Theorem \ref{t:parcn1} for
$B$, then Theorem \ref{t:parcn1} easily also holds for $A$. So we
may assume that the lowest order terms of $A$, denoted as above by
$A_0 (w)=A_0(w^1,\ldots,w^n)$, are generically independent.
Furthermore, by taking every $A^j$ as defined in \eqref{e:Aj} to
some appropriate power, we may assume that $k_1 = \dots = k_n =:
\ell_0\geq 0$. In addition, since after a linear change of
coordinates $\widetilde w=(\widetilde w^1,\ldots,\widetilde
w^n)={\mathcal L}w$, $\det A_0$ may be normalized in the
$\widetilde w^n$ direction, we claim that we may assume that such
a property already holds in the original $w$-coordinates. Indeed,
if in the $\widetilde w$-coordinates we have found a
parametrization $\Psi$ as given in \eqref{e:formalstuff} whose
homogeneous parts  satisfy all the conditions of Theorem
\ref{t:parcn1}, then it suffices to set $\psi
(b,\lambda):={\mathcal L}^{-1}\cdot \Psi (b,{\mathcal L}\cdot
\lambda)$ for $(b,\lambda)\in \holmaps \times \glnc$ to get the
right parametrization  (whose homogeneous parts  satisfy all the
requirements of Theorem \ref{t:parcn1}) in the original
$w$-coordinates. This proves the claim.

After all these reductions, we may start with the formal part of
the proof and first construct a sequence $\tilde{p}_\ell \colon
\C^n\times J_0^{\ell_0+\ell}(\Cn)\to \Cn$ of polynomial maps,
$\ell \geq 1$, such that for every germ of a holomorphic map
$b\colon (\C^n,0)\to \C^n$ if there exists a biholomorphism
$u\colon (\Cn,0)\to (\Cn,0)$ with $u^\prime(0)$ being the identity
and $A(u(z)) =b(z)$, then $u_\ell (z) = \tilde{p}_\ell
(z,\jetm{0}{\ell_0 + \ell} b)$ for all $\ell \geq 1$. Here we
recall that $u_\ell$ denotes the homogeneous part of order $\ell$
of $u$. In what follows we also use the notation $\theta
(z)|_\ell$ to denote the homogeneous part of order $\ell$ of a
formal power series $\theta$. Moreover, for every $k\in \N^*$, we
introduce coordinates $\Lambda^k=(\Lambda_{\alpha})_{|\alpha|\leq
k}$, $\alpha \in \N^n$, in the jet space
 $J_0^k(\Cn)$ associated to a choice of coordinates
 $z=(z_1,\ldots,z_n)$ in $\Cn$. By Lemma \ref{l:multmap}, we may choose
  for each integer $\ell \geq 2$ a linear map $B_\ell \colon\poly^n_{\ell_0 +
\ell} \to \poly^n_\ell$, which is a left inverse to multiplication
by $A_0^\prime$. We set $\tilde{p}_1 (z): = z$ and define
inductively the ${\tilde p}_\ell$, $\ell \geq 2$, as follows
\begin{equation}
\begin{cases}
\begin{aligned}
  \tilde{p}_\ell (z,\Lambda^{\ell_0+\ell}): =&
  B_\ell\left(\sum_{|\alpha|=\ell_0 + \ell}\frac{\Lambda_\alpha}{\alpha !} z^\alpha -
  A(P^{\ell - 1} (z,\Lambda^{\ell_0 + \ell - 1} ))|_{\ell_0 + \ell}  \right),\cr \label{e:pelldef}
P^{\ell -1}(z,\Lambda^{\ell_0+\ell -1}):=& \sum_{1\leq k\leq \ell
-1}{\tilde p}_k(z,\Lambda^{\ell_0+k}).
\end{aligned}
\end{cases}
\end{equation}
It is easy to see that we also have
\begin{equation}\label{e:theendwhen?}
\tilde{p}_\ell (z,\Lambda^{\ell_0+\ell}) =
B_\ell\left(\sum_{|\alpha|=\ell_0 +
\ell}\frac{\Lambda_\alpha}{\alpha !} z^\alpha -
  \sum_{\nu =0}^{\ell -1}A_\nu(P^{\ell - 1} (z,\Lambda^{\ell_0 + \ell - 1} ))|_{\ell_0 + \ell}  \right)
\end{equation}
for $\ell\geq 2$, and hence by induction from
\eqref{e:theendwhen?} and \eqref{e:pelldef} that the ${\tilde
p}_\ell$ are indeed polynomial mappings of their arguments and
homogeneous of degree $\ell$ in their first factor. Now assume
that $b\colon (\C^n,0)\to \C^n$ is a germ of a holomorphic map and
that $u\colon (\Cn,0)\to (\Cn,0)$ is a germ of a biholomorphism
satisfying $A(u(z))= b(z)$ with $ u^\prime(0)$ being the identity.
We claim that for every $\ell\geq 1$, $u_l(z)={\tilde
p}(z,j_0^{\ell_0+\ell}b)$. First note that $u_1={\tilde p}_1$.
Suppose now that for $\ell \geq 2$, we have $u_k(z)={\tilde
p}_k(z,j_0^{\ell_0+k}b)$ for $k<\ell$ and set $U_\ell (z) =
\sum_{1\leq k\leq \ell} u_\ell (z)$. We then have
\begin{equation}  \label{e:uprogr}
  \begin{aligned}
    b_{\ell_0 + \ell} (z) & = A(u(z))|_{\ell_0 + \ell} \\
    & = \sum_{\nu=0}^{\ell-1} A_\nu (u(z))|_{\ell_0 + \ell} \\
    & =\sum_{\nu=0}^{\ell-1} A_\nu (U_\ell(z))|_{\ell_0 + \ell}.\\
  \end{aligned}
\end{equation}
Noticing that $A_0(U_\ell(z))|_{\ell_0+\ell}=A_0^\prime (z)\cdot
u_\ell (z)+A_0(U_{\ell -1}(z))|_{\ell_0+\ell}$ and that for $1\leq
\nu \leq \ell -1$, $A_\nu(U_\ell(z))|_{\ell_0+\ell}=A_\nu(U_{\ell
-1 }(z))|_{\ell_0+\ell}$, we obtain that
\begin{equation}
 b_{\ell_0+\ell}(z)= A_0^\prime (z)\cdot u_\ell (z) +
\sum_{\nu=1}^{\ell-1} A_\nu (U_{\ell -1}(z))|_{\ell_0 + \ell},
\end{equation}
and since $u_k(z)={\tilde p}_k(z,j_0^{\ell_0+k}b)$ for $k<\ell$ by
assumption, we get
$$A_0^\prime(z)\cdot u_\ell(z) = b_{\ell_0 + \ell}
(z) - \sum_{\nu=0}^{\ell-1} A_\nu (P^{\ell -1}(z,j_0^{\ell_0+\ell
-1}b))|_{\ell_0 + \ell},$$ i.e.
\begin{equation}\label{e:trip}
 u_\ell = B_\ell \left(b_{\ell_0 + \ell}
(\cdot) - \sum_{\nu=0}^{\ell-1} A_\nu (P^{\ell
-1}(\cdot,j_0^{\ell_0+\ell -1}b))|_{\ell_0 + \ell}\right).
\end{equation}
and consequently $u_\ell (z) = \tilde{p}_\ell
(z,j_0^{\ell_0+\ell}b)$. This proves the claim.

We are now going to prove the convergence statement for the
$\tilde{p}_\ell$. In this context, this means that there exists a
polydisc $\Delta \subset \Cn$ centered at the origin such that
given any $C>0$, there exists $K>0$ such that whenever $b\colon
(\Cn,0)\to \Cn $ is a germ of holomorphic satisfying
$\deltanorm{b_k} \leq C^k$ for all $k\in\N$ where
$b(z)=\sum_{k\geq 0}b_k(z)$ is the decomposition of $b$ into
homogeneous terms, then
$\deltanorm{\tilde{p}_\ell(\cdot,j_0^{\ell_0+\ell}b)} \leq K^\ell$
for all $\ell\in\N$.  Let $\Delta=\Omega \subset \Cn$ be the
polydisc given by Lemma \ref{l:multmap} and $D>0$ also given by
the same lemma (and depending only on $A_0^\prime$) such that the
norm of each linear map $B_\ell \colon\poly^n_{\ell_0 + \ell} \to
\poly^n_\ell$ does not exceed $1/D$, for all $\ell \geq 2$. Let
$b$ be as above satisfying $\deltanorm{b_k}\leq C^k$ for all $k\in
\N$. In what follows we assume $b$ to be fixed and then suppress
the dependence on $b$ when writing ${\tilde p}_\ell$ and $P^{\ell
-1}$. Since $A$ is convergent, there exist $R>0$ and $a>0$ such
that $\onenorm{A_j} \leq R a^j$ for all $j\in \N$. We also set
$m_0: = \ell_0 +1\geq 1$. Let $\ell \geq 2$. By \eqref{e:trip}, we
have
\begin{equation}\label{e:sunshine}
\deltanorm{{\tilde p}_{\ell}}\leq
D^{-1}\left(\deltanorm{b_{\ell_0+\ell}}+\deltanorm{ (A_0 \circ
P^{\ell -1})|_{\ell_0 + \ell}}+\deltanorm{\sum_{\nu=1}^{\ell-1}
(A_\nu \circ P^{\ell -1})|_{\ell_0 + \ell}}\right).
\end{equation}
Since $A_0(P^\ell(z))|_{\ell_0+\ell}=A_0^\prime (z)\cdot {\tilde
p}_\ell (z)+A_0(P^{\ell -1}(z))$, \eqref{e:sunshine} implies
\begin{equation}
\deltanorm{{\tilde p}_{\ell}}\leq
D^{-1}\left(C^{\ell_0+\ell}+\deltanorm{(A_0\circ
P^l)|_{\ell_0+\ell}-A_0^\prime \cdot {\tilde p}_\ell
}+\deltanorm{\sum_{\nu=1}^{\ell-1} (A_\nu \circ P^{\ell
-1})|_{\ell_0 + \ell}}\right).
\end{equation}
Applying \eqref{e:homcompest1} and \eqref{e:homcompest2} of Lemma
\ref{l:homcompest} in our context and using the fact that ${\tilde
p}_1(z)=z$ yields
\begin{equation}\label{e:out}
\begin{aligned}
\deltanorm{{\tilde p}_{\ell}} &\leq
D^{-1}\left(C^{\ell_0+\ell}+\sum \frac{(m_0+d)!}{k_1!\ldots
k_{\ell -1}!}\onenorm{A_d}\deltanorm{{\tilde p}_1}^{k_1}\ldots
\deltanorm{{\tilde p}_{\ell -1}}^{k_{\ell -1}}\right)\\
&\leq D^{-1}\left(C^{m_0+\ell-1}+\sum \frac{(m_0+d)!}{k_1!\ldots
k_{\ell -1}!}Ra^d\deltanorm{{\tilde p}_1}^{k_1}\ldots
\deltanorm{{\tilde p}_{\ell -1}}^{k_{\ell -1}}\right),
\end{aligned}
\end{equation}
where $k_1 + \dots + k_{\ell-1}=m_0+d$ and the sum goes over all
nonnegative integers $k_1,\dots, k_{\ell-1}$ with $\sum_{1\leq
j\leq \ell -1} j k_j = m_0 + \ell -1$. We shall now use
Lemma~\ref{l:fdbcomp} with $\gamma_1: =\deltanorm{\tilde{p}_1}$
(that is independent of $b$), $R$, $a$, $C$, $D$ as above and
$s=m_0=\ell_0+1$. Let $(\gamma_\ell)_{\ell \geq 0}$ be given by
Lemma \ref{l:fdbcomp} and associated to the above data. We claim
that for all $\ell \geq 1$, $\deltanorm{{\tilde p}_\ell}\leq
\gamma_\ell$. Indeed, by construction this is the case for $\ell
=1$. Let us  assume, therefore, that $\deltanorm{\tilde{p}_k}\leq
\gamma_k$ for $k<\ell$. Then \eqref{e:out} and \eqref{e:fdbcomp}
imply that
\begin{equation}
  \deltanorm{ \tilde{p}_{\ell}  }\leq D^{-1}\left(C^{m_0+\ell-1}+\sum \frac{(m_0+d)!}{k_1!\ldots
k_{\ell -1}!}Ra^d{\gamma_1}^{k_1}\ldots \gamma_{\ell -1}^{k_{\ell
-1}}\right) = \gamma_\ell,
  \label{e:thmpfest2}
\end{equation}
which proves the claim by induction. Since the series $\sum
\gamma_k t^k$ converges by Lemma \ref{l:fdbcomp}, we see that
there exists $K>0$ such that $\deltanorm{\tilde{p}_\ell}\leq
K^\ell$ for $\ell \in \N^*$.

We shall now construct the $p_\ell$ from the $\tilde p_\ell$ for
all $\ell \geq 1$. For any linear map $N\colon \Cn\to \Cn$ and for
every positive integer $k$ we denote by $N^k_*\colon J_0^k(\Cn)\to
J_0^k(\Cn)$ the (unique linear) map satisfying $j_0^k(f\circ
N)=N^k_*(j_0^kf)$ for every formal power series mapping $f\colon
(\Cn,0)\to \Cn$. (Note that by identifying jets of order $k$ with
Taylor polynomials of degree $k$, one may also view the mapping
$N^k_*$ as acting on such polynomials.) For every $\ell \geq 1$
and every $(z,\genmat,\Lambda^{\ell_0+\ell})\in \Cn\times \glnc
\times J_0^{\ell_0+\ell}(\Cn)$, we set
\begin{equation}\label{e:ddd} r_\ell
(z,\genmat,\Lambda^{\ell_0+\ell}) := \tilde{p}_\ell(\genmat z ,
(\genmat^{-1})^{\ell_0+\ell}_{*}(\Lambda^{\ell_0+\ell} )).
\end{equation}
It is easy to see that each $r_\ell$ is a homogeneous polynomial
map of degree $\ell$ in $z$, whose components are rational
functions of the other variables. In fact, even more can be said:
by Cramer's rule we see that there exists a nonnegative integer
$\kappa(\ell)$ such that $(\det \genmat)^{\kappa(\ell)} r_{\ell}
(z,\genmat, \Lambda^{\ell_0 + \ell})=:p_\ell (z,\genmat,
\Lambda^{\ell_0 + \ell})$ is a polynomial in all its variables.
Note that $r_1(z,\lambda,\Lambda^{\ell_0+1})=\genmat z$ by
construction which proves \eqref{e:7}. Moreover, for every germ of
a holomorphic map $b \colon (\Cn,0)\to \C^n$, we claim that if
$u\colon (\Cn,0)\to (\Cn)$ is a germ of a biholomorphism
satisfying $A(u(z))=b(z)$, then necessarily for every $\ell \geq 1
$, $u_\ell(z)=r_\ell (z,u^\prime (0),j_0^{\ell_0+\ell}b)$. Indeed,
set $\genmat:=u^\prime(0)$, and note that $v(z):=u(\genmat^{-1}z)$
is a biholomorphism tangent to the identity which satisfies
$A(v(z))=b(\genmat^{-1}z)$. Hence for every $\ell \geq 1$, the
homogeneous part of order $\ell$ of $v$ denoted $v_\ell$ satisfies
$v_\ell (z)={\tilde p}_{\ell} (z,j_0^{\ell_0+\ell}(b\circ
\genmat^{-1}))={\tilde p}_{\ell}
(z,(\genmat^{-1})^{\ell_0+\ell}_*(j_0^{\ell_0+\ell}b))$ and
therefore $v(\genmat z)=u(z)$ implies that $u_\ell (z)={\tilde
p}_{\ell} (\genmat
z,(\genmat^{-1})^{\ell_0+\ell}_*(j_0^{\ell_0+\ell}b))$ which is
the required conclusion in view of \eqref{e:ddd}. To complete the
proof of Theorem \ref{t:parcn1}, it remains to check the
convergence statement for the $r_\ell$. So let $\Delta\subset \Cn$
be the polydisc given above and $L\subset\glnc$ be a given compact
subset. Let $C>0$ and $b\colon (\Cn,0)\to \Cn$ be a germ of
holomorphic map satisfying $\deltanorm{b_k}\leq C^k$ for all $k\in
\N$. Since there exists  $C_1>0$ such that
$\vnorm{\genmat^{-1}}\leq C_1$ for all $\lambda \in L$ and since
$(\genmat^{-1})^k_{*} (b_k)$ is the homogeneous part of order $k$
of $b\circ \genmat^{-1}$, we have $\deltanorm{(\genmat^{-1})^k_{*}
(b_k) } \leq \vnorm{b_k}_{C_1 \Delta}$ where $C_1\Delta$ is the
polydisc of multiradius $C_1\rho$ if $\rho$ is the multiradius of
$\Delta$. Using the homogeneity of each $b_k$ we therefore get
$\deltanorm{(\genmat^{-1})^k_{*}(b_k)}\leq C_1^k
\deltanorm{b_k}\leq C_1^k C^k =:\tilde{C}^k$, for all $\genmat\in
L$ and $b$ as above. Hence $\deltanorm{(b\circ
\genmat^{-1})|_k}\leq {\tilde C}^k$ for all $k\in \N$ and by the
first part of the proof, there exists $\tilde{K}>0$ such that
$$\deltanorm{\tilde{p}_\ell (\cdot,j_0^{\ell
+ \ell_0} (b\circ \genmat^{-1}))}= \deltanorm{\tilde{p}_\ell
(\cdot,(\genmat^{-1})^{\ell +\ell_0}_{*} (j_0^{\ell + \ell_0}
b))}\leq \tilde{K}^\ell.$$ Furthermore there exists $C_2>0$ such
that $\vnorm{\genmat}\leq C_2$ for all $\genmat\in L$, so we have
\begin{equation*}\begin{split}
  \deltanorm{r_{\ell} (\cdot,\genmat,\jetm{0}{\ell_{0}+\ell }b) } &=
\deltanorm{\tilde{p}_{\ell}(\genmat \cdot,(\genmat^{-1})^{\ell+\ell_0}_{*}(\jetm{0}{\ell+\ell_0} b)) } \\
&\leq
\vnorm{\tilde{p}_{\ell}(\cdot,(\genmat^{-1})^{\ell+\ell_0}_{*}(\jetm{0}{\ell+\ell_0}
b))}_{C_2 \Delta}\\
&\leq
C_2^{\ell}\deltanorm{\tilde{p}_{\ell}(\cdot,(\genmat^{-1})^{\ell+\ell_0}_{*}(\jetm{0}{\ell+\ell_0}
b )) } \leq K^\ell,
\end{split}
\end{equation*}
where $K= \tilde{K} C_2$. The proof of Theorem \ref{t:parcn1} is
complete.

\section{Proof of Theorem \ref{t:parcnfull}}\label{s:holomorphicity}

\subsection{Holomorphy of $\Psi(b,\lambda)$ }\label{ss:end?}
To each {\em fixed} choice of polynomial mappings $(p_\ell)$ and integers $(\kappa (\ell))$ satisfying the conclusions of Theorem~\ref{t:parcn1}, we consider the associated  map $\Psi \colon \holmaps \times
\glnc \to \biholmaps$ given in \eqref{e:formalstuff} and we show here that $\Psi$ is
holomorphic furnishing the proof of Theorem~\ref{t:parcnfull}.
Here we recall that  holomorphy of a map $f\colon X\supset \Omega\to Y$ between locally
convex topological complex vector spaces $X,Y$ with $\Omega$ open in $X$ means that $f$ is holomorphic if
it is continuous and its composition with any continuous linear
functional on $Y$ is holomorphic along finite dimensional affine
subspaces of $X$ intersected with $\Omega$ (this last property
being called {\em G\^ateaux-holomorphy}). Note that the fact
that $\Psi$ takes its values indeed in $\biholmaps$ follows from
Theorem \ref{t:parcn1} and the following well-known elementary
fact.

\begin{lem}
  \label{l:normsconv} Let $\phi=\phi(z)$ be a formal power series in
  $z=(z_1,\ldots,z_n)$. Consider the decomposition of $\phi$ into
  homogeneous polynomial terms and in Taylor series $\phi(z)=\sum_{k\geq 0} \phi_k(z) = \sum_{\alpha \in \N^n} \phi_\alpha z^\alpha$.
  Then the following are equivalent:
  \begin{enumerate}
  \item[{\rm (i)}] $\phi$ is convergent;
  \item[{\rm (ii)}] There exists $C>0$ such that for all $\alpha \in \N^n$, $|\phi_\alpha|\leq C^{|\alpha|}$; \label{cest}
  \item[{\rm (iii)}] There exists a polydisc $\Delta\subset \Cn$ centered at the origin and
  $D>0$ such that for every positive integer $k$, $\deltanorm{\phi_k}\leq D^{k}$; \label{dest}
  \item[{\rm (iv)}] For every polydisc $\Delta \subset \Cn$ centered at the origin, there exists $K>0$
    such that for every positive integer $k$, $\deltanorm{\phi_k}\leq
  K^{k}$.\label{kest}
  \end{enumerate}
  All the power series $\phi$ satisfying either of the inequalities in
  {\rm (ii)}, {\rm (iii)} and {\rm (iv)}  have the property that their region of convergence contains
  some neighbourhood of the origin only depending on $C$ or $D$
  or $K$, respectively.
\end{lem}

For $r>0$, we let $\holmapsn_r$ be the Fr\'echet space of
holomorphic maps $b: \Delta_r\to\Cn$, where $\Delta_r = \{z\in\Cn
\colon |z_j|<r \text{ for all } j\}$, endowed with the topology of
uniform convergence on compacts subsets of $\Delta_r$. It is in
our case generated by the seminorms
\[\forall f=(f^1,\ldots,f^n)
\in \holmapsn_r,\quad \vnorm{f}_s= \max_{1\leq j\leq n}
\sup_{z\in\Delta_s} |f^j(z)|, \quad s<r, \]
If we choose a
sequence $r_j>0$ converging monotonely to $0$, the space
$\holmapsn$ is a locally convex space which is the
inductive limit of the sequence of
Fr\'echet spaces $\holmapsn_{r_j}$, or, equivalently, of the 
Banach spaces $H^{\infty} (\Delta_{r_j})^n$. We denote by $\holmapsn_r^0$
(resp.\ $\holmapsn^0$) the closed subspace of $\holmapsn_r$ (resp.\
of $\holmapsn$) consisting of those elements preserving the origin.
We also let $\biholmaps_r$ be the open subset of $ \holmapsn^0_r$
consisting of those maps $b$ satisfying $\det b^\prime (0) \neq
0$. For each $r>0$, we have a continuous linear injection
$\iota_r\colon \holmapsn_r \to \holmapsn$, through which we identify
$\holmapsn_r$ with a subspace of $\holmapsn$. Note though that
the image of $\iota_r$ in $\holmapsn_s$ is {\em not closed}
for $s<r$; hence, $\holmapsn$ is not a strict inductive limit
of the spaces $\holmapsn_{r_j}$. However,  in our setting, 
the inclusion mappings for the inductive limit are
{\em compact} (by Montel's Theorem). For this special case of an inductive 
limit
some 
powerful results from functional analysis are available; the main one is that this 
inductive limit is {\em regular} (cited as Fact 1 below, which implies Fact 3) and that the
topology induced on it as an inductive limit of locally convex spaces agrees with 
the topology induced on it as an inductive limit of topological spaces (this implies
Fact 2 below). In what follows, when we identify
subsets of $\holmapsn_r$ (resp.\ of $\holmapsn^0_r$) as subsets of
$\holmapsn$ (resp.\ of $\holmapsn^0$), we always use the canonical
embedding.
These results can be found in e.g. \cite[\S 23, p. 132]{FlWlLN} or \cite[p. 26-27]{Kom} 
(and in this particular setting seem to originate from work by e Silva \cite{eSilva}
and Raikov \cite{Raikov}). For the convenience of the reader, 
we summarize the facts we are going 
to use so we can refer to them later:
\begin{enumerate}
  \item[Fact 1:] \cite[2.2. Satz, p. 136]{FlWlLN}
    A subset $B\subset\holmapsn$ is bounded if and only if there exists $r>0$ such that $B\subset\holmapsn_r$ and
    $B$ is bounded in $\holmapsn_r$.
  \item[Fact 2:] A mapping $f:\holmapsn\to Y$, where $Y$ is a topological space, is continuous if and only if it is
    sequentially continuous.
  \item[Fact 3:] \cite[2.10. Korollar 3, p.139]{FlWlLN} 
    A sequence $(y_n)$ converges to $y$ in $\holmapsn$ if and only if there exists an $r>0$ such that
    $y_j,y\in\holmapsn_r$ and the sequence converges there.
\end{enumerate}

Fact 2 is an immediate consequence of Fact 1, Fact 3 and the fact that $\holmaps_r$ is Fr\'echet and thus 
metrizable.
We can now prove the following result concerning the parametrization
$\Psi$, which essentially reduces the study of $\Psi$ as a map between spaces of germs
to that of a map between the Fr\'echet spaces.

\begin{lem}\label{l:addendum}
  With the assumptions of Theorem~{\rm \ref{t:parcn1}} and in the above mentioned setting, let $\Psi\colon \holmapsn \times \glnc \to\biholmaps$ be the
  operator given by  {\rm \eqref{e:formalstuff}}. 
  Then there exists a family $({\mathcal V}_R)_{R>0}$, where for each $R>0$, 
  ${\mathcal V}_R$ is an open subset of $\holmapsn_R$, satisfying the following properties:
  \begin{enumerate}
\item[{\rm (i)}] for any given bounded set $B\subset \holmapsn$, there exists
   $S=S(B)>0$ such
  that $B\subset \mathcal{V}_{S}$;
  \item[{\rm (ii)}] for every $R>0$ and for each relatively compact open subset
  $\Omega\subset \glnc$  there exists an $\widetilde R >0$ such that
  \[ \Psi|_{\mathcal{V}_R \times \Omega} \colon {\mathcal{V}_R\times\Omega} \to \biholmaps_{\widetilde R}
  \subset \holmapsn_{\widetilde R} \]
  is holomorphic.
\end{enumerate}
\end{lem}

\begin{proof}
 For every $R>0$, let $\mathcal{V}_R$ be defined by
 \[ \mathcal{V}_R = \{ b\in\holmapsn_R \colon \vnorm{b-b(0)}_{R/2} < 1,\ 
 \vnorm{b(0)}_{R/2}<\frac{1}{R} \}.\]
  This is clearly an open subset of $\holmapsn_R$. From the usual Cauchy estimates together with Fact 1,
  it is easy to see that every bounded set
  is contained in a set $\mathcal{V}_R$ for some $R>0$.

  Next, for every $R>0$, by the Cauchy estimates, we easily have
  \[ \forall b\in {\mathcal V}_R,\ \forall d\in \N^*,\quad
  \vnorm{b_d}_{R/2} \leq \left( \frac{nR}{2} \right)^d\ \left({\rm and}\ \vnorm{b_0}_{R/2}<\frac{1}{R}\right).\]
  (Recall here that $b_d$ is the homogeneous part of order $d$ of $b$.)
  By Theorem~\ref{t:parcn1} (ii) (which holds automatically for all polydiscs centered at the origin),
  there exists $K>0$ (depending on $R$) such that
  \[ \vnorm{\Psi(b,\lambda)|_d}_{R/2}\leq K^d, \quad \forall b\in \mathcal{V}_R, \quad \forall \lambda\in\Omega,\
  \forall d\in \N^*.\]
  Choosing $\widetilde{R}\leq R/4K$, we see that for each $(b,\lambda)\in\mathcal{V}_R\times\Omega$,
  we have
  \begin{equation}\label{+}
  \vnorm{\Psi(b,\lambda)}_{\widetilde R} \leq 1,
  \end{equation}
  and $\Psi(b,\lambda)\in \holmapsn_{\widetilde R}$. Moreover,
  setting
  $\Psi^{k}(b,\lambda):=\sum_{\ell \leq k}\Psi(b,\lambda)|_{\ell}$ for every $k\geq 1$ and $(b,\lambda)\in {\mathcal V}_{R}
  \times \Omega$, we
  have by the above choices that $(\Psi^{k})$ converges uniformly to $\Psi$ in ${\mathcal H}_{\widetilde R}$ i.e. that
  for all $r\leq \widetilde R$, $\sup_{(b,\lambda)\in {\mathcal V}_{R}\times \Omega}\vnorm{\Psi^{k}(b,\lambda)-\Psi(b,\lambda)}_{r}\to
  0$ when $k\to +\infty$. Denote by ${\mathcal M}_n(\C)$ be the space of $n\times n$ matrices and let
  $E\subset {\holmapsn}_{R}\times {\mathcal M}_n(\C)$ be any affine
  subspace and $L\colon {\holmapsn}_{\widetilde R}\to \C$ be any continuous linear
  functional. In view of the explicit form of the $\Psi^{k}$ given by Theorem \ref{t:parcn1},
  we see that $L\circ \Psi^{k}$ is clearly holomorphic (and rational) along $E\cap(\mathcal{V}_R\times\Omega)$,
  and since $L\circ \Psi^k$ uniformly converges to $L\circ \Psi$ on ${\mathcal V}_R\times
  \Omega$, it follows that $L\circ \Psi$ is holomorphic along $E\cap(\mathcal{V}_R\times\Omega)$. Hence,
  $\Psi|_{\mathcal{V}_R \times \Omega} \colon {\mathcal{V}_R\times\Omega} \to \holmapsn_{\widetilde R}$ is G\^ateaux-holomorphic.
  To conclude  that it is holomorphic,
  it remains to check that it is continuous. In fact, in view of
  \cite[Example 3.8(g)]{Dineenbook}, since $\holmapsn_R$ is metrizable,
  we need only to show that
  for each continuous linear functional
  $L\colon {\holmapsn}_{\widetilde R}\to \C$ the map $L\circ \Psi$ is
  holomorphic. Since this latter is already G\^ateaux-holomorphic, in view
  of \cite[Proposition~3.7]{Dineenbook}
  it is enough to show that $L\circ \Psi$ is bounded on $\mathcal{V}_R\times \Omega$.
  Since $\Psi (\mathcal{V}_R\times \Omega)$ is clearly a bounded
  set in ${\holmapsn}_{\widetilde R}$ in view of \eqref{+}, we
  reach the desired conclusion. The proof of the lemma is
  complete.
\end{proof}

\begin{proof}[Completion of the proof of Theorem {\rm \ref{t:parcnfull}}]
It is now not difficult to derive Theorem \ref{t:parcnfull} from
Lemma \ref{l:addendum}. Firstly, it follows from the
lemma that $\Psi\colon {\holmapsn}\times \glnc \to {\biholmaps}\subset {\holmapsn}$ is G\^ateaux-holomorphic. It remains
therefore to check that $\Psi$ is continuous (between the right
topological spaces). But Lemma \ref{l:addendum} together with Fact
2 and Fact 3 from above implies that $\Psi$ is indeed continuous:
If $(y_n)\subset \holmapsn\times\glnc$ is a convergent sequence
with limit $y$, there exists $R>0$ and a relatively compact open
set $\Omega\subset \glnc$ such that $(y_n)\subset
\mathcal{V}_R\times\Omega$ (and
$y\in\mathcal{V}_{R}\times\Omega$). The lemma then implies that
$\Psi(y_n) \to \Psi(y)$ and we are done.\end{proof}

\subsection{A further application of Theorems \ref{t:parcnfull} and
\ref{t:parcn1}}\label{s:onemore}
We want here to give
the following continuous version of Corollary~\ref{c:corparcn1}.

\begin{thm}
  \label{c:corparcn2}   Let $A:(\Cn,0)\to \Cn$ be a germ of a holomorphic map
 of generic rank $n$, $X$ a topological space and
  $b=b(z,\omega)$ be a $\Cn$-valued continuous map defined
  on an open neighbourhood $V$ of $\{0\}\times {X}\subset \Cn\times {X}$ such that
  \begin{enumerate}
    \item[{\rm (1)}] $z\mapsto b(z,\omega)$ is holomorphic on $V_{\omega} = \{ z\in\Cn \colon (z,\omega) \in V \}$ for each $\omega\in X$;
    \item[{\rm (2)}] For each compact set $K\subset\Cn$ and for every point $\omega\in X$ with $K\times\{\omega\}\subset V$ there exists
      an open neighbourhood $U$ of $\omega$ such that $b$ is defined and uniformely bounded on $K\times U$ $($this is satisfied if e.g. every
      point in $X$ has a compact neighbourhood basis$)$.
  \end{enumerate}
   Then there exists a continuous map $\Gamma=\Gamma(z,\genmat,\omega)\colon \Cn\times\glnc\times {X}\to \Cn$,
defined on an open neighbourhood $\Omega$ of $\{0\}\times
\glnc\times {X}$ $($and holomorphic in the first two factors$)$
  such that if $u:(\Cn,0)\to(\Cn,0)$ is a germ of a biholomorphism satisfying
  $A(u(z)) = b(z,\omega_{0})$ for some $\omega_{0}\in {X}$, then necessarily
  $u(z) = \Gamma(z,u^\prime(0),\omega_{0})$. Furthermore, the map $\Gamma$
  has the properties {\rm (i)} and {\rm (ii)} given in Corollary~{\rm \ref{c:corparcn1}}.
\end{thm}

\begin{proof}
  We define the neighbourhood $\Omega$ in the following way:
  For each point $\omega\in X$ we choose subsets of the form
  $ \Delta_R\times U$, where $U$ and $R$ are small enough so that
  for $\omega\in U$, $b(\cdot,\omega)\in {\mathcal V}_R$ (with
  the notation of \S \ref{s:holomorphicity}). We claim that this is possible by
  Assumption (2): indeed, start by finding a neighbourhood $U_0$ of $\omega$ and an $R_{0}>0$ such that
  such that $b$ is defined on $\Delta_{R_0}\times U_0$. Choose an $R_1 < R_0$ and consider the compact set
  $\partial \Delta_{R_1} \subset \Cn$.
  By assumption, there exists a neighbourhood $U_1\subset U_0$ of $\omega$
  such that $\vnorm{b(\cdot,\omega)}_{\partial \Delta_{R_1}}\leq C$ for some constant $C$, for all
  $\omega \in U_1$. The Cauchy estimates imply that
  \[ \vnorm{b_\alpha (0,\omega)} \leq \frac{C}{ R_{1}^{|\alpha|}}, \quad \omega\in U_{1}.\]
  From these estimates it is clear that if we choose $0<R<R_1$ small enough,
  the desired claim follows.

  Next, for every relatively compact open set $W\subset \glnc$
  we may choose by Lemma \ref{l:addendum} an $\widetilde R>0$ such that
  $\Psi(b,\lambda)\in {\holmapsn}_{\widetilde R}$ for every $(b,\lambda)\in \mathcal{V}_R \times W$. We let
  $\Omega$ be the union of the open sets $\Delta_{\widetilde R}\times W\times U$, as $\omega$ ranges over $X$.
  The mapping $\Gamma$ is then defined by
  \[ \Gamma(z,\lambda,\omega) = \Psi(b(\cdot,\omega),\lambda) (z).\]
  We will check that $\Gamma$ is continuous by checking that it is continuous
  on sets of the form  $\Delta_{\widetilde R}\times W\times U$.
  To do that, we will show the following two facts: firstly, the mapping $\omega\mapsto b(\cdot,\omega)$ is continuous from $U$ to
  $\holmapsn_R$; and secondly, the evaluation mapping $(f,z)\mapsto f(z)$ is holomorphic from $\holmapsn_{\widetilde R}\times \Delta_{\widetilde R}$ to $\Cn$.
  These both follow from our assumptions using the Cauchy formula for polydiscs, as we will show; combining all these facts we get the result.

  Let us start with the proof of the first claim. Let $\omega\in X$ and $U(=U_\omega)$ be as in the above construction of $\Omega$.
  We check that $\omega\mapsto b(\cdot,\omega)$, as a map from $U$ to
  $\holmapsn_R$, is continuous at every point $\omega_0\in U$. Let $\epsilon > 0$ and  $0<r< R$.
  We need to show that there exists a neighbourhood $U_\epsilon$ of $\omega_{0}$ such that for
  $\omega\in U$,
  \[ \vnorm{ b(\cdot,\omega) - b(\cdot,\omega_{0}) }_r < \epsilon.\]
  We choose $s$ with $r < s < R$. By assumption (2), there exists a neighbourhood $V\subset U$ of $\omega_0$ such that $b$ is
  bounded on $\partial \Delta_s \times V$. The function
  \[ \rho (\omega) = \int_{\partial\Delta_s} \vnorm{ b(z,\omega) - b(z,\omega_0) } \, d|z| \]
  is continuous on $V$, and $\rho(\omega_0) = 0$. The Cauchy formula implies that
  \[ b(z,\omega) - b(z,\omega_0) =\frac{1}{(2\pi i)^n} \int_{\partial\Delta_s}
  \frac{ b(\zeta,\omega) - b(\zeta,\omega_0)}{\zeta - z} \, d\zeta, \quad z\in \Delta_r,\quad \omega \in V.\]
  Estimating, we get
  \[ \vnorm{ b(\cdot,\omega) - b(\cdot,\omega_0) }_{r} \leq \frac{1}{(2\pi)^n}(s - r)^{-n} \rho(\omega), \]
  and the desired estimate follows.

  The proof of the second claim is similar to what we did for proving the holomorphy of the
  solution mapping in \S \ref{s:holomorphicity}. It is clear that the evaluation mapping is holomorphic
  along finite dimensional affine subspaces. To complete the proof, we just have to show that the mapping is
  locally bounded. We leave it to the reader to check that this latter fact also follows from the Cauchy
  formula. The proof of the corollary is, therefore, complete.
\end{proof}

\section{Results for linear singular analytic systems}\label{s:linearcase}
 We consider here linear systems of the
form $\Theta(z)\cdot w=b(z)$, $w\in \C^m$, where $\Theta (z)$ is a
$m\times m$ matrix with holomorphic coefficients near the origin
in $\Cn$, of generic rank $m$, and $b\colon (\Cn,0)\to \C^m$ is a
germ of a holomorphic map. We are interested in getting a
holomorphic parametrization of all solutions $w=w(z)$ of such
systems, as well as other systems that are variations of the
original one. In fact, for such systems and contrary to the
nonlinear case, there is a unique potential obvious holomorphic
solution of the system; such a solution has to be the quotient of
the Weierstrass division of $\widehat \Theta\cdot b$ by $\det \Theta$,
where $\widehat \Theta$ denotes the transpose of the comatrix of $\Theta$ as in the proof of 
Lemma~\ref{l:multmap}. By analyzing the proof of the Weierstrass
division theorem, it is very likely that such a quotient would be
a suitable parametrization of the system satisfying the required
conditions of Proposition~\ref{t:lineareqn} below. 

Since we were not able to find in the literature the precise result given by Proposition~\ref{t:lineareqn} 
below, we choose instead to follow the approach 
developed in the previous section (for the non-linear case)  to get the parametrization needed for
our purposes here. The obtained parametrization may be a priori different from the
one given by Weierstrass division, since the algorithm proposed here 
to compute the parametrization is different from the algorithm given
by Weierstrass division. In what follows we keep the notation
defined in previous sections.

We start with the following:

\begin{prop}
  \label{t:lineareqn} Let $\Theta$ be a $m\times m$ matrix with
  holomorphic coefficients near the origin in $\Cn$, $m,n\geq 1$,
  such that $\Theta$ is of generic rank $m$. Then there exists $\ell_0\in
  \N$ and for every nonnegative integer $\ell$ polynomial mappings
  $q_\ell\colon \Cn\times J_0^{\ell_0+\ell}(\Cn,\C^m)\to \C^m$ homogeneous of degree $\ell$
  in their first variable such
  that for every germ of a holomorphic map $b\colon (\Cn,0)\to
  \Cm$,
  if $u\colon(\Cn,0)\to\Cm$ is a germ of holomorphic map satisfying $ \Theta(z)\cdot u(z) = b(z)$,
  then  for every $\ell\geq 0$,
  $u_{\ell} (z) = q_\ell (z, \jetm{0}{\ell_0 + \ell} b)$, where $u_\ell$ denotes the homogeneous
  part of order $\ell$ of $u$.  Furthermore, the $q_{\ell}$ can be chosen in such a way that
  there exists a polydisc $\Delta\subset \Cn$ centered at the origin for which the following holds:
  For all $C>0$ there exists $K>0$ with the property that
  whenever $b\colon (\Cn,0)\to \Cm$ is germ of a holomorphic map
  written into homogeneous terms $b(z)=\sum_{k\geq 0}b_k(z)$ and
  satisfying $\deltanorm{b_k}\leq C^k$ for all $k\in \N$, then
 $\deltanorm{q_k(\cdot,j_0^{\ell_0+\ell}b)}\leq K^k$ for all $k\in \N$.
\end{prop}

The proof of Proposition \ref{t:lineareqn} following our algorithm
is given in \S \ref{ss:lineareqn}. For our purposes, in addition to
the preceding result,
we need also the following refined version 
as well as a version with parameters given by
Proposition~\ref{t:lineqn2} below which follows immediately from Lemma~\ref{l:normsconv} and 
Proposition~\ref{t:lineqn}.

\begin{prop}
  \label{t:lineqn}
  Let $\Theta$ be a $m\times m$ matrix with
  holomorphic coefficients near the origin in $\Cn$, $m,n\geq 1$,
  such that $\Theta$ is of generic rank $m$. Then there exists $\ell_0\in
  \N$ and for every $\ell \in \N$ nonnegative integers $k(\ell)$ and polynomial mappings
  $p_\ell\colon \Cn\times J_0^{\ell_0+\ell}(\Cn)\times J_0^{\ell_0+\ell}(\Cn,\C^m)\to \C^m$ homogeneous
  of degree $\ell$ in their first variable such
  that for every germ of a holomorphic map $b\colon (\Cn,0)\to
  \Cm$ and every germ of a biholomorphic map $c\colon (\Cn,0)\to
  (\Cn,0)$,  if $u\colon(\Cn,0)\to\Cm$ is a germ of holomorphic map satisfying $ \Theta(c(z))\cdot u(z) = b(z)$,
  then  for every $\ell\geq 0$,
  $$u_{\ell} (z) = \frac{p_\ell (z, j_0^{\ell_0+\ell}c,\jetm{0}{\ell_0 + \ell} b)}{(\det j_0^1c)^{k(\ell)}},$$
  where $u_\ell$ denotes the homogeneous
  part of order $\ell$ of $u$.  Furthermore, the $p_{\ell}$ can be chosen in such a way that
  there exists a polydisc $\Delta\subset \Cn$ centered at the origin for which the following holds:
  For all $C>0$ and any compact subset $L\subset\glnc$ there exists $D>0$ with the property that
  whenever $b\colon (\Cn,0)\to \Cm$ is germ of a holomorphic map
  and $c\colon (\Cn,0)\to (\Cn,0)$ is a germ of a bihomorphism
  $($both written into homogeneous terms $b(z)=\sum_{k\geq 0}b_k(z)$, $c(z)=\sum_{k\geq
  0}c_k(z))$  satisfying $\deltanorm{b_k}\leq C^k$, $\deltanorm{c_k}\leq C^k$, for all $k\in
  \N$, and $c_1\in L$, then for every positive integer $\ell$, $$\deltanorm{\frac{p_\ell(\cdot,j_0^{\ell_0+\ell}c,j_0^{\ell_0+\ell}b)}
  {(\det j_0^1c)^{k(\ell)}}}\leq D^\ell.$$
\end{prop}

\begin{prop}
  \label{t:lineqn2}
  Let $\Theta$ be a $m\times m$ matrix with
  holomorphic coefficients near the origin in $\Cn$, $m,n\geq 1$,
  such that $\Theta$ is of generic rank $m$. Let $X$ be a complex manifold, and assume that
  $c\colon \C^n\times X \to \C^n$ and $b\colon \Cn\times X \to\Cm$ are  holomorphic maps defined on some
  neighbourhood of $\{ 0 \}\times X$. Assume that $c$ satisfies
  \[ c(0,\omega) = 0, \quad \det c_z (0,\omega) \neq 0, \quad {\rm for}\ {\rm every}\ \omega\in X. \]
  Then there exists a holomorphic map $\Gamma\colon \Cn\times X\to \Cm$ defined on a neighbourhood of $\{0\} \times X$
  such that if $u\colon (\Cn,0)\to\Cm$ is a germ of a holomorphic map satisfying $ \Theta(c(z,\omega_0)) \cdot u(z) = b(z,\omega_0)$
  for some $\omega_0\in X$, then
  $u(z)= \Gamma(z,\omega_{0})$. Furthermore, if we write the Taylor expansions
  $$ b(z,\omega)=\sum_{\alpha \in \N^n}b_{\alpha}(\omega)z^\alpha,\
  c(z,\omega)=\sum_{\alpha \in \N^n}c_{\alpha}(\omega)z^\alpha,$$
  then there exists an integer
  $\ell_0$ and for every $\alpha \in \N^n$ nonnegative integers $k_\alpha$ and polynomial mappings
  $p_\alpha \colon J_0^{\ell_0+|\alpha|}(\Cn, \Cm) \times J_0^{\ell_0+|\alpha|}(\Cn)\to \C^n$
   such that
  \begin{equation}
    \Gamma (z,\omega) = \sum_{\alpha \in \N^n}\frac{p_\alpha ((b_\beta(\omega),c_{\beta}(\omega))_{|\beta|\leq \ell_0+|\alpha|})}
    {(\det c_z(0,\omega))^{k_\alpha}}\, z^{\alpha}.
\end{equation}
\end{prop}

The proof of Proposition \ref{t:lineqn} is given in \S
\ref{ss:lineqn}.

\subsection{Proof of Proposition \ref{t:lineareqn}}\label{ss:lineareqn}

\subsubsection{Reducing the linear singular system to a simpler
one}\label{ss:reducelinear}
 A process similar to that used in \S \ref{ss:reduce} is done
 to bring the original $m\times m$ matrix $\Theta$ to another which 
 has a convenient form.
Given our $m\times m$ matrix $\Theta(z)$ with holomorphic
coefficients near the origin in $\Cn$, for each $j=1,\ldots,m$ we
write $\Theta^j (z)$ for its $j$-th row and decompose the rows
into homogeneous terms as $\Theta^j (z) = \sum_{k_j \leq k}
\Theta^j_k (z)$ where $\Theta^j_k$ is a homogeneous polynomial map
$\Cn\to \Cm$ of degree $k$. For $\nu \in \N$, we define
$\Theta_\nu$ to be the $m\times m$ matrix whose $j$-th row is
given by $\Theta^j_{k_j + \nu}$. Hence the $j$-th row of
$\Theta_\nu$ is a homogeneous polynomial map of degree $k_j+\nu$.
We say that $\Theta_0$ is the matrix formed from $\Theta$ by
taking its lowest order homogeneous terms. We have the following
reduction, analogous to that done in Lemma \ref{l:simplify1}.

\begin{lem}
  \label{l:simplify2} Let $\Theta(z)$ be a $m\times m$ matrix with holomorphic coefficients near the origin in
  $\Cn$ and assume that $\Theta$ is generically of full rank. Then there exists a $m\times m$ matrix $Q(z)$ of
  polynomials such that $N(z): = Q(z)\cdot \Theta(z)$ has the property that the matrix $N_0$ formed from $N$ by taking the
  lowest order homogeneous terms in each row is
generically of full rank.
\end{lem}
\begin{proof}
  Let $e$ be the order of vanishing (at the origin) of $\det \Theta$ and define $D := e - \sum_{1\leq j \leq m} k_j$.
  First note that $\det \Theta_0 $ is nonzero if and only if $D = 0$. In case $D=0$, we are therefore done. Hence we may
  assume that $D> 0$. We claim that there is
  an $m\times m$ matrix $Q^{(1)} (z)$ of polynomials such that
  if we denote $\tilde k_j$ the lowest order homogeneous in the $j$th row of the matrix
  $\tilde \Theta (z) := Q^{(1)} (z) \Theta(z)$ and if we denote by $\tilde e$ the order of vanishing of
  $\det \tilde \Theta$, then $\tilde D:=\tilde e-\sum_{1\leq j\leq m}$ is strictly smaller than $D$.
  Indeed, if $\det \Theta_0 = 0$, there exist polynomials $p_1,\dots,p_m\in \C [z]$ not all identically zero such
  that
  \begin{equation}\label{e:moredetails}
  p_1 \Theta^1_{k_1 } + \dots + p_m \Theta^m_{k_m} = 0.
  \end{equation}
   Moreover, since each $\Theta^j_{k_j}$ is homogeneous
  of degree $k_j$, we may choose
  each polynomial $p_j$ to be homogeneous of degree $\delta - k_j$ for some positive integer $\delta$. Without loss of
  generality, we may also assume that that $p_1 \neq 0$. Define
  \[Q^{(1)}:=\begin{pmatrix}
    p_1 & p_2 &  \dots & p_m \\
    0   & 1   &  \dots & 0   \\
        &     &  \ddots &    \\
    0   & \dots&       & 1
  \end{pmatrix},\]
  so that the order of vanishing of $\det \tilde \Theta $ is $e + \delta - k_1$ since $\det \tilde \Theta =p_1\, \det \Theta$.
  Note that for $j=2,\ldots,m$, $\tilde k_j=k+j$ and that it follows from \eqref{e:moredetails} that the lowest
  order homogeneous term of $\tilde \Theta^1$ is of
  order $\tilde k_1>\delta$. Hence $\tilde D = e + \delta - k_1 - \sum_{1\leq j\leq m} \tilde k_j < e - \sum_{1\leq j\leq m}
  k_j = D$. This proves the claim. Therefore,
  after a finite number of these procedures, 
  we arrive at a matrix $N$ satisfying the property claimed in the lemma.
\end{proof}

\subsubsection{Completion of the proof of Proposition~{\rm \ref{t:lineareqn}}} Let $\Theta_0$ be as defined  in
\ref{ss:reducelinear}. By Lemma~\ref{l:simplify2}, there exists a
$m\times m$ matrix $Q(z)$ of polynomials such that the matrix
$N_0(z)$ formed from $N(z):=Q(z)\cdot \Theta (z)$ by taking the
lowest order homogeneous terms in each row is generically of full
rank. Now observe that if we prove Proposition \ref{t:lineareqn}
for the matrix $N$ the result easily follows for the original
matrix $\Theta$. So we may assume in what follows that the matrix
$\Theta_0$ is generically of full rank. Furthermore, by
multiplying the rows with appropriate homogeneous polynomials, we
may assume that the rows of $\Theta_0$ are homogeneous of the same
degree, say $\ell_0\geq 0$. By Lemma \ref{l:multmap}, we may
choose a polydisc $\Delta$ centered at the origin and a constant
$D>0$ and for all $\ell\geq 0$ a left inverse $B_\ell$ to the
operator given by multiplication with $\Theta_0$ (which takes
$\poly^m_\ell$ to $\poly^m_{\ell_0 + \ell}$), and which has norm
at most $D^{-1}$ (when $\poly^m_\ell$ and $\poly^m_{\ell_0 +
\ell}$ are equipped with $\deltanorm{\cdot}$). We define
inductively the polynomials maps $(q_\ell)_{\ell \geq 0}$ as
follows:
\begin{equation}
\begin{cases}
q_0(\Lambda^{\ell_0}):=B_0\left(\sum_{|\alpha|=\ell_0}\displaystyle
\frac{\Lambda_{\alpha}}{\alpha!}z^{\alpha}\right),\cr
q_{\ell}(z,\Lambda^{\ell_0+\ell}) := B_\ell \left(
\sum_{|\alpha|=\ell_0+\ell}\displaystyle
\frac{\Lambda_\alpha}{\alpha!}z^{\alpha} - \sum_{j=0}^{\ell-1}
\Theta_{\ell-j}(z)\cdot q_j(z,\Lambda^{\ell_0+j}) \right),\ \ell
\geq 1.
  \label{e:pelldeflin}
  \end{cases}
\end{equation}
It is easy to show by induction on $\ell$ that each $p_\ell\colon
\C^n \times \njetsp{0}{\ell_0 + \ell}(\Cn,\C^m)\to \C^m$ is a
polynomial map that is homogeneous of degree $\ell$ in its first
variable $z$. We claim that given a germ of a holomorphic map
$b\colon (\Cn,0)\to \Cm$ if $u\colon (\Cn,0)\to \C^m$ is a germ of
a holomorphic map satisfying
\begin{equation}\label{e:stuff}
\Theta (z)\cdot u(z)=b(z),
\end{equation}
then necessarily $u_\ell (z)=q_\ell (z,j_0^{\ell_0+\ell}b)$.
Indeed first note that by considering homogeneous terms of order
$\ell_0$ in \eqref{e:stuff}, we get $\Theta_0(z)\cdot
u_0(z)=b_{\ell_0}(z)$ and hence $u_0(z)=q_0(j_0^{\ell_0}b)$. Now
assume that $u_k(z)=q_k(z,j_0^{\ell_0+k}b)$ for all $k<\ell$.
Considering homogeneous terms of order $\ell+\ell_0$ in
\eqref{e:stuff}, we obtain
$$\Theta_0(z)\cdot u_\ell (z)+\sum_{\nu=1}^{\ell}\Theta_{\nu}(z)\cdot u_{\ell -\nu}(z)=b_{\ell_0+\ell}(z),$$
which easily gives the required conclusion by the induction
assumption.

 The convergence statement can now be derived by a
similar but much simpler argument as in the proof of
Theorem~\ref{t:parcn1}. Indeed, let $C>0$ and let $b\colon
(\Cn,0)\to \Cm$ be a germ of holomorphic map satisfying
$\deltanorm{b_k}\leq C^k$ for all $k\in \N$ where
$b(z)=\sum_{k\geq 0}b_k(z)$ is the decomposition of $b$ into
homogeneous terms; without loss of generality we may assume that
$\ell_0 > 0$. Since $\Theta$ is convergent, there exists $a>0$
such that $\deltanorm{\Theta_j}\leq a^j$. We first claim that for
any $\sigma_0\geq 0$, the series $\varphi(t): = \sum_{j\in \N}
\sigma_j t^j$ converges in a neighbourhood of $0$, where
$\sigma_j$ is defined for $j\geq 1$ by
\begin{equation}\label{e:i}
D\sigma_j = \sum_{k<j} a^{j-k} \sigma_k + C^{j+\ell_0}.
\end{equation}
This follows from the fact that $\varphi (t)$ is a solution of the
equation
\[ D \varphi(t) - D\sigma_0  = \frac{at}{1-at} \varphi(t) + \frac{tC^{\ell_0+1}}{1-Ct}.\]
Now, we set $\sigma_0:=C^{\ell_0}D^{-1}$ and $(\sigma_j)_{j\geq
1}$ as given by \eqref{e:i}. Then by using \eqref{e:pelldeflin}
and \eqref{e:i} it is easy to show by induction that for all $\ell
\in \N$, $\deltanorm{q_\ell(\cdot,j_0^{\ell_0+\ell}b)}\leq
\sigma_\ell$. Since $\sum_{j\geq 0}\sigma_jt^j$ is convergent,
there exists $K>0$ (independent of $b$) such that $\sigma_\ell
\leq K^\ell$ for all $\ell \in \N$. This completes the proof of
Proposition \ref{t:lineareqn}.

\subsection{Proof of Proposition~\ref{t:lineqn}}\label{ss:lineqn}
Actually, this proposition is a simple consequence of Proposition
\ref{t:lineareqn} and the inverse function Theorem. We will, of
course, need a version of the latter which is compatible with the
desired estimates. So we first turn to this. The reader will
notice that the results we derive here are by no means new. 
Nevertheless we think it is important to put the proofs of the
precise desired results here. We first have the following lemma.

\begin{lem}\label{l:invfun}
  For every positive integer $\ell$, there exist a polynomial map
  $\rho_\ell \colon \Cn\times\njetsp{0}{\ell} (\Cn)\to
  \Cn$, homogeneous of degree $\ell$ in its first variable,
  such that for any germ of a biholomorphic map $c\colon (\Cn,0)\to
  (\Cn,0)$ written into homogeneous terms $c(z)=\sum_{j\geq
  1}c_j(z)$ is given by
  \[ c^{-1} (z) = \sum_{\ell >0 }
  \frac{\rho_{\ell}(z,\jetm{0}{\ell} c)}{(\det c_1)^{\ell}}.
  \]
  Furthermore, there exists a polydisc $\Delta \subset \Cn$ centered at the origin such that
  for any $C>0$ and any compact subset $L\subset \glnc$
  there exists $K>0$ such that for every $c\colon (\Cn,0)\to (\Cn,0)$ as above satisfying $c_1\in L$ and $\deltanorm{c_k}
  \leq C^k$ for every $k\in \N^*$, then for every $\ell \in \N^*$
  $$\deltanorm{\frac{\rho_{\ell}( \cdot , \jetm{0}{\ell} c )}{(\det c_1)^\ell}} \leq
  K^\ell.$$
\end{lem}
\begin{proof}
We set $\widehat\rho_1(z,\Lambda^1):=(\Lambda^1)^{-1}\cdot z$ and
define the polynomials $\rho_\ell$, $\ell \geq 2$, inductively as
follows (here we use again the notation $u_\ell (z) = u(z)|_\ell$):
\begin{equation}\label{e:fedup}
\begin{cases}
 \widehat\rho_\ell (z, \Lambda^\ell): = - (\Lambda^1)^{-1}
  \left( \sum_{1 < j\leq d} \sum_{|\alpha|=j}\displaystyle \frac{\Lambda_\alpha}{\alpha!} \left((R^{\ell-1} (z,\Lambda^{\ell -1}))^\alpha\right)|_{\ell}
  \right),\cr
  R^{\ell -1}(z,\Lambda^{\ell -1}):=\sum_{1\leq \nu \leq \ell -1}\widehat\rho_\nu
  (z,\Lambda^{|\nu|}).
  \end{cases}
 \end{equation}
It is easy to see that for every $\ell \geq 1$,
$\rho_\ell(z,\Lambda^\ell):=(\det c_1)^\ell\, \widehat \rho_\ell
(z,\Lambda^\ell)$ satisfies the required conditions of Lemma
\ref{l:invfun}. It remains therefore to check the estimates for
the $\widehat \rho_\ell$. So let $\Delta$ be the unit polydisc in
$\Cn$, $L$ a compact subset of $\glnc$, $C>0$ and suppose that
$c\colon (\Cn,0)\to (\Cn,0)$ is a germ of a biholomorphic map
satisfying $\deltanorm{c_k}\leq C^k$ for all $k>0$. We set
  \[ A := \min_{\genmat \in L} \frac{1}{\vnorm{\genmat^{-1}}}, \quad B_1 := \max_{\genmat\in L} \vnorm{\genmat^{-1}}. \]
(Here, as in Lemma \ref{l:multmap}, the matrix norm is adapted to
the maximum norm on $\Cn$.) Let $\widetilde C>0$ (independent of
$c$) such that $\onenorm{c_j}\leq \widetilde C^j$ for all $j>0$.
Consider the one-dimensional equation
  \begin{equation}\label{e:unique}
   A \phi(t) - \frac{\widetilde C^2 \phi (t)^2}{1 - \widetilde C \phi (t)} = A B_1 t,
   \end{equation}
  It is then easy to check that \eqref{e:unique} has a unique
  convergent solution
  \[ \phi (t) = \sum_{\ell>0} B_{\ell}\,  t^{\ell} \]
  vanishing at the origin, which furthermore satisfies $B_\ell>0$
  for all $\ell>0$. Moreover, Proposition~\ref{p:faadebruno} implies that
  the $B_{\ell}$ satisfy the following relation for ${\ell}\geq 2$:
  \begin{equation}\label{e:relation}
   A B_{\ell} = \sum  \frac{d! C^d} {k_1! \dots k_{{\ell}-1}!} \, B_{1}^{k_1} \dots B_{{\ell}-1}^{k_{{\ell}-1}},
   \end{equation}
  where $d=k_1 + \dots + k_{{\ell}-1}\geq 2$, and the sum goes over all nonnegative integers
  $k_1,\dots,k_{{\ell}-1}$ with $k_1 + 2 k_2 + \dots + ({\ell}-1) k_{{\ell}-1}=\ell$.
  We claim that for all $\ell\geq 1$, $\deltanorm{\widehat \rho_\ell(\cdot, \jetm{0}{\ell} c)} \leq B_j$.
  For $\ell=1$ this is
  clear by the choice of $B_1$ and the definition of $\widehat \rho_1$.
  Suppose now that $\deltanorm{\widehat \rho_d(\cdot,  \jetm{0}{d} c)} \leq B_d$ for
  $d < \ell$. Using \eqref{e:fedup} and Lemma~\ref{l:homcompest},
  we get the following estimates
  \begin{align*}
    \deltanorm{\widehat \rho_{\ell} (\cdot, \jetm{0}{\ell} c )} &\leq \vnorm{c_{1}^{-1}} \deltanorm{ \sum_{1 < j\leq \ell}  \left(c_j
   (R^{\ell-1} (z,j_0^{\ell -1}c))\right)|_{\ell} } \\
    &\leq \vnorm{c_{1}^{-1}} \sum_{2\leq j\leq \ell}
    \frac{j!}{k_{1} \dots k_{\ell - 1}!} \onenorm{c_j} \deltanorm{\widehat \rho_{1}}^{k_{1}} \dots \deltanorm{\widehat \rho_{\ell-1}}^{k_{\ell-1}}
    \end{align*}
    where $j=k_1 + \dots + k_{{\ell}-1}$, and the sum goes over all nonnegative integers
  $k_1,\dots,k_{{\ell}-1}$ with $k_1 + 2 k_2 + \dots + ({\ell}-1)
  k_{{\ell}-1}=\ell$. Hence, by \eqref{e:relation} and the induction assumption we get
    \begin{align*}
    \deltanorm{\widehat \rho_{\ell} (\cdot, \jetm{0}{\ell} c )} &\leq  \vnorm{c_{1}^{-1}} \sum_{2\leq j\leq \ell}
    \frac{j!}{k_{1} \dots k_{\ell - 1}!} C^d B_{1}^{k_{1}} \dots B_{\ell-1}^{k_{\ell-1}} \\
    & =   \vnorm{c_{1}^{-1}} A B_\ell \leq B_\ell,
  \end{align*}
  by our choice of $A$, which proves the claim and also completes the proof of the lemma.
\end{proof}

\begin{rem}\label{r:1}
Note that Lemma \ref{l:invfun} also holds for {\em every} polydisc
$\Delta \subset \Cn$ centered at the origin.
\end{rem}
We also need the following lemma about the composition of analytic
functions with estimates.

\begin{lem}\label{l:compan}
  Let $K>0$ and assume that $\Delta\subset\CN$ and $\Delta^{\prime}\subset\Cn $
  are polydiscs centered at $0$. There exists $C>0$ such that whenever
  $f\colon (\Cn,0)\to \Cm$ and $g\colon (\CN,0)\to (\Cn,0)$ are
  germs of holomorphic mappings written into homogeneous terms $f=\sum_{j\geq
  0}f_j$, $g=\sum_{j\geq 1}g_j$ and  satisfying $\deltanorm{g_j}\leq K^{j}$  and $\vnorm{f_j}_{\Delta^\prime}\leq K^j$
  for every positive integer $j$, then
  the composition $h= f\circ g=\sum_{j\geq 0}h_j$ satisfies
   $\vnorm{h_j}_{\Delta}\leq C^j$ for $j>0$.
\end{lem}
Before we begin with the proof, let us point out that there are easier proofs for 
this lemma than the one we give; we prefer to give a proof 
using majorization by a conveniently chosen
power series because this ties in nicely with the proofs of the other convergence
estimates given before.
\begin{proof}
    We use the geometric sum to compute
  \begin{equation}\label{e:exam}
   \frac{K_1 \frac{K t}{1- K t}}{1-K_1 \frac{K t}{1-K t}} = \frac{K_1 K t}{1- K(1+ K_1) t} = \sum_{j>0} K_1 K^j (1 + K_1)^{j-1} t^j.
   \end{equation}
  Now choose $K_1$ large enough (depending only on $K$ and $\Delta'$) so that $\onenorm{f_j}  \leq K_1^j$ for all $j>0$,
  for all $f$ satisfying $\vnorm{f_j}_{\Delta^{\prime}} \leq K^j$.
  We claim that $\deltanorm{h_j} \leq  K_1 K^j(1+ K_1)^{j-1}$ for $j>0$. This is proved using Lemma~\ref{l:homcompest} and
  Faa di Bruno's formula (Proposition~\ref{p:faadebruno}). For every $j>0$, we have
  the following estimates:
  \begin{align*}
    \deltanorm{h_j} &\leq \sum_{d=1}^{j} \deltanorm{(f_{d}\circ g)|_{j}} \\
    &\leq \sum \frac{d! \, \onenorm{f_d}}{k_1!\dots k_j!} \, \deltanorm{g_{1}}^{k_1}\dots
    \deltanorm{g_{j}}^{k_j},
    \end{align*}
  where $d=k_1+\ldots+k_j\geq 1$ and the sum goes over all nonnegative
  integers $k_1,\ldots,k_{j}$ such that $k_1+2k_2+\ldots+jk_j=j$.
  Using \eqref{e:exam}, we get
    \begin{align*}
        \deltanorm{h_j} &\leq \sum  \frac{d! \, K_1^d}{k_1!\dots k_j!} \, K^{{k_1}+ 2 k_2 + \dots j{k_j}} \\
    & = K_1 K^j (1+K_1)^{j-1}.
  \end{align*}
 Setting $C :=  K(1+K_1)$ yields the desired result.
\end{proof}

\begin{proof}[Completion of the proof of Proposition~{\rm \ref{t:lineqn}}]
We can now prove Proposition~\ref{t:lineqn}. Let $\ell_0$
given by Proposition \ref{t:lineqn} and for every $k \in \N$, let
$q_k$ be the polynomial map given by the same proposition. Define
for every nonnegative integer $\ell $ the map $\tilde p_\ell$ by
setting
\begin{equation}\label{e:soccer}
 \tilde{p}_\ell (z, \jetm{0}{\ell + \ell_0} c,
\jetm{0}{\ell+\ell_0} b) = \left(\sum_{k\in \N}q_{k} (c(z),
\jetm{0}{\ell + \ell_0} (b\circ c^{-1}))\right)|_{\ell},
\end{equation}
for all germs of holomorphic maps $b$ and $c$ as in Proposition
\ref{t:lineqn}, where, as used before, the notation $|_\ell$ on the right hand side of \eqref{e:soccer} denotes the homogeneous part of order $\ell$ of the corresponding power series mapping. It is clear from Lemma~\ref{l:invfun} and the
chain rule that for every $\ell \in \N$, there exists $k(\ell)\in
\N$ such that for every $b$ and $c$ as above we have
\[{(\det \jetm{0}{1}c)^{k(\ell)}}\, \tilde{p}_\ell (z, \jetm{0}{\ell + \ell_0} c, \jetm{0}{\ell+\ell_0} b)
=: {{p}_\ell (z, \jetm{0}{\ell + \ell_0} c, \jetm{0}{\ell+\ell_0}
b)} \] for some polynomial map satisfying the required properties.
Finally, we leave the reader to check that the needed estimates
follow from the estimates given in Lemma~\ref{l:invfun} and Lemma
\ref{l:compan} and from Remark \ref{r:1}. The proof of Proposition
\ref{t:lineqn} is complete.
\end{proof}

\begin{rem}
Let ${\mathcal H}(\C^n,\C^m)$ be the topological vector space of germs of holomorphic maps $(\C^n,0)\to \C^m$ (endowed with the topology of uniform convergence on compact neighbourhoods of the origin in $\C^n$) and $\Theta =\Theta (z)$ be a $m\times m$ matrix with entries in ${\mathcal H}(\C^n,\C)$ with $\Theta$ of generic rank $m$. Define the holomorphic  map $\Theta_{*}\colon \biholmaps \times {\mathcal H}(\C^n,\C^m)\to \biholmaps \times {\mathcal H}(\C^n,\C^m)$ by setting $\Theta_{*}(c,u)=(c,\left(\Theta \circ c\right) \cdot u)$ for all $u\in {\mathcal H}(\C^n,\C^m)$ and $c\in \biholmaps$. Then we leave to the reader to check that, following the lines of the proof of Theorem~\ref{t:parcnfull} given in \S \ref{s:holomorphicity}, the statement given by Proposition~\ref{t:lineqn} yields the existence of a holomorphic left inverse $\Phi\colon \biholmaps \times {\mathcal H}(\C^n,\C^m) \to \biholmaps \times {\mathcal H}(\C^n,\C^m)$ to the map $\Theta_{*}$. 
\end{rem}

%\part{Application to the parametrization of groups of automorphisms of real-analytic CR manifolds}

\section{A class of real-analytic generic submanifolds}\label{s:class}
We introduce here a class of real-analytic generic submanifolds,
more general than the class of essentially finite generic
submanifolds, and state the most general parametrization theorem
of this paper (Theorem \ref{t:further1}) for such a class; all the
results stated in \S \ref{s:highcodim} will then follow from that
theorem. We also discuss a few properties of this class of
manifolds, compare it to some well-known other ones and also give
several examples of manifolds in this class that are not essentially
finite.

\subsection{Nondegeneracy conditions for real-analytic generic
submanifolds}\label{ss:nondegeneracy} We start by introducing a
nondegeneracy condition for real-analytic CR submanifolds of $\CN$
and then compare it to some well-known ones. For more background
concerning CR structures we refer the reader e.g.\ to the books
\cite{BERbook, boggess}.

Let $(M,p)$ be a germ of a real-analytic generic submanifold of $\CN$ 
 i.e.\ satisfying $T_pM+J(T_pM)=T_p\CN$ where $J$ is the complex
structure map of $\CN$ and $T_pM$ (resp.\ $T_p\CN$) denotes the
tangent space of $M$ (resp.\ of $\CN$) at $p$. Let $n$ and $d$ be
respectively the CR dimension and the codimension of the manifold
$M$ so that $N=n+d$. Let $\rho =(\rho_1,\ldots,\rho_d)$ be a
real-analytic vector valued defining equation for $M$ in some
neighbourhood $U$ of $p$ in $\CN$ satisfying $\partial
\rho_1\wedge \ldots \wedge \partial \rho_d\not =0$ in $U$.  For
every point $q\in \CN$ sufficiently close to $p$, recall that the
Segre variety attached to $q$ is the $n$-dimensional complex
submanifold of $\CN$ given by $$S_q:=\{Z\in U:\rho
(Z,\overline{q})=0 \},$$ where the original real-analytic map
$\rho$ has been complexified and $U$ may have been shrinked so that
$\rho(Z,\zeta)$ is now convergent on $U\times \overline{U}$. 
On the other hand, the
complexification of $M$ is the germ through the point $(p,\bar p)\in \C^{2N}$ of the
$2n+d$-dimensional complex submanifold of $\C^{2N}$ given by
\begin{equation}\label{complexification}
{\mathcal M}:=\{(Z,\zeta)\in (\C^{2N},(p,\bar p)):\rho (Z,\zeta)=0\}.
\end{equation}
Here and throughout the paper, given any holomorphic map $\eta$ defined on some open subset $\Omega$ of $\C_t^k$
and given any point $a\in \Omega$, we denote by $\{t\in (\C^k,a): \eta (t)=0\}$ the germ at $a$ of the complex-analytic set $\{t\in \Omega: \eta (t)=0\}$.

For every integer $k$ and for $u\in \CN$, let $E^{k,n}_u(\CN)$ be
the jet space of order $k$ at $q$ of $n$-dimensional complex
submanifolds of $\CN$ passing through $u$. Then $\cup_{u\in
\CN}E_u^{k,n}(\CN)$ carries a natural fiber bundle structure over
$\CN$. Following \cite{DW1, BJT}, for every $q\in M$ sufficiently
close to $p$, we consider the 
antiholomorphic map $\pi_q^k$ defined as
follows:
\begin{equation}\label{e:segredef} \pi_q^k \colon S_q \to E^{k,n}_q(\CN), \quad \pi_q^k (\xi) = j^k_q
S_\xi,
\end{equation}
where $j^k_q S_\xi$ denotes the $k$-jet at $q$ of the submanifold
$S_\xi$ (see e.g.\ \cite{GG}).

\begin{defin}
{\em We say that the germ $(M,p)$ belongs to the class $\mcl$ if
the antiholomorphic map $\pi^k_p$ is of generic rank $n$ for $k$ large
enough $($or equivalently, is generically immersive for $k$ large
enough$)$}. \end{defin} For a germ $(M,p)\in \mcl$, we denote by
$\kappa_M(p)$ the smallest integer $k$ for which the map $\pi_p^k$
is of generic rank $n$. Clearly, the integer $\kappa_M(p)$ is a
biholomorphic invariant attached to the germ $(M,p)$. Furthermore, for $q\in M$
sufficiently close to $p$, we see that, shrinking the neighbourhood $U$ if necessary, that
 $S_q$ depends anti-holomorphically on $q$ (i.e. can be parametrized by a holomorphic map depending anti-holomorphically on $q$; see e.g.\ \S \ref{ss:propertiesmcl} below for more details if needed). By unique continuation it follows that if $(M,p)\in {\mathcal C}$, then 
the map $\pi^{\kappa_M(p)}_q$ is
still  generically immersive for 
$q\in M$ close enough to $p$ i.e.\ that $(M,q)$ is in $\mcl$. Hence,
there exists a neighbourhood $V$ of $p$ in $M$ such that
$\kappa_M$ is well defined on $V$ and from the definition we have
that $\kappa_M$ is {\em upper-semicontinuous}.

It is interesting to compare the above local open condition with
other well-known ones. We recall the following nondegeneracy
conditions for a germ $(M,p)$ of a real-analytic generic
submanifold of $\CN$:
\begin{enumerate}
\item[(I)] the submanifold $M$ is {\em Levi-nondegenerate} at $p$ if the map $\pi_p^1$ is an immersion at $p$;
\item[(II)]  the submanifold $M$ is said to be {\em finitely
nondegenerate} at $p$ if the  map $\pi_p^k$ is an
immersion at $p$ for $k$ sufficiently large; in this case, $M$ is said be
$k_p$-nondegenerate at $p$ if $k_p$ is the smallest $k$ for which
the above condition holds (see \cite{BHR1, BERbook}).
\item[(III)] the submanifold $M$ is said to be {\em essentially
finite} at $p$ if the antiholomorphic map $\pi_p^k$ is finite for $k$
large enough; in this case, the smallest of such $k$'s is called
the {\em essential type} of $M$ at $p$ and denoted ${\rm
esstype}_M(p)$ (see \cite{BJT}).
\item[(IV)] the submanifold $M$ is said to be {\em holomorphically nondegenerate}
at $p$ (in the sense of Stanton \cite{S2}) if for a generic point
$q\in M$ sufficiently close to $p$ the antiholomorphic map $\pi_q^k$
is generically immersive for $k$ large enough. (Equivalently, for
a generic point $q\in M$ sufficiently close to $p$, $(M,q)\in
\mcl$.)
\end{enumerate}
The definition of holomorphic nondegeneracy given in (IV) does not
correspond here to the original one but it is not difficult to
show by using Stanton's criterion of holomorphic nondegeneracy
(see Proposition~\ref{p:stanton} below) that it is in fact equivalent to the
original one. It is well-known that
(I)$\Rightarrow$(II)$\Rightarrow$(III)$\Rightarrow$(IV) (but
converses do not hold). It is also clear from the above
definitions that if $(M,p)\in \mcl$ then $(M,p)$ is
holomorphically nondegenerate and that if $(M,p)$ is essentially
finite then $(M,p)\in \mcl$. Hence, the essentially finite
submanifolds form one general subclass of
submanifolds that are in the class $\mcl$. Another subclass of the
class $\mcl$ that is worth to point out consists of those germs
$(M,p)$ of generic real-analytic submanifolds that are rigid and
holomorphically nondegenerate (see \S \ref{ss:propertiesmcl} for
the definition and details). We summarize the above in the
following:

\begin{prop}\label{t:prop1} Consider the following conditions on a germ $(M,p)$ of a real-analytic generic submanifold
of $\CN$:
  \begin{enumerate} \item[{\rm (i)}] $(M,p)$ is essentially finite;
    \item[{\rm (ii)}] $(M,p)$ does not contain any nontrivial complex analytic subvariety;
    \item[{\rm (iii)}] $(M,p)$ is rigid and holomorphically nondegenerate.
  \end{enumerate}
  If $(M,p)$ satisfies one of the above three conditions then
  $(M,p)\in\mcl$.
\end{prop}
In \S \ref{ss:propertiesmcl}, we shall provide the reader with
an elementary example of a germ of a submanifold belonging to
the class $\mcl$ that is not essentially finite nor rigid
holomorphically nondegenerate and also give a more explicit
criterion (Proposition \ref{p:classnormal}) that allows one to
decide when a given germ of a submanifold is in $\mcl$ or not.

Another nondegeneracy condition on a germ $(M,p)$ we also need to
recall is the minimality condition due to Tumanov \cite{TU2}:
$(M,p)$ is called minimal if there is no proper CR submanifold
$S\subset M$ passing through $p$ with the same CR dimension as
that of $M$. As is well-known (see \cite{BERbook}), for real-analytic generic submanifolds, 
this condition is equivalent to the finite type condition of Kohn \cite{Ko1} and
Bloom-Graham \cite{BG1} on the CR vector fields of $M$.

\subsection{The main result of the paper}\label{ss:maincrap}
The most general parametrization theorem we shall prove in this
paper is the following.

\begin{thm}
  \label{t:further1}  Let  $M$ be a real-analytic generic submanifold of $\CN$ of codimension $d$.
  Let $p\in M$ and assume that $(M,p)$ is minimal and belongs to
  the class $\mcl$ and set $\ell_p:=(d+1)\kappa_M(p)$.
  Then there exists an open subset $\Omega \subset \CN \times G_{p,p}^{\ell_p}(\CN)$ and a real-analytic map
  $\Psi (Z,\Lambda) \colon \Omega\to \CN$ holomorphic in the first factor, such that the following hold :
  \begin{enumerate}
  \item[{\rm (i)}] for any $H\in\autMp $ the point
  $(p,j_p^{\ell_p}H)$ belongs to  $\Omega$ and the following identity holds:
  \[H(Z) = \Psi(Z,\jetm{p}{\ell_p} H)\text{ for all } Z\in \CN \text{ near } p;\]
  \item[{\rm (ii)}] the map $\Psi$ has the following formal Taylor
  expansion
  \begin{equation}
  \Psi (Z,\Lambda)=\sum_{\alpha \in \N^N}\frac{P_\alpha (\Lambda,\overline{\Lambda})}
  {(D(\Lambda^1)^{s_\alpha}(\overline{D(\Lambda^1)})^{r_\alpha}}\, (Z-p)^{\alpha},
  \end{equation}
  where for every $\alpha \in \N^N$, $s_\alpha$ and $r_\alpha$ are nonnegative
  integers, $P_\alpha$ and $D$ are polynomials in their arguments, 
  $\Lambda^1$ denotes the linear part of the jet $\Lambda$ and $D(j_0^1H)\not =0$ for every $H\in \autMp$.
  \end{enumerate}\end{thm}

Let us remark here that the set $\Omega$ appearing in the statement of 
Theorem~\ref{t:further1} is only determined by 
some algebraic conditions on the first order
jets, and these conditions are automatically fulfilled for 
jets of CR automorphisms of $M$.

   Quite analogously to \S \ref{s:highcodim}, we obtain by following some arguments from \cite{BERbook, BER4} 
   the following result from Theorem \ref{t:further1}.

  \begin{thm}
   \label{t:further2} Let  $M$ be a real-analytic generic submanifold of $\CN$ of codimension $d$.
  Let $p\in M$ and assume that $(M,p)$ is minimal and belongs to
  the class $\mcl$ and set $\ell_p:=(d+1)\kappa_M(p)$. Then the jet mapping
  \begin{equation}\label{e:promote}
   \jetm{p}{\ell_p} \colon \autMp \to G_{p,p}^{\ell_p}(\CN)
   \end{equation}
  is a continuous group homomorphism that is a homeomorphism onto a real-algebraic
  Lie subgroup of $G_{p,p}^{\ell_p}(\CN)$.
\end{thm}

Another special class of real-analytic generic submanifolds of
$\CN$ that has been much studied in recent years is that of
finitely nondegenerate ones.  For a germ of a $k_p$-nondegenerate
real-analytic generic submanifold $(M,p)$, Theorem
\ref{t:further1} (and hence Theorem \ref{t:further2}) has been
established by Baouendi, Ebenfelt and Rothschild in \cite{BER4}
(see also \cite{Z2} for a previous weaker version) with a bound
$\ell_p'=(d+1)k_p$. We would like to point out that it follows
from the discussion in \S \ref{ss:nondegeneracy} that we always
have $\ell_p\leq \ell_p'$ with $\ell_p$ given by Theorem
\ref{t:further1}. Furthermore,  we give in \S
\ref{ss:propertiesmcl} a number of examples of finitely
nondegenerate submanifolds for which the above inequality is in
fact strict (see Example \ref{e:hitit}). This shows that Theorem
\ref{t:further1} improves also the known results for finitely
nondegenerate submanifolds.

The proof of Theorem \ref{t:further1} is given throughout \S
\ref{s:jets}--\S \ref{s:jetpar} and is completed in \S
\ref{s:conclusion} together with the proof of Theorem
\ref{t:further2}.

\subsection{Some properties of the class $\mcl$ and comparison between known invariants}\label{ss:propertiesmcl}
Our goal here is to give a characterization (Proposition
\ref{p:classnormal}) by using local holomorphic coordinates of
germs of real-analytic generic submanifolds that belong to $\mcl$.

Let $(M,p)$ be a germ of a real-analytic generic submanifold of
$\CN$. Following \cite{BERbook}, we use a system of {\em normal
coordinates} $(z,w)\in\Cn\times\Cd$ for $(M,p)$, which means that
in these coordinates $p=0$ and that there is a defining equation
for $M$ near the origin of the form
\begin{equation}
  \imag w = \varphi ( z, \bar{z},\real w ), \quad \varphi(z,0,s) = \varphi(0,\bar{z},s) = 0,
  \label{e:normalreal}
\end{equation}
where $\varphi \colon \Cn\times\Rd \to \Rd$ is a germ of a
real-analytic map at the origin. By the implicit function theorem,
\eqref{e:normalreal} is also equivalent to a defining equation of
the following form:
\begin{equation}
  w = Q ( z, \bar{z}, \bar{w} ),
  \label{e:normalcomplex}
\end{equation}
with $Q(z,0,\bar{w}) = Q(0,\bar{z},\bar{w}) = \bar{w}$; here
$Q\colon\Cn\times\Cn\times\Cd\to\Cd$ is a germ of a holomorphic
map at the origin. We write and use the following expansions:
\begin{equation}
  \bar{Q} (\chi,z,w) = \sum_{\alpha \in \N^n} \bar{q}_\alpha ( z,w)
  \chi^\alpha,\ \bar{Q}=(\bar{Q}^1,\ldots,\bar{Q}^d),\
  \bar{q}_\alpha=(\bar{q}^1_\alpha,\ldots,\bar{q}^d_\alpha).
  \label{e:normalexp}
\end{equation}
Here and throughout the rest of the paper, for any holomorphic
function $h\colon \C^r\to \C$ defined near the origin in $\C^r$,
we denote by $\overline{h}\colon \C^r\to \C$ the holomorphic
function obtained from $h$ by taking the complex conjugates of the
coefficients of its Taylor series. In the above normal
coordinates, the complexification ${\mathcal M}$ is given by
\begin{equation}\label{e:complexcoord}
{\mathcal M}:=\{((z,w),(\chi, \tau))\in (\CN\times \CN,0):\tau
=\bar{Q}(\chi,z,w)\},
\end{equation}
and for every integer $k$, the mapping $\pi_p^k$ defined in
\eqref{e:segredef} may be identified with the conjugate of the local holomorphic
map $z\mapsto (\alpha !\, \overline{q}_\alpha (z,0))_{|\alpha|\leq
k}$. As a consequence, we have the following characterization.

\begin{prop}
  \label{p:classnormal} Let $(M,p)$ be a germ of a real-analytic generic submanifold in $\CN$ and let $(z,w)$ be normal
  coordinates for $(M,p)$ so that $M$ is given by \eqref{e:normalcomplex}. Then the
  following are equivalent:
  \begin{enumerate} \item[{\rm (i)}] $(M,p)\in\mcl$;
    \item[{\rm (ii)}]\label{bush} for some integer $k$, the local holomorphic map $z \mapsto \left( \bar{q}_\alpha (z,0) \colon | \alpha | \leq k \right)$
    is generically of full rank.
  \end{enumerate}
  Furthermore, if any of these conditions is satisfied, the integer $\kappa_M(p)$
  is the smallest integer $k$ for which condition {\rm (ii)} holds.
\end{prop}

As already mentioned in \S \ref{ss:nondegeneracy}, the first
natural class of manifolds that belong to $\mcl$ is given by the
class of essentially finite ones. Another one is given by the
class of rigid holomorphically nondegenerate submanifolds. Recall
that $(M,p)$ is called {\em rigid} if there exists normal
coordinates $(z,w)\in \Cn \times \Cd$ such that $M$ is given by
equation of the form \eqref{e:normalreal} with the additional
requirement that $\varphi$ does not depend on $s$ i.e.\ $\varphi
(z,\bar{z},s)=\varphi (z,\bar{z})$ (see \cite{BERbook}). The fact
that the class of rigid holomorphically nondegenerate
real-analytic generic submanifolds is contained in $\mcl$ follows
from Proposition \ref{p:classnormal} and Stanton's criterion of
holomorphic nondegeneracy whose statement we now recall.

\begin{prop}\label{p:stanton} $($\cite{S2, BERbook}$)$
Let $(M,p)$ be a germ of a real-analytic generic submanifold in $\CN$ and let $(z,w)$ be normal
  coordinates for $(M,p)$ so that $M$ is given by \eqref{e:normalcomplex}. Then the
  following are equivalent:
  \begin{enumerate} \item[{\rm (i)}] $(M,p)$ is holomorphically nondegenerate;
    \item[{\rm (ii)}] for some integer $k$, the local holomorphic map $(z,w) \mapsto \left( \bar{q}_\alpha (z,w) \colon | \alpha | \leq k \right)$
    is of generic rank $N$.
\end{enumerate}
\end{prop}

We now
exhibit examples of generic submanifolds that belong to the class
$\mcl$ without being essentially finite nor rigid holomorphically
nondegenerate.

\begin{exa}
  \label{ex:notessfin} Let $M\subset\C^2_z\times\C_w$ be defined by
  \[  \imag w = |z_1|^2 + (\real{w})\, |z_2|^2,\]
Let us check that $M$ is neither essentially finite nor rigid; by
using Proposition \ref{p:classnormal} one sees that it does belong to $\mcl$,
though. Moreover, we have
$\kappa_{M}(0)=2$. Indeed, rewriting in complex form, the defining equation 
of $M$ is given by
\[ w = \bar w \frac{1 + i |z_2|^2}{1- i |z_2|^2} +
2i \frac{|z_1|^2}{1-i|z_2|^2}= Q(z,\bar z, \bar w).\]
Thus, $Q(z,\bar z,0)$ is given by
\[ Q(z,\bar z,0) = 2i \sum_{j=0}^{\infty} i^j|z_1|^2 |z_2|^{2j}.\]
Using this expansion, we can now apply Proposition~\ref{p:classnormal} to see that
$M\in \mathcal{C}$; however, $M$ is not essentially finite since the
map $\chi\mapsto(\bar q_\alpha (\chi)) $  is--up to constant factors and
uninteresting components--just the map
\[ \chi \mapsto (\chi_1, \chi_1 \chi_2, \chi_1 \chi_2^2,\dots), \]
which is not finite.
To check that $M$ is not rigid, note that if $M$ is rigid, there exists a 
(germ of a real-analytic) curve
$(\R,0)\ni r\mapsto C(r) = (z(r),w(r))= (z(r),s(r) + i t(r))$ in $M$ 
with $C (0) = 0$ and $s'(0)\neq 0$ 
such that the germs $(M,0)$ and $(M,C(r))$  are biholomorphically equivalent for all $r$ close to $0$.  
Such a curve is obtained as follows: since $M$ is assumed to be rigid, in suitable local holomorphic 
coordinates $(\zeta,\eta)\in \C^2\times \C$, 
$M$ is given by an equation of the form $\imag \eta = \varphi(\zeta,\bar \zeta)$ with $\varphi$ real-analytic near the origin; in these
coordinates the claimed curve is just the line $(0,r)$, $r\in (\R,0)$. On the other hand,
 $(M,0)$ contains the complex line 
$(0,z_2,0)$ while for $(z,s+it)\in \C^2\times \C$ close to $0$  with $s\neq 0$ the germ $(M,(z,s+it))$ does not contain any complex curves
through $(z,s+it)$  since a calculation
in the given coordinates shows that $M$ is actually strongly pseudoconvex at all
such points. The germs $(M,0)$ and $(M,C(r))$ can therefore 
not be biholomorphically equivalent for $r\not =0$ close to $0$, which shows that $M$ is not rigid.
\end{exa}

Suppose now that $(M,0)$ is germ of a real-analytic
$k_0$-nondegenerate generic submanifold (as defined in \S
\ref{ss:nondegeneracy}), then obviously $(M,0)\in \mcl$ and one
has $\kappa_M(0)\leq k_0$. In the following example, we show that
this inequality may be strict. In fact, this example shows for
every $k \geq 2$ there exists a $k$-nondegenerate real-analytic
hypersurface $(M_k,0)$ in $\C^3$ for which $\kappa_{M_k}(0)=1<k$.
\begin{exa}\label{e:hitit}
  Let $(M_k,0) \subset \C^3_{(z_1,z_2,w)}$ be defined by
  \[ \imag w = \real \left( (z_1 + z_2) \bar{z_2}^k \right) + |z_1|^2.\]
\end{exa}
Moreover, analogous examples of real-analytic hypersurfaces or
generic submanifolds may be constructed in $\CN$ for arbitrary
$N$; this is left to the reader as an easy exercise. Similarly,
for essentially finite submanifolds, we have the following.

\begin{exa}
  For every $k\geq 3$, let $(\widetilde M_k, 0 )   \subset \C^3_{(z_1,z_2,w)}$ be defined by
  \[ \imag w = |z_1|^2 + |z_1 z_2|^2  + |z_2|^{2k}. \]
  Then $\widetilde M_k$ is essentially finite at the origin (but not finitely nondegenerate)
  and ${\rm esstype}_{\widetilde M_k}(0)=k$ while $\kappa_{M_k}(0) = 2<{\rm esstype}_{\widetilde M_k}(0)$.
\end{exa}

\section{Spaces of holomorphic maps in jet spaces}\label{s:jets}

Recall
that for any positive integer $k$, we denote the space of jets at the
origin of order $k$ of holomorphic mappings from $\C^N$ into
itself and fixing the origin by $J_{0,0}^k(\C^N)$. Given
coordinates $Z = (Z_1,\dots,Z_N)$ in $\CN$, we identify a jet
${\mathcal J}\in J_{0,0}^{k} (\CN)$ with a polynomial map of the
form
\begin{equation}\label{jet}
{\mathcal J}={\mathcal J} (Z)=\sum_{\alpha \in \N^r,\
1\leq|\alpha|\leq k}\frac{\Lambda^k_{\alpha}}{\alpha !}\, Z^{\alpha},
\end{equation}
where $\Lambda^k_{\alpha}\in \CN$. We thus have for a jet $
{\mathcal J}\in J_{0,0}^k(\C^N)$ coordinates
\begin{equation}\label{e:jetnotation}
\Lambda^k:=(\Lambda^k_{\alpha})_{1\leq|\alpha|\leq k}
\end{equation} as given in \eqref{jet}. Given a germ of a
holomorphic map $h\colon (\CN,0)\to (\CN,0)$, $h=h(t)$, we denote
by $j_t^kh$ the $k$-jet of $h$ at $t$ for $t$ sufficiently small
and also use the following splitting $j_t^kh=:(h(t),\widehat
j_t^kh)$. Moreover, since $h(0)=0$, we may also identify $j^k_0h$
with $\widehat j^k_0h$, that we will freely do in the sequel.

For later purposes, given a splitting $(z,w)\in\Cn\times\Cd$
(where $N = n+d$) of coordinates in $\CN$ (as for example induced
by a choice of normal coordinates \eqref{e:normalcomplex}),
 we introduce a special component of any jet $\Lambda^k
\in J_{0,0}^k(\C^N)$. We denote the set of multiindices of length
one having $0$ from the $n+1$-th to the $N$-th digit by $S$, and
the projection onto the first $n$ coordinates  by ${\proj}_1\colon
\CN\to \Cn$ (that is, ${\proj}_1(z,w)=z$). Then
\begin{equation}\label{e:crapagain}
\widetilde \Lambda^1:=({\proj}_1(\Lambda_{\alpha}))_{\alpha \in
S}.
\end{equation}
Note that for any local holomorphic map $(\C^n \times \Cd,0) \ni
(z,w) \mapsto H(z,w)=(f(z,w),g(z,w))\in (\Cn \times \Cd,0)$, if
$\jetm{0}{k}H=\Lambda^k$, then $\widetilde
\Lambda^1=(\frac{\partial f}{\partial z}(0))$. We can therefore
identify $\widetilde \Lambda^1$ with an $n\times n$ matrix or
equivalently with an element of $J_{0,0}^1(\Cn)$. Throughout the
paper, given any other notation for a jet $\lambda^k\in
J_{0,0}^k(\CN)$, $\widetilde \lambda^1$ will always denote the
component of $\lambda^k$ defined above. We denote by ${\mathcal
G}_0^k(\C^N)$ the connected open subset of $J_{0,0}^k(\C^N)$
consisting of those jets $\Lambda^k$ for which ${\det}\ \widetilde
\Lambda^1\not = 0$ i.e.\ for which $\widetilde \Lambda^1\in \glnc
$.

For later purposes, it is convenient to introduce the following
spaces of functions. For any positive integers $k,s$, let
${\mathcal F}_s^k = \mathcal{O} (\{0\}\times {\mathcal G}^1_0(\CN)\times {\mathcal G}^k_0(\CN))$ 
be the ring 
of germs of holomorphic functions on 
$\{0\}\times {\mathcal G}^1_0(\CN)\times {\mathcal G}^k_0(\CN) \subset
\C^s \times {\mathcal G}^1_0(\CN)\times {\mathcal G}^k_0(\CN)$; recall
that this is just the space of holomorphic functions which are
defined on some connected open subset (which may depend on 
the function) of $\C^s \times
{\mathcal G}^1_0(\CN)\times {\mathcal G}^k_0(\CN)$ containing
$\{0\}\times {\mathcal G}^1_0(\CN)\times {\mathcal G}^k_0(\CN)$ where we 
identify two such functions if there exists an open connected subset 
containing $\{0\}\times {\mathcal G}^1_0(\CN)\times {\mathcal G}^k_0(\CN)$ 
on which 
they are both defined and agree. 
For any $\Theta=\Theta
(t,\lambda^1,\Lambda^k)\in {\mathcal F}_s^k$, we write the Taylor
expansion
$$ \Theta (t,\lambda^1,\Lambda^k)=\sum_{\nu \in \N^s}\Theta_{\nu}
(\lambda^1,\Lambda^k)t^{\nu}.$$
 Let ${\mathcal E}^k_s$ denote the subspace of ${\mathcal F}^k_s$
consisting of those germs $\Theta\in {\mathcal F}^k_s$ for
which $\Theta_{\nu}(\lambda^1,\Lambda^k)$ can be written, for any
$\nu\in \N^s$, in the following form
$$\frac{B_{\nu}(\lambda^1,\Lambda^k)}{({\det}\, \widetilde \Lambda^1)^{c_{\nu}}
({\det}\, \widetilde \lambda^1)^{a_{\nu}}},$$ for some
polynomial $B_{\nu}$ and some nonnegative integers
$c_{\nu},a_\nu$.

In what follows, for $x\in \C^q$ and $y\in \C^r$, we denote by $\C\{x\}[y]$ the ring of polynomials in $y$ with coefficients that are convergent power series in $x$.

Throughout the paper, we will need the following elementary facts about the space
${\mathcal E}^k_s$ that we shall freely use.

\begin{lem}\label{l:heavystuff} Let $\Theta_1,\ldots,\Theta_r,\Theta_{r+1},\ldots, \Theta_q$
be $q$ holomorphic functions belonging to ${\mathcal E}^k_s$ such
that $\Theta_{r+1}
(0,\lambda^1,\Lambda^k)=\ldots=\Theta_q(0,\lambda^1,\Lambda^{k})=0$
for all $(\lambda^1,\Lambda^k)\in {\mathcal G}^1_0(\CN)\times
{\mathcal G}^k_0(\CN)$. Then the following hold.
\begin{enumerate}
\item[{\rm (i)}] If $V(t,T;y)\in \C\{t,T\}[y]$, $(t,T,y)\in \C^s\times \C^{q-r}\times \C^r$, 
then the holomorphic function
\[V(t,\Theta_{r+1}(t,\lambda^1,\Lambda^k),\ldots,\Theta_{q}(t,\lambda^1,\Lambda^k);
\Theta_{1}(t,\lambda^1,\Lambda^k),\ldots,\Theta_{r}(t,\lambda^1,\Lambda^k))\]
also belongs to ${\mathcal E}^k_s$.
\item[{\rm (ii)}] If $C(t;y)\in \C\{t\}[y]$, 
  $D(t;y)\in \C\{t\}[y]$, $(t,y)\in \C^s\times \C^r$, and for some
nonnegative integers $a,b$
  $D(0;\lambda^1,\Lambda^k)= ({\det}\, \widetilde
\Lambda^1)^a\ ({\det}\, \widetilde \lambda^1)^b$, then for every nonnegative integer $m$
the ratio
$$\frac{C(t;
\Theta_{1}(t,\lambda^1,\Lambda^k),\ldots,\Theta_{r}(t,\lambda^1,\Lambda^k))}{(D(t;
\Theta_{1}(t,\lambda^1,\Lambda^k),\ldots,\Theta_{r}(t,\lambda^1,\Lambda^k)))^m}$$
defines an element of ${\mathcal E}^k_s$.
\end{enumerate}
\end{lem}

%%%%%%%%%%%%%%%%%%%%%%%%%%%%%%%%%%%%%%%%%%%%%%%%%%%%%%%%%%%%%%%%%%%%%%%%%%%%%%%%%%%%%%%
%%%%%%%%%%%%%%%%%%%%%%%%%%%%%%%%%%%%%%%%%%%%%%%%%%%%%%%%%%%%%%%%

\section{Reflection identities}\label{s:reflectionid}

Throughout the rest of this paper, we let $(M,0)$ denote a germ of
a real-analytic generic submanifold of $\CN$ through the origin.
This section is the starting point of the proof of Theorem
\ref{t:further1}. We shall assume, without loss of generality,
that $M$ is given in normal coordinates $Z=(z,w)\in \Cn \times
\Cd$ as in \eqref{e:normalreal} and \eqref{e:normalcomplex}.
Recall also that the complexification $\M \subset \CN \times \CN$
is the complex submanifold given through the origin by
\eqref{e:complexcoord}. In what follows, we shall keep the
notation and terminology introduced in \S \ref{s:class}-\S
\ref{s:jets}.

This section is devoted to the first of the three main steps of
the proof of Theorem \ref{t:further1}. We establish here some
general identities on the complexification $\M$ of $M$ (see
Propositions \ref{p:reflection1} and \ref{p:reflection3}), which
are valid without any condition on the submanifold $M$.

Recall that for any positive integer $k$, $\Lambda^k$ denotes the
coordinates in $J_{0,0}^{k}(\CN)$ as explained in \S \ref{s:jets},
and $\widetilde \Lambda^1$ is the component of $\Lambda^1$ given
by \eqref{e:crapagain}, corresponding to our choice of normal
coordinates $(z,w)$. In what follows, for every $H\in \autM$,
following the splitting of the coordinates $Z=(z,w)$, we split the
map $H$ as $H=(f,g)\in \Cn\times \Cd$.

For the map $\bar{Q}$ previously defined, we will use the
following notation to simplify our formulas:
\[\bar{Q}_Z = \begin{pmatrix}
  \bar{Q}^1_{Z_1} &\dots &   \bar{Q}^d_{Z_1} \\
  \vdots & \ddots & \vdots \\
  \bar{Q}^1_{Z_N} &\dots &   \bar{Q}^d_{Z_N}
\end{pmatrix},
\quad
 \bar{Q}_{\chi^\alpha,Z}=, \begin{pmatrix}
  \bar{Q}^1_{\chi^\alpha,Z_1} &\dots &   \bar{Q}^d_{\chi^\alpha,Z_1} \\
  \vdots & \ddots & \vdots \\
  \bar{Q}^1_{\chi^\alpha,Z_N} &\dots &   \bar{Q}^d_{\chi^\alpha,Z_N}
\end{pmatrix},\]
where $\alpha \in \N^n$ and
$\bar{Q}=(\bar{Q}^1,\ldots,\bar{Q}^d)$.

 We start with the
following well-known computational lemma (see e.g. \cite[Equation (3.4.14)]{BER4}).

\begin{prop}\label{p:reflection1} There exists a polynomial $\D=\D(Z,\zeta, \Lambda^1)\in \C \{Z,\zeta\}[\Lambda^1]$ which is universal $($i.e. only
  depends on $M)$, and, for every $\alpha \in \N^n\setminus \{0\}$, another
universal $\Cd$-valued polynomial map
$\PP_{\alpha}=\PP_{\alpha}(Z,\zeta,\Lambda^{|\alpha|})$ whose components are in the ring   $\C\{Z,\zeta\}[\Lambda^{|\alpha|}]$  such that for any $H\in \autM$, the following holds:
\begin{enumerate}
\item[{\rm (i)}] $\D (0,0,\Lambda^1)={\det}\, \widetilde \Lambda^1$;
\item[{\rm (ii)}] $\D (0,0,\jetm{0}{1} H)\not =0$;
\item[{\rm (iii)}] for all $(Z,\zeta)\in \M$ near $0$,
\begin{equation}\label{e:fundamental1}
(\D(Z,\zeta,\widehat j_{\zeta}^{1} \overline{H} ))^{2|\alpha|-1}\,
\bar{Q}_{\chi^{\alpha}}(\bar{f}(\zeta),H(Z))=\PP_{\alpha}(Z,\zeta,\widehat
j_{\zeta}^{\left| \alpha \right|} \overline{H}).
\end{equation}
\end{enumerate}
\end{prop}

Our next identity is concerned with pure transversal derivatives
of every element of $\autM$.

\begin{prop}\label{p:reflection2}
For any $\mu \in \N^d\setminus \{0\}$ and $\alpha \in
\N^n\setminus\{0\}$, there exist a universal $\Cd$-valued
polynomial map ${\mathcal T}_{\mu,\alpha}(Z,\zeta,Z',\zeta',\lambda^{|\mu|-1},\Lambda^{|\mu|})$ whose components belong
to the ring $\C \{Z,\zeta,Z',\zeta'\}[\lambda^{|\mu|-1},\Lambda^{|\mu|}]$ and another universal $\Cd$-valued
polynomial map ${\mathcal Q}_{\mu, \alpha}(Z,\zeta,\Lambda^{|\alpha|+|\mu|})$ whose components are in the ring $\C \{Z,\zeta\}[\Lambda^{|\alpha|+|\mu|}]$, such that for any $H\in \autM$ and for any
$(Z,\zeta)\in \M$ close to the origin, the following relation
holds:
\begin{equation}\label{e:fundamental2}
\frac{\partial^{|\mu|}H}{\partial w^\mu}(Z)\cdot
\bar{Q}_{\chi^{\alpha},Z}(\bar{f}(\zeta),H(Z)) =(*)_1+(*)_2,
\end{equation}
where \begin{equation}\label{e:*1} (*)_1:= {\mathcal T}_{\mu,
\alpha}\left(Z,\zeta,H(Z),\overline{H}(\zeta), \widehat
j_{Z}^{\left| \mu \right| -1} H , \widehat j_{\zeta}^{\left| \mu
\right|} \overline{H} \right)
\end{equation} and
\begin{equation}\label{e:*2}
(*)_2:=\frac{{\mathcal Q}_{\mu,\alpha}(Z,\zeta, \widehat
j_{\zeta}^{|\alpha| + |\mu|} \overline{H})}{(\D(Z,\zeta, \widehat
j_{\zeta}^{1} \overline{H} ))^{2|\alpha|+|\mu|-1}}
\end{equation}
\end{prop}

\begin{proof} Consider the following holomorphic vector fields
tangent to $\M$:
\begin{equation}\label{e:Rj}
R_j=\frac{\partial}{\partial w_j}+\bar{Q}_{w_j}(\chi,z,w)\cdot
\frac{\partial}{\partial \tau}= \frac{\partial}{\partial
w_j}+\sum_{k=1}^d\bar{Q}^k_{w_j}(\chi,z,w)
\frac{\partial}{\partial \tau_k},\ j=1,\ldots,d.
\end{equation}

We prove the proposition by induction on the length of $\mu$. For
$\mu \in \N^d$ of length one, say, without loss of generality, for
$\mu =(1,0,\ldots,0)$, \eqref{e:fundamental2} follows easily after
applying $R_1$ to \eqref{e:fundamental1}. Suppose now that
\eqref{e:fundamental2} holds for all multiindeces $\mu$ of a
certain length $k$. Then if $\mu'\in \N^d$ is of length $k+1$, say
again $\mu'=(1,0,\ldots,0)+\mu$ for some $\mu\in \N^d$ of length
$k$, we apply again $R_1$ to \eqref{e:fundamental2}, and one
easily sees from the chain rule that \eqref{e:fundamental2} is
satisfied with $\mu$ replaced by $\mu'$. We leave the details to
the reader.
\end{proof}

The next lemma deals with the transversal derivatives of the normal
component $g$ of any $H=(f,g)\in \autM$.

\begin{lem}\label{l:gderivative}
 For any $\mu \in \N^d\setminus \{0\}$, there exist a universal $\Cd$-valued polynomial map $W_\mu={W}_{\mu}(Z,\zeta,Z',\zeta',\lambda^{|\mu|-1},\Lambda^{|\mu|})$ whose components belong
to the ring $\C \{Z,\zeta,Z',\zeta'\}[\lambda^{|\mu|-1},\Lambda^{|\mu|}]$  such that for every $H=(f,g)\in
  \autM$ and for any $(Z,\zeta)\in \M$ sufficiently close to the
  origin,
\begin{equation}\label{e:gmu}
\frac{\partial^{|\mu|}g}{\partial
w^\mu}(Z)=\frac{\partial^{|\mu|}f}{\partial w^\mu}(Z)\cdot
Q_z(f(Z),\overline{H}(\zeta))+(*)_3
\end{equation}
with
\begin{equation}\label{e:*3}
(*)_3:=W_\mu \left(Z,\zeta,H(Z),\overline{H}(\zeta), \widehat
j_{Z}^{\left| \mu \right| - 1}H , \widehat j_{\zeta}^{\left| \mu
\right| } \overline{H} \right).
\end{equation}
\end{lem}

\begin{proof} We again prove the lemma by induction on $|\mu|$.

Let $\mu\in \N^d$ of length one, say $\mu =(1,0,\ldots,0)$.
Applying the vector field $R_1$ (as defined in \eqref{e:Rj}) to
the equation $g(Z)=Q(f(Z),\overline{H}(\zeta))$, for $(Z,\zeta)\in
\M$, one easily sees that the existence of $W_\mu$ satisfying
\eqref{e:gmu} follows from the chain rule.

Suppose now that \eqref{e:gmu} holds for all multiindeces $\mu$ of
a certain length $k$. Then if $\mu'\in \N^d$ is of length $k+1$,
say again $\mu'=(1,0,\ldots,0)+\mu$ for some $\mu\in \N^d$ of
length $k$, we apply again $R_1$ to \eqref{e:gmu}, and the
existence of $W_{\mu'}$ satisfying the desired equality follows
again from the chain rule. We leave the details to the reader.
\end{proof}

Combining Proposition \ref{p:reflection2} and Lemma
\ref{l:gderivative}, one gets the following.

\begin{prop}\label{p:reflection3}
For any $\mu \in \N^d\setminus \{0\}$ and $\alpha \in
\N^n\setminus\{0\}$, there exist a universal $\Cd$-valued
polynomial map ${\mathcal T}'_{\mu,\alpha}(Z,\zeta,Z',\zeta',\lambda^{|\mu|-1},\Lambda^{|\mu|})$ whose components belong
to the ring $\C \{Z,\zeta,Z',\zeta'\}[\lambda^{|\mu|-1},\Lambda^{|\mu|}]$  such that for any $H=(f,g)\in \autM$ and
for any $(Z,\zeta)\in \M$ close to the origin, the following
relation holds:
\begin{equation}\label{e:fundamental3}
\frac{\partial^{|\mu|}f}{\partial w^\mu}(Z)\cdot
\left(\bar{Q}_{\chi^{\alpha},z}(\bar{f}(\zeta),H(Z))+
Q_z(f(Z),\overline{H}(\zeta))\cdot
\bar{Q}_{\chi^{\alpha},w}(\bar{f}(\zeta),H(Z))\right)=(*)_1'+(*)_2,
\end{equation}
where $(*)_2$ is given by \eqref{e:*2} and $(*)_1'$ is given by
\begin{equation}\label{e:*1'}
(*)_1':={\mathcal
T}_{\mu,\alpha}'\left(Z,\zeta,H(Z),\overline{H}(\zeta), \widehat
j_{Z}^{\left| \mu \right| - 1} H, \widehat j_{\zeta}^{\left| \mu
\right|}{\overline{H}} \right).
\end{equation}
\end{prop}

\begin{rem}
In \eqref{e:fundamental3}, one should see $\displaystyle
\frac{\partial^{|\mu|}f}{\partial w^\mu}(Z)$ as a $1\times n$
matrix, $\bar{Q}_{\chi^{\alpha},z}(\bar{f}(\zeta),H(Z))$ as a
$n\times d$ matrix, $Q_z(f(Z),\overline{H}(\zeta))$ as a $n\times
d$ matrix and
$\bar{Q}_{\chi^{\alpha},w}(\bar{f}(\zeta),H(Z))=(\bar{Q}^j_{\chi^{\alpha},w_i}(\bar{f}(\zeta),H(Z)))_{1\leq
i,j\leq d}$ as a $d\times d$ matrix.
\end{rem}

%%%%%%%%%%%%%%%%%%%%%%%%%%%%%%%%%%%%%%%%%%%%%%%%%%%%%%%%%%%%%%%%%%%%%%%%%%%%%%%%%%%%

\section{Parametrization of $\autM$ along Segre sets}\label{s:jetpar}
This section is devoted to the second part of the proof of Theorem
\ref{t:further1}, which involves two main ingredients: the Segre
sets and mappings introduced by Baouendi, Ebenfelt and Rothschild
\cite{BER2} and the parametrization theorems for singular analytic
systems developed
in the first part of the paper. We shall now prove
that for a real-analytic generic submanifold $M\in\mcl$ (the class
$\mcl$ is defined in \S \ref{s:class}), the elements of $\autM$,
restricted to any Segre set, are parametrized by their jets at the
origin. Such a property is in fact first shown to be true on the
Segre variety attached to the origin (in Proposition
\ref{p:firststep}) and then we establish that this parametrization
property ``propagates" to higher order Segre sets (this is done in
Proposition \ref{p:propagation}). For both of these results, we
now use the fact that $M$ belongs to the class $\mcl$.

\subsection{Segre sets and mappings}

We start by recalling the Segre mappings as introduced in
\cite{BER2}. They play a fundamental role in the study of mappings
between real-analytic or real-algebraic CR-manifolds of arbitrary
codimension \cite{BERbook, BERbull, Rsurvey}. In what follows, for
any positive integer $s$, $t^s:=(t^s_1,\ldots,t^s_n)$ will denote
an $n$-dimensional complex variable.

For $j\in \N^*$, we define the Segre mappings
$v^j\colon (\C^{nj},0)\to (\CN,0)$ inductively: For $j = 1$,
$v^1(t^1):=(t^1,0)$; and for $j\geq 1$,

\begin{equation}\label{e:segresetsmaps}
v^{j+1}(t^1,\ldots,t^j):=(t^{j+1},Q(t^{j+1},\bar{v}^{j}(t^1,\ldots,t^j))).
\end{equation}
By definition, given a sufficiently small neighbourhood $U_j$ of
the origin in $\C^{nj}$, the Segre set is the image of the
neighborhood $U_j$ under the map $v^j$ (for a thorough discussion,
see e.g. \cite{BERbook}). Note that for any $j\geq 1$, the map
$$(\C^{n(j+1)},0)\ni (t^1,\ldots,t^j,t^{j+1})\mapsto
(v^{j+1}(t^1,\ldots,t^{j+1}),\bar{v}^j(t^1,\ldots,t^j))\in
(\CN\times \CN,0)$$ takes in fact its values on the
complexification $\M$. In what follows, for any $j\geq 1$, we
shall use the following useful notation
\begin{equation}\label{e:usefulnotation}
t^{[j]}:=(t^1,\ldots,t^j).
\end{equation}

\subsection{Jet parametrization of CR automorphisms along the
first Segre set}\label{ss:firstorder}
 In what follows, we keep the
notation introduced at the beginning of this section as well as in \S
\ref{s:reflectionid}.

The next proposition gives  a precise parametrization property by
jets at the origin of all elements of $\autM$ and of their
derivatives along the first Segre set. Note that for any $H\in
\autM$ and for any positive integers $k,s$, $(0,\jetm{0}{1}
\overline{H},\jetm{0}{k} H)$ and
$(0,\jetm{0}{1}H,\jetm{0}{k}\overline{H})$ are in the domain of
definition of any holomorphic function belonging to the space
${\mathcal F}_s^k$ as defined in \S \ref{s:jets}.

\begin{prop}\label{p:firststep}
 In the above setting and with the above notation,
assume furthermore that $(M,0)$ belongs to the class $\mcl$  and
denote by $\kappa_0:=\kappa(M,0)$ the invariant integer attached
to $(M,0)$ as defined in \S \ref{s:class}. Then for any multiindex
$\beta \in \N^N$, there exists a $\CN$-valued holomorphic map
\begin{multline}
\Psi_1^{\beta}=\Psi_1^{\beta}\left(
t^1,\lambda^1,\Lambda^{\kappa_0+|\beta|}\right):=\\
\left( \widetilde\Psi_1^{\beta}\left(
t^1,\lambda^1,\Lambda^{\kappa_0+|\beta|}\right),
\widehat\Psi_1^{\beta}\left(t^1,\lambda^1,\Lambda^{\kappa_0+|\beta|}\right)\right)
\in
\C^n\times \Cd,
\end{multline}
whose components belong to ${\mathcal E}_{n}^{\kappa_0+|\beta|}$,
satisfying the following properties:
\begin{enumerate}
\item[{\rm (i)}] If we denote by $S\subset \N^N$ the set of multiindices of length one
having $0$ from the $(n+1)$-th to the $N$-th component, then for any
$(\lambda^1,\Lambda^{\kappa_0})\in {\mathcal G}_0^1(\CN)\times
{\mathcal G}^{\kappa_0}_0(\CN)$ we have
\begin{equation}\label{e:conditions}
\Psi_1^0(0,\lambda^1,\Lambda^{\kappa_0})=0,\ (\widetilde
\Psi_1^{\beta}(0,\lambda^1,\Lambda^{\kappa_0}))_{\beta \in S}=
\frac{\partial \widetilde \Psi_1^0}{\partial
t^1}(0,\lambda^1,\Lambda^{\kappa_0})=\widetilde \lambda^1.
\end{equation}
\item[{\rm (ii)}] For every $H\in \autM$, we have
\begin{equation}
((\partial^{\beta}H)\circ v^1) (t^1)=
\Psi_1^{\beta}(t^1,\jetm{0}{1}H,j_0^{\kappa_0+|\beta|}\overline{H}),\quad
\end{equation}
for all $t^1$ sufficiently close to the origin in $\Cn$.
\end{enumerate}
\end{prop}

 We prove the proposition by induction on $|\beta|$. The induction start
 $\beta=0$ is treated in \S \ref{ss:first0} and the induction step is
 carried out in \S \ref{ss:firstinduction}.

\subsubsection{Proof of Proposition \ref{p:firststep} for
$|\beta|=0$}\label{ss:first0} For $\beta=0$, we first note that it
follows from the normality of the chosen coordinates $(z,w)$ that
for any $H=(f,g)\in \autM$, one has $(g\circ v^1)(t^1)\equiv
0$. We may therefore define $\widehat
\Psi_1^0(t^1,\lambda^1,\Lambda^{\kappa_0}):=0$, and it is enough
to construct
$\widetilde \Psi_1^0(t^1,\lambda^1,\Lambda^{\kappa_0})$.
Setting $Z=v^1(t^1)$ and $\zeta=0$ in \eqref{e:fundamental1}
implies that for any $H\in \autM$
and for any $\alpha \in \N^n\setminus \{0\}$,
\begin{equation}
\bar{Q}_{\chi^{\alpha}}(0 ,(H\circ
v^{1})(t^1))=\frac{\PP_{\alpha}(v^{1}(t^1),0, \widehat j_0^{\left|
\alpha \right|} \overline{H} )} {(\D(v^{1}(t^1),0, \widehat
j_{0}^{1} \overline{H} ))^{2|\alpha|-1}},
\end{equation}
for all $t^1$ sufficiently close to the origin. Using
\eqref{e:normalexp} and the fact that $g\circ v^1\equiv 0$, we
obtain equivalently
\begin{equation}\label{e:use1}
\alpha!\, \bar{q}_{\alpha}(f\circ
v^{1}(t^1),0)=\frac{\PP_{\alpha}(v^{1}(t^1),0, \widehat
j_{0}^{\left| \alpha \right|} \overline{H}) } {(\D(v^{1}(t^1),0,
\widehat j_{0}^{1} \overline{H} ))^{2|\alpha|-1}}.
\end{equation}
In what follows, we use \eqref{e:use1} only for multiindices
$\alpha\in \N^n\setminus \{0\}$ with $|\alpha|\leq \kappa_0$. For
such $\alpha$'s, we set
\begin{equation}\label{e:use2}
p_{\alpha}(t^1,\Lambda^{\kappa_0}):=\PP_{\alpha}(v^{1}(t^1),0,\Lambda^{|\alpha|})
,\quad \Delta_{\sigma}(t^1,
\Lambda^1):=(\D(v^{1}(t^1),0,\Lambda^1))^{\sigma},\ \sigma \in \N,
\end{equation}
and
\begin{equation}\label{e:use3}
r_{\alpha}(t^1,\Lambda^{\kappa_0}):=\frac{p_{\alpha}(t^1,\Lambda^{\kappa_0})}{\Delta_{2|\alpha|-1}(t^1,\Lambda^{1})}.
\end{equation}
It follows from Lemma \ref{l:heavystuff} and Proposition
\ref{p:reflection1} that the map
$r_\alpha=(r^1_{\alpha},\ldots,r^d_{\alpha})$ defines a
$\Cd$-valued holomorphic map in a neighborhood of $\{0\} \times
{\mathcal G}^{\kappa_0}_0(\CN)$ whose components belong to
${\mathcal E}_n^{\kappa_0}$. Moreover, since $M\in\mcl$ Proposition \ref{p:classnormal} implies that
we may choose $n$
integers $i_1,\ldots,i_n\in \{1,\ldots,d\}$ and $n$ multiindeces
$\alpha^{(1)},\ldots,\alpha^{(n)}\in \N^n\setminus \{0\}$ of
length $\leq \kappa_0$ so that $z\mapsto
(\bar{q}^{i_1}_{\alpha^{(1)}}(z,0),\ldots,\bar{q}^{i_n}_{\alpha^{(n)}}(z,0))$
is of generic rank $n$. Consider the following system of
complex-analytic equations in the unknown $T^1\in \C^{n-1}$
\begin{equation}\label{e:system1}
\alpha^{(\nu)} !\,
\bar{q}^{i_{\nu}}_{\alpha^{(\nu)}}(T^1,0)=r^{i_\nu}_{\alpha^{(\nu)}}(t^1,\Lambda^{\kappa_0}),\
\nu=1,\ldots,n.
\end{equation}
By Corollary \ref{c:corparcn1}, there exists a holomorphic map
\begin{equation}\label{e:gamma}
\Gamma_1^0=\Gamma_1^0
(t^1,\widetilde\lambda^1,\Lambda^{\kappa_0})\colon \C^n\times
\glnc \times {\mathcal G}_0^{\kappa_0}(\CN)\to \C^n
\end{equation} defined in a neighbourhood of $\{0\}\times \glnc \times {\mathcal G}_0^{\kappa_0}(\CN)$, satisfying
$$\Gamma_1^0 (0,\widetilde\lambda^1,\Lambda^{\kappa_0})=0,\quad
\frac{\partial \Gamma_1^0}{\partial
t^1}(0,\widetilde\lambda^1,\Lambda^{\kappa_0})=\widetilde
\lambda^1$$
 such that if $u\colon (\Cn,0)\to (\Cn,0)$ is a germ of
a biholomorphism satisfying \eqref{e:system1} for some
$\Lambda^{\kappa_0}\in {\mathcal G}_0^{\kappa_0}(\CN)$, then
$u(t^1)=\Gamma_1^0 (t^1,\jetm{0}{1}u,\Lambda^{\kappa_0})$. Hence, since
for any $H\in \autM$, $t^1\mapsto (f\circ v^1)(t^1)$ is a germ at
0 of a biholomorphism satisfying \eqref{e:system1} with
$\Lambda^{\kappa_0}=j_0^{\kappa_0}\overline{H}$ in view of
\eqref{e:use1}, \eqref{e:use2} and \eqref{e:use3}, we have
\begin{equation}\label{e:bored1}
(f\circ v^1)(t^1)=\Gamma_1^0 (t^1,j^1_0(f\circ v^1
),j_0^{\kappa_0}\overline{H})=\Gamma_1^0 \left(t^1,\left(\frac{\partial
f}{\partial z }(0)\right),j_0^{\kappa_0}\overline{H}\right).
\end{equation}
Therefore if we set $\widetilde
\Psi_1^0(t^1,\lambda^1,\Lambda^{\kappa_0}):=\Gamma_1^0(t^1,\widetilde
\lambda^1,\Lambda^{\kappa_0})$, we are left to prove that the
components of $\widetilde
\Psi_1^0(t^1,\lambda^1,\Lambda^{\kappa_0})$ belong to ${\mathcal
E}_n^{\kappa_0}$. For this, we write the taylor expansion
\begin{equation}\label{e:close}
\widetilde
\Psi_1^0(t^1,\lambda^1,\Lambda^{\kappa_0})=\Gamma_1^0(t^1,\widetilde
\lambda^1,\Lambda^{\kappa_0})=\sum_{\gamma \in
\N^n}\Gamma^0_{1,\gamma}(\widetilde
\lambda^1,\Lambda^{\kappa_0})(t^1)^{\gamma}. \end{equation}
 From Corollary \ref{c:corparcn1} (ii) we know that for each $\gamma \in \N^n$,
there exist positive integers $e_\gamma$ and $s_\gamma$ such that
$\Gamma^0_{1,\gamma}(\widetilde \lambda^1,\Lambda^{\kappa_0})$ can
be written in the form
\begin{equation}\label{e:taylor1}
\frac{\rho_\gamma \left(\widetilde
\lambda^1,\left((\partial_{t^1}^{\nu}
r_{\alpha})(0,\Lambda^{\kappa_0})\right)_{\substack{1\leq|\alpha|\leq
\kappa_0\\ |\nu|\leq s_\gamma}} \right)}{({\det}\, \widetilde
\lambda^1)^{e_{\gamma}}}
\end{equation}
where $\rho_\gamma$ is a $\Cn$-valued polynomial map of its
arguments. Since each component of the map $r_{\alpha}$,
$|\alpha|\leq \kappa_0$, belongs to ${\mathcal E}_{n}^{\kappa_0}$,
it follows from \eqref{e:taylor1} and \eqref{e:close} that the
same fact holds for every component of $\widetilde \Psi_1^0$. This
completes the proof of the Proposition \ref{p:firststep} in the
case $\beta =0$.

\subsubsection{Induction for Proposition \ref{p:firststep}}\label{ss:firstinduction}
Assume now that Proposition \ref{p:firststep} holds for all
multindices $\beta\in \N^N$ of length $|\alpha|\leq k$, $k\in \N$,
and let $\delta=(\delta_1,\ldots,\delta_N)\in \N^N$ with
$|\delta|=k+1$. Suppose first that there exists $m\in
\{1,\ldots,n\}$ such that $\delta_m>0$. Then by setting
\begin{equation}\label{e:cheat}
\Psi_1^{\delta}(t^1,\lambda^1,\Lambda^{\kappa_0+k+1}):=\frac{\partial
\Psi_1^{\beta}}{\partial
t^1_m}(t^1,\lambda^1,\Lambda^{\kappa_0+k})
\end{equation} where $t^1=(t^1_1,\ldots,t^1_n)\in \C^n$ and $\beta =
(\beta_1,\ldots,\beta_N)$, $\beta_i=\delta_i$ for $i\not =m$ and
$\beta_m=\delta_m-1$, we see that the derivative of order $\delta$
of every $H\in\autM$ is parametrized in the right way and we are
done in this case. Equation \eqref{e:cheat} also shows that the
parametrizations $\Psi_1^{\delta}$ with $\delta\in S$  are
constructed in such a way that the second part of
\eqref{e:conditions} holds. We may now therefore assume that
$\delta_1=\ldots=\delta_n=0$ (and $|\delta|=k+1$). We set $\mu
=(\delta_{n+1},\ldots,\delta_N)\in \N^d$ and we note that
$|\mu|=|\delta|=k+1$. Using \eqref{e:fundamental3} and
\eqref{e:*1'} in which we set $Z=v^1(t^1)$ and $\zeta=0$ and using
the fact that $Q(z,0,0)=0$ (normality of the coordinates), we
obtain the following equality for every $\alpha \in \N^n\setminus
\{0\}$ and every $H=(f,g)\in \autM$ and every $t^1\in \Cn$
sufficiently close to the origin
\begin{multline}\label{e:formula1}
\alpha!\, ((\partial_{w}^{\mu}f)\circ v^1)\cdot  \frac{\partial
\bar{q}_{\alpha}}{\partial z}(f\circ v^1,0) =\frac{{\mathcal
Q}_{\mu,\alpha}(v^1,0, \widehat j_{0}^{|\alpha| + k + 1}
\overline{H}) }{ \Delta_{2|\alpha|+k}(t^1, \widehat j_{0}^{1}
\overline{H})
}\\
+{\mathcal T}_{\mu,\alpha}'\left(v^1,0,H \circ v^1,0,
((\partial^{\beta} H)\circ v^1 \right)_{1\leq |\beta|\leq
k},\jetm{0}{k+1} \overline{H}).
\end{multline}
We have deliberately dropped the variable $t^1$ in
\eqref{e:formula1} to simplify a bit the formula. In what follows,
we restrict ourselves to multiindices $\alpha\in \N^{d}$ such that
$|\alpha|\leq \kappa_0$. For such $\alpha$'s, we set
\begin{equation}\label{e:xi}
\Xi_{\mu,\alpha} (t^1,\Lambda^{\kappa_0+k+1}):=\frac{{\mathcal
Q}_{\mu,\alpha}(v^1(t^1),0,\Lambda^{|\alpha|+k+1})}{
\Delta_{2|\alpha|+k}(t^1,\Lambda^1)},
\end{equation}
\begin{multline}\label{e:upsilon}
\Upsilon_{\mu,\alpha}
(t^1,\lambda^1,\Lambda^{\kappa_0+k+1}):=\\{\mathcal T}_{\mu,\alpha}'
\left(v^1(t^1),0,\Psi_1^0(t^1,\lambda^1,\Lambda^{\kappa_0}),0,\Psi^{(k)}_1(t^1,\lambda^1,\Lambda^{\kappa_0+k}),\Lambda^{k+1}\right).
\end{multline}
In the last equation we use the notation
\[ \Psi^{(k)}_1(t^1,\lambda^1,\Lambda^{\kappa_0+k}):=
\left( \Psi^{\beta}_1(t^1,\lambda^1,\Lambda^{\kappa_0+|\beta|})
\right)_{1\leq |\beta|\leq k}.\] It follows from \eqref{e:use2},
\eqref{e:conditions}, Lemma \ref{l:heavystuff} and Proposition
\ref{p:reflection1} that the mappings $\Xi_{\mu,\alpha}$ and
$\Upsilon_{\mu,\alpha}$ define holomorphic mappings in an open
neighborhood of $\{0\}\times {\mathcal G}^1_0(\CN)\times {\mathcal
G}_0^{\kappa_0+k+1} (\CN)$ in $\Cn \times J^1_{0,0}(\CN)\times
J_{0,0}^{\kappa_0+k+1}(\CN)$, whose components belong to
${\mathcal E}_n^{\kappa_0+k+1}$. Consider the following linear
singular system in $X\in \Cn$:
\begin{equation}\label{e:system2}
\alpha^{(\nu)}!\ X\cdot \frac{\partial
\bar{q}^{i_\nu}_{\alpha^{(\nu)}}}{\partial z}(\widetilde
\Psi^0_1(t^1,
\lambda^1,\Lambda^{\kappa_0}),0)=\Xi_{\mu,\alpha^{(\nu)}}(t^1,\Lambda^{\kappa_0+k+1})+\Upsilon_{\mu,\alpha^{(\nu)}}
(t^1,\lambda^1,\Lambda^{\kappa_0+k+1}),
\end{equation}
for $\nu=1,\ldots,n$. By \eqref{e:formula1} and the induction
assumption, for any $H\in \autM$, $X=(\partial_{w}^{\mu}f)\circ
v^1$ satisfies the above system with $\lambda^1=\jetm{0}{1} H$ and
$\Lambda^{\kappa_0+k+1}=\jetm{0}{\kappa_0+k+1}\overline{H}$.  Now,
the holomorphic map $z\mapsto
(\bar{q}^{i_1}_{\alpha^{(1)}}(z,0),\ldots,\bar{q}^{i_n}_{\alpha^{(n)}}(z,0))$
is of generic rank $n$,
$\Psi_1^0(0,\lambda^1,\Lambda^{\kappa_0})=0$, and $\displaystyle
\frac{\partial \widetilde \Psi_1^0}{\partial
t^1}(0,\lambda^1,\Lambda^{\kappa_0})$ is invertible for all
$(\lambda^1,\Lambda^{\kappa_0})\in {\mathcal G}^1_0(\CN)\times
{\mathcal G}^{\kappa_0}_0(\CN)$ by \eqref{e:conditions}. Hence we
may apply Proposition \ref{t:lineqn2} to conclude that there
exists a holomorphic map
\begin{equation}\label{e:gammamu}
\Gamma_1^\mu=\Gamma_1^\mu
(t^1,\lambda^1,\Lambda^{\kappa_0+k+1})\colon \C^n\times {\mathcal
G}^1_0(\CN)\times {\mathcal G}_0^{\kappa_0+k+1}(\CN)\to \C^n
\end{equation} defined in a neighbourhood of $\{0\}\times
{\mathcal G}^1_0(\CN)\times {\mathcal G}_0^{\kappa_0+k+1}(\CN)$
 which satisfies that if $X\colon (\Cn,0)\to (\Cn,0)$ is a solution
 of \eqref{e:system2} for some $(\lambda^1,\Lambda^{\kappa_0+k+1})\in
 {\mathcal G}^1_0(\CN)\times {\mathcal G}_0^{\kappa_0+k+1}(\CN)$, then
$X=X(t^1)=\Gamma_1^\mu (t^1,\lambda^1,\Lambda^{\kappa_0+k+1})$. We
set
\begin{equation}\label{e:tired1}
\widetilde
\Psi_1^\delta(t^1,\lambda^1,\Lambda^{\kappa_0+k+1}):=\Gamma_1^\mu
(t^1,\lambda^1,\Lambda^{\kappa_0+k+1}),
\end{equation}
\begin{equation}\label{e:tired2}
\widehat
\Psi_1^\delta(t^1,\lambda^1,\Lambda^{\kappa_0+k+1}):=W_\mu
\left(v^1(t^1),0,\Psi_1^0(t^1,\lambda^1,\Lambda^{\kappa_0}),0,
\Psi_1^{(k)}
(t^1,\lambda^1,\Lambda^{\kappa_0+k}),\Lambda_0^{k+1}\right),
\end{equation}
where $W_\mu$ is the $\Cd$-valued polynomial map given by Lemma
\ref{l:gderivative}. In view of Lemma \ref{l:gderivative}, the
constructed map $\Psi_1^\delta$ is the desired parametrization for
the $\delta$-th derivative of all elements of $\autM$. It remains
to check that its components belong to the space ${\mathcal
E}_n^{\kappa_0+k+1}$. In view of the induction assumption, Lemma
\ref{l:heavystuff} (i) and \eqref{e:tired2}, this fact is clear
for all the components of $\widehat \Psi_1^\delta$. To deal with
the components of $\widetilde \Psi_1^{\delta}$, we first expand
$\widetilde \Psi_1^\delta$ in a Taylor series:
\begin{equation}\label{e:form0}
\widetilde
\Psi_1^\delta(t^1,\lambda^1,\Lambda^{\kappa_0+k+1})=\Gamma_1^\mu
(t^1,\lambda^1,\Lambda^{\kappa_0+k+1})=\sum_{\gamma \in
\N^n}\Gamma_{1,\gamma}^\mu
(\lambda^1,\Lambda^{\kappa_0+k+1})(t^1)^{\gamma}. \end{equation}
Note that from Proposition \ref{t:lineqn2},
we know that each $\Gamma_{1,\gamma}^\mu
(\lambda^1,\Lambda^{\kappa_0+k+1})$ is of the following form
\begin{equation}\label{e:form}
\frac{C_{\gamma}\left[\left((\partial_{t^1}^{\nu}
\Xi_{\mu,\alpha})
%(0,\lambda^1,\Lambda^{\kappa_0+k+1})
,
(\partial_{t^1}^{\nu}\Upsilon_{\mu,\alpha})
%(0,\lambda^1,\Lambda^{\kappa_0+k+1})
,
(\partial_{t^1}^{\nu}\widetilde
\Psi_1^0)
%(0,\lambda^1,\Lambda^{\kappa_0+k+1})
\right)_{\substack{1\leq
|\alpha|\leq \kappa_0\\ |\nu|\leq s_\gamma}}\right]}{\left[{\rm
det}\, \left(\displaystyle \frac{\partial \widetilde
\Psi_1^0}{\partial
t^1}(0,\lambda^1,\Lambda^{\kappa_0})\right)\right]^{e_\gamma}},
\end{equation}
 where $C_\gamma$ is a $\Cn$-valued polynomial map of its arguments, 
 $s_\gamma$ and $e_\gamma$ are positive integers, and
 the numerator is evaluated at $(0,\lambda^1,\Lambda^{\kappa_0+k+1})$. 
 Since each component
of the $\widetilde \Psi_1^0$, $\Xi_{\mu,\alpha}$ and
$\Upsilon_{\mu,\alpha}$ for $|\alpha|\leq \kappa_0$ belongs to
${\mathcal E}_n^{\kappa_0+k+1}$ and since $\displaystyle
\frac{\partial \widetilde \Psi_1^0}{\partial
t^1}(0,\lambda^1,\Lambda^{\kappa_0})=\widetilde \lambda^1$ by
\eqref{e:conditions}, it follows clearly from \eqref{e:form0} and
\eqref{e:form} that each component of $\widetilde \Psi_1^0$
belongs to the space ${\mathcal E}_n^{\kappa_0+k+1}$. This
completes the proof of Proposition \ref{p:firststep}.

\subsection{Iteration along higher order Segre sets}\label{ss:iteration}

We now turn to constructing a parametrization for the elements of the
stability group $\autM$ along higher order Segre sets analogous to
that obtained along the first Segre set in \S \ref{ss:firstorder}.
For any positive integer $k$, we recall the notation
$t^{[k]}=(t^1,\ldots,t^k)\in \C^{nk}$ as introduced in
\eqref{e:usefulnotation}. We will prove the following.

\begin{prop}\label{p:propagation}
With the previous notation and under the assumptions of
Proposition \ref{p:firststep}, the following holds. For any $k\in
\N^*$ and any multiindex $\beta \in \N^N$, there exists a
$\CN$-valued holomorphic map
\begin{multline}
  \Psi_k^{\beta}=\Psi_k^{\beta}(t^{[k]},\lambda^1,\Lambda^{k\kappa_0+|\beta|}):=
\\  (\widetilde\Psi_k^{\beta}(t^{[k]},\lambda^1,\Lambda^{k\kappa_0+|\beta|}),\widehat\Psi_k^{\beta}(t^{[k]},\lambda^1,\Lambda^{k\kappa_0+|\beta|}))\in
\C^n\times \Cd,
\end{multline}
  whose components belong to ${\mathcal
E}_{kn}^{k\kappa_0+|\beta|}$, satisfying the following properties:
\begin{enumerate}
\item[{\rm (i)}] If we denote by $S\subset \N^N$ the set of multiindices of length one
having $0$ from the $(n+1)$-th to the $N$-th component, then for any
$(\lambda^1,\Lambda^{k\kappa_0})\in {\mathcal G}_0^1(\CN)\times
{\mathcal G}^{k\kappa_0}_0(\CN)$ we have
\begin{equation}\label{e:conditions++}
  \begin{split}
\Psi_k^0(0,\lambda^1,\Lambda^{k\kappa_0})&=0,\\ (\widetilde
\Psi_k^\beta (0,\lambda^1,\Lambda^{k\kappa_0+1}))_{\beta \in
S}&=\frac{\partial \widetilde \Psi_k^0}{\partial
t^k}(0,\lambda^1,\Lambda^{k\kappa_0})=\begin{cases} \widetilde
\lambda^1,\ {\rm if}\ k\ {\rm is\ odd}\cr \widetilde \Lambda^1,\
{\rm if}\ k\ {\rm is\ even}
\end{cases}\end{split}
\end{equation}
\item[{\rm (ii)}] For every $H\in \autM$,
\begin{equation}
((\partial^{\beta}H)\circ v^k) (t^{[k]})=
\Psi_k^{\beta}(t^{[k]},\jetm{0}{1} H,\jetm{0}{{k\kappa_0+|\beta|}} \overline{H}),\quad
{\rm if}\ k\ {\rm is}\ {\rm odd},
\end{equation}

\begin{equation}
((\partial^{\beta}H)\circ v^k) (t^{[k]})=
\Psi_k^{\beta}(t^{[k]},\jetm{0}{1} \overline{H},\jetm{0}{{k\kappa_0+|\beta|}}{H}),\quad
{\rm if}\ k\ {\rm is}\ {\rm even},
\end{equation}
for all $t^{[k]}$ sufficiently close to the origin in $\C^{kn}$.
\end{enumerate}
\end{prop}

We prove Proposition \ref{p:propagation} by induction on $k$. For
$k=1$, the result is exactly the content of Proposition
\ref{p:firststep}. Suppose now that the result holds for all $k\in
\{1,\ldots,m\}$ and let us prove the proposition for $k=m+1$. We
assume that $m+1$ is even, and leave the odd case to the reader.
We shall construct the maps $(\Psi_{m+1}^{\beta})_{\beta\in \N^N}$
by induction on $|\beta|$.

\subsubsection{Construction of the parametrization $\Psi_{m+1}^0$}
We use \eqref{e:fundamental1} in which we set
$Z=v^{m+1}=v^{m+1}\left(t^{[m+1]}\right)$ and $\zeta
=\bar{v}^{m}=\bar{v}^{m}\left(t^{[m]}\right)$. We get for every $\alpha \in
\N^n\setminus \{0\}$
\begin{equation}\label{e:morning}
\bar{Q}_{\chi^{\alpha}}(\bar{f}\circ \bar{v}^m ,H\circ v^{m+1}
)=\frac{\PP_{\alpha}(v^{m+1},\bar{v}^m,
((\partial^{\beta}\overline{H})\circ \bar{v}^m)_{1\leq |\beta|\leq
|\alpha|})
} {(\D(v^{m+1},\bar{v}^m, ((\partial^\beta \overline{H})\circ
\bar{v}^m)_{|\beta|=1})
)^{2|\alpha|-1}},
\end{equation}
for all $t^{[m+1]}$ sufficiently close to the origin in
$\C^{(m+1)n}$. In what follows, we will only use \eqref{e:morning}
for multiindices $\alpha\in \N^n\setminus \{0\}$ such that
$|\alpha|\leq \kappa_0$. For such $\alpha's$, we set
\begin{multline}\label{e:afternoon}
\varrho_\alpha (t^{[m+1]},\lambda^1,\Lambda^{(m+1)\kappa_0})
:=\\ \frac{\PP_{\alpha}\left( v^{m+1}\left(t^{[m+1]}\right),\bar{v}^m\left(t^{[m]}\right),
\overline{\Psi_m^{(|\alpha|)}}(t^{[m]},\lambda^1,\Lambda^{m\kappa_0+|\alpha|})
\right) } {\left(
\D(v^{m+1}\left(t^{[m+1]}\right),\bar{v}^m\left(t^{[m]}\right),
\overline{\Psi_m^{(1)}}(t^{[m]},\lambda^1,\Lambda^{m\kappa_0+1})
\right)^{2|\alpha|-1}}.
\end{multline}
As before, for any integer $k$ we set
\[
{\Psi_m^{(k)}}(t^{[m]},\lambda^1,\Lambda^{m\kappa_0+k})= \left(
{\Psi_m^{\beta}}(t^{[m]},\lambda^1,\Lambda^{m\kappa_0+|\beta|})
\right)_{1\leq |\beta|\leq k }\] and note that in
\eqref{e:afternoon} the bar is applied to the parametrizations
$\Psi_m^{\beta}$ viewed as functions of the variable $t^{[m]}$
only. It follows from Proposition \ref{p:reflection1} (i) and
\eqref{e:conditions++} (with $k=m$) that
\begin{equation}\label{e:night}
\begin{aligned}
\D(0,0,\overline{\Psi_m^{(1)}}(0,\lambda^1,\Lambda^{m\kappa_0+1}))
&= {\det}\,
((\overline{\Psi_m^{\beta}}(0,\lambda^1,\Lambda^{m\kappa_0+1}))_{\beta
\in S})\\
&= {\det}\, \widetilde \lambda^1.&\\
\end{aligned}
\end{equation}
This shows that for each $\alpha \in \N^n\setminus \{0\}$ with
$|\alpha|\leq \kappa_0$,
$\varrho_{\alpha}:=(\varrho^1_\alpha,\ldots,\varrho^d_{\alpha})$
defines a $\Cd$-valued holomorphic map near $\{0\} \times
{\mathcal G}^1_0(\CN)\times {\mathcal
G}_0^{m\kappa_0+1}(\CN)\subset \C^{(m+1)n}\times
{\mathcal G}^1_0(\CN)\times {\mathcal G}_0^{m\kappa_0+1}(\CN)$.
Moreover, it also follows from Lemma \ref{l:heavystuff} (ii) and
the induction assumption that the components of each
$\varrho_\alpha$ belong to the space ${\mathcal
E}_{(m+1)n}^{(m+1)\kappa_0}$.  Next, we notice that we  have
\begin{equation}\label{e:add}
H\circ v^{m+1}=v^{m+1}(f\circ v^{m+1},\bar{f}\circ
\bar{v}^m,\ldots,\bar{f}\circ \bar{v}^1). \end{equation}
 Since $M$
is assumed to be in $\mcl$, by Proposition \ref{p:classnormal}, we
may choose $n$ integers $i_1,\ldots,i_n\in \{1,\ldots,d\}$ and $n$
multiindices $\alpha^{(1)},\ldots,\alpha^{(n)}\in \N^n\setminus
\{0\}$ of length $\leq \kappa_0$ so that
$z\mapsto (\bar{q}^{\, i_1}_{\alpha^{(1)}}(z,0),\ldots,\bar{q}^{\, i_n}_{\alpha^{(n)}}(z,0))$
is of generic rank $n$. Consider the following singular analytic
system in the unknowns $(T^1,\ldots,T^{m+1})\in \C^n\times \C^n
\ldots \times \C^n= \C^{(m+1)n}$:
\begin{equation}\label{e:system3}
\begin{cases}
  \begin{aligned}
\bar{Q}^{i_\nu}_{\chi^{\alpha^{(\nu)}}}(T^m,v^{m+1}(T^1,&\ldots,T^{m+1}))=
\varrho^{i_\nu}_{\alpha^{(\nu)}}
(t^{[m+1]},\lambda^1,\Lambda^{(m+1)\kappa_0}),&& \nu=1,\ldots,n \\

T^k&=\overline{\widetilde
\Psi_k^0}(t^{[k]},\lambda^1,\Lambda^{k\kappa_0}), &&\hspace{-2cm}k\in
\{1,\ldots,m\},\ k\ {\rm odd}\\
T^k&=\widetilde
\Psi_k^0(t^{[k]},\lambda^1,\Lambda^{k\kappa_0}), &&\hspace{-2cm}k\in
\{1,\ldots,m\},\ k\ {\rm even}
\end{aligned}
\end{cases}
\end{equation}
Note that by \eqref{e:morning}, \eqref{e:afternoon}, \eqref{e:add}
and  the induction assumption, for any $H=(f,g)\in \autM$, the map
\begin{equation}\label{e:vartheta}
(T^1,T^2,\ldots,T^m,T^{m+1})=\vartheta_f
\left(t^{[m+1]}\right):=(\bar{f}\circ \bar{v}^1,f\circ
v^2,\ldots,\bar{f}\circ \bar{v}^m,f\circ v^{m+1})
\end{equation}
is a solution of the system \eqref{e:system3} with
$\lambda^1:=\jetm{0}{1}\overline{H}$ and
$\Lambda^{(m+1)\kappa_0}:=j_0^{(m+1)\kappa_0}H$. Moreover, it
follows from the normality of the chosen coordinates $(z,w)$ that
the Jacobian matrix of the map $\vartheta_f$ is the following
$(m+1)n\times (m+1)n$ block triagular matrix
\begin{equation}\label{e:matrix}
\left(\begin{array}{ccccc} \frac{\partial \bar{f}}{\partial
\chi}(0)&0&0&\ldots&0\\
0&\frac{\partial {f}}{\partial z}(0)&0&\ldots&0\\
0&0& \frac{\partial \bar{f}}{\partial \chi}(0)&0&0\\
\vdots&\vdots&0&\ddots&0\\
0&0&0&0&\frac{\partial {f}}{\partial z}(0)
\end{array} \right)
\end{equation}
Hence the map $\vartheta_f$ is a local biholomorphism of
$\C^{(m+1)n}$ fixing the origin. In order to apply Corollary
\ref{c:corparcn1} so that to get a parametrization of the
invertible solutions of the system \eqref{e:system3} we need to
check that the holomorphic map
$$(T^1,\ldots,T^m,T^{m+1})\mapsto \left(T^1,\ldots,T^m,
\left(\bar{Q}^{i_\nu}_{\chi^{\alpha^{(\nu)}}}(T^m,v^{m+1}(T^1,\ldots,T^{m+1}))\right)_{\nu=1,\ldots,n}\right)$$
is of generic rank $(m+1)n$. Since the Jacobian of this system with respect to $(T^1,\dots,T^{m+1})$ is
equal to the the Jacobian of the map
$$(T^1,\ldots,T^m,T^{m+1})\mapsto
\left(\bar{Q}^{i_\nu}_{\chi^{\alpha^{(\nu)}}}(T^m,v^{m+1}(T^1,\ldots,T^{m+1}))\right)_{\nu=1,\ldots,n}$$
with respect to $T^{m+1}$, it is enough to check that this map
is of generic rank $n$ for generic $T^1,\dots,T^m$ (or, equivalently, is of generic rank $n$ for one particular point $T^1_0,\dots T^{m}_0$).
This is indeed the case since, by the
above choice, the holomorphic map
$$T^{m+1}\mapsto
\left(\bar{Q}^{i_\nu}_{\chi^{\alpha^{(\nu)}}}(0,T^{m+1},0)\right)_{\nu=1,\ldots,n}
= \left( \bar{q}^{\, i_1}_{\alpha^{(1)}}(T^{m+1},0),\ldots,\bar{q}^{\, i_n}_{\alpha^{(n)}}(T^{m+1},0) \right)$$
is of generic rank $n$. We may, therefore, apply Corollary
\ref{c:corparcn1} to conclude that there exists a holomorphic map
\begin{equation}\label{e:gamma+}
\Gamma_{m+1}^0\colon \C^{(m+1)n}\times \glc \times
{\mathcal G}_0^1(\CN)\times {\mathcal G}_0^{(m+1)\kappa_0}(\CN) \to \C^{(m+1)n}
\end{equation} defined in a neighbourhood of $\{0\}\times \glc \times {\mathcal G}_0^1(\CN)\times {\mathcal G}_0^{(m+1)\kappa_0}(\CN)$, satisfying
\begin{equation}\label{e:satisfaction}
\Gamma_{m+1}^0 (0,A,\lambda^1,\Lambda^{(m+1)\kappa_0})=0,\quad
\frac{\partial \Gamma_{m+1}^0}{\partial
t^{[m+1]}}(0,A,\lambda^1,\Lambda^{(m+1)\kappa_0})=A \end{equation}
 for all
$(A,\lambda^1,\Lambda^{(m+1)\kappa_0})\in \glc \times {\mathcal
G}_0^1(\CN)\times {\mathcal G}_0^{(m+1)\kappa_0}(\CN)$
 such that if $u\colon (\C^{(m+1)n},0)\to (\C^{(m+1)n},0)$ is a germ of
a biholomorphism satisfying \eqref{e:system3} for some
$(\lambda^1,\Lambda^{(m+1)\kappa_0})\in {\mathcal
G}^1_0(\CN)\times {\mathcal G}_0^{\kappa_0}(\CN)$, then
\[u\left(t^{[m+1]}\right)=\Gamma_{m+1}^0
(t^{[m+1]},\jetm{0}{1} u,\lambda^1,\Lambda^{(m+1)\kappa_0}).\] Consider
the following map $\varpi\colon {\mathcal G}^1_0(\CN)\times
{\mathcal G}^1_0(\CN)\to \glc$ given by
\begin{equation}\label{e:setup0}
(\lambda^1,\Lambda^1)\mapsto \varpi(\lambda^1,\Lambda^1):=
\left(\begin{array}{ccccc} \widetilde \lambda^1&0&0&\ldots&0\\
0&\widetilde \Lambda^1&0&\ldots&0\\
0&0& \widetilde \lambda^1&0&0\\
\vdots&\vdots&0&\ddots&0\\
0&0&0&0&\widetilde \Lambda^1
\end{array} \right).
\end{equation}
We write
$\Gamma_{m+1}^0=(\Gamma_{m+1}^{0,1},\ldots,\Gamma_{m+1}^{0,m+1})\in
\C^n\times \C^n \ldots \times \C^n=\C^{(m+1)n}$ and set
\begin{equation}\label{e:setup1}
\widetilde
\Psi_{m+1}^0(t^{[m+1]},\lambda^1,\Lambda^{(m+1)\kappa_0}):=
\Gamma_{m+1}^{0,m+1}(t^{[m+1]},\varpi(\lambda^1,\Lambda^1),\lambda^1,\Lambda^{(m+1)\kappa_0}).
\end{equation}
Since for every $(\lambda^1,\Lambda^{(m+1)\kappa_0})\in
{\mathcal G}_0^{1}(\CN)\times {\mathcal G}_0^{(m+1)\kappa_0}$,
$\widetilde\Psi_{m+1}^0(0,\lambda^1,\Lambda^{(m+1)\kappa_0})=0$ and ${\Psi_{m}^0}(0,\lambda^1,\Lambda^{m\kappa_0})=0$
by \eqref{e:satisfaction}, \eqref{e:setup1} and the induction
assumption, we may also set
\begin{multline}\label{e:setup2}
\widehat
\Psi_{m+1}^0(t^{[m+1]},\lambda^1,\Lambda^{(m+1)\kappa_0}):=\\Q(\widetilde
\Psi_{m+1}^0(t^{[m+1]},\lambda^1,\Lambda^{(m+1)\kappa_0}),\overline{\Psi_{m}^0}(t^{[m]},\lambda^1,\Lambda^{m\kappa_0})).
\end{multline}
Then $\Psi^0_{m+1}=(\widetilde \Psi^0_{m+1},\widehat
\Psi^0_{m+1})$ defines a $\CN$-valued holomorphic map in a
neighbourhood of $\{0\}\times {\mathcal G}_0^{1}(\CN)\times
{\mathcal G}_0^{(m+1)\kappa_0}(\CN)$ in $\C^{(m+1)n}\times
{\mathcal G}_0^{1}(\CN)\times {\mathcal G}_0^{(m+1)\kappa_0}(\CN)$
and satisfies
$\Psi_{m+1}^0(0,\lambda^1,\Lambda^{(m+1)\kappa_0})=0$ for all
$(\lambda^1,\Lambda^{(m+1)\kappa_0})\in {\mathcal
G}_0^{1}(\CN)\times {\mathcal G}_0^{(m+1)\kappa_0}(\CN)$.
Furthermore, in view of \eqref{e:matrix}, \eqref{e:setup0} and the
above construction, for any $H=(f,g)\in \autM$, we have that for
$t^{[m+1]}$ sufficiently small in $\C^{(m+1)n}$,
\begin{equation}
(H\circ
v^{m+1})\left(t^{[m+1]}\right)=\Psi_{m+1}^0(t^{[m+1]},\jetm{0}{1}\overline{H},j_0^{(m+1)\kappa_0}H).
\end{equation}
From \eqref{e:satisfaction}, we have $\displaystyle \frac{\partial
\Gamma_{m+1}^0}{\partial
t^{[m+1]}}(0,\varpi(\lambda^1,\Lambda^1),\lambda^1,\Lambda^{(m+1)\kappa_0})=\varpi(\lambda^1,\Lambda^1)$
and hence from \eqref{e:setup0} and \eqref{e:setup1}, we get that
$$\frac{\partial
\widetilde \Psi_{m+1}^0}{\partial
t^{m+1}}(0,\lambda^1,\Lambda^{(m+1)\kappa_0})=\widetilde
\Lambda^1$$ for all $(\lambda^1,\Lambda^{(m+1)\kappa_0})\in
{\mathcal G}^1_0(\CN)\times {\mathcal G}^{(m+1)\kappa_0}_0(\CN)$.
This proves part of \eqref{e:conditions++}. To finish the case
$\beta=0$, we need to check that each of the components of the
constructed parametrization $\Psi_{m+1}^0$ belong to the space
${\mathcal E}^{(m+1)\kappa_0}_{(m+1)n}$. It is actually enough to
prove that each component of $\widetilde \Psi_{m+1}^0$ is in
${\mathcal E}^{(m+1)\kappa_0}_{(m+1)n}$ by \eqref{e:setup2}, Lemma
\ref{l:heavystuff} (i) and the induction assumption. Consider the
map $\Gamma_{m+1}^0$ given in \eqref{e:gamma+} and obtained from
Corollary \ref{c:corparcn1} applied to the system \eqref{e:system3}.
Then from that result, we know that if we write the Taylor
expansion
$$\Gamma_{m+1}^0(t^{[m+1]},A,\lambda^1,\Lambda^{(m+1)\kappa_0}):=\sum_{\gamma \in \N^{(m+1)n}}
\Gamma^0_{m+1,\gamma}(A,\lambda^1,\Lambda^{(m+1)\kappa_0})\left(t^{[m+1]}\right)^{\gamma},$$
then each $
\Gamma^0_{m+1,\gamma}(A,\lambda^1,\Lambda^{(m+1)\kappa_0})$ can be
written in the following form
\begin{equation}\label{e:component}
\frac{\Omega_\gamma (A,\lambda^1,\Lambda^{(m+1)\kappa_0})}{({\rm
det}\ A)^{e_\gamma}} \end{equation}
 where $e_\gamma$ is a nonnegative
integer and $\Omega_\gamma (A,\lambda^1,\Lambda^{(m+1)\kappa_0})$
is a polynomial map in the following arguments
\begin{equation}\label{e:argument1}
\left((\partial^{s_\gamma}_{t^{[m+1]}}\varrho_{\alpha})
(0,\lambda^1,\Lambda^{(m+1)\kappa_0})\right)_{|\alpha|\leq
\kappa_0}
\end{equation}
\begin{equation}\label{e:argument2}
\left((\partial^{s_\gamma}_{t^{[m+1]}}\Psi^0_k)(0,\lambda^1,\Lambda^{m\kappa_0})\right)_{1\leq
k\leq m}
\end{equation}
\begin{equation}\label{e:argument3}
\left((\partial^{s_\gamma}_{t^{[m+1]}}\overline{\Psi^0_k})(0,\lambda^1,\Lambda^{m\kappa_0})\right)_{1\leq
k\leq m} \end{equation}
 and in $A$. Here $s_\gamma$ denotes some positive integer. By the
 induction assumption and the discussion after \eqref{e:night}, each
 component of each term in \eqref{e:argument1}, \eqref{e:argument2} and in
 \eqref{e:argument3} can be written as a polynomial in $\lambda^1$
 and $\Lambda^{m\kappa_0}$ over a product of powers of ${\det}\, \widetilde
 \lambda^1$ and ${\det}\, \widetilde \Lambda^{1}$. Now,
 if we write the Taylor expansion
 $$\widetilde \Psi_{m+1}^{0}(t^{[m+1]},\lambda^1,\Lambda^{(m+1)\kappa_0}):=\sum_{\gamma \in \N^{(m+1)n}}
\widetilde
\Psi^0_{m+1,\gamma}(\lambda^1,\Lambda^{(m+1)\kappa_0})\left(t^{[m+1]}\right)^{\gamma},$$
then in view of \eqref{e:setup1}, each component of $\widetilde
\Psi^0_{m+1,\gamma}(\lambda^1,\Lambda^{(m+1)\kappa_0})$ is a
component of \eqref{e:component} with $A=\varpi(\lambda^1,
\Lambda^1)$. From \eqref{e:setup0} and the above discussion, we
see that each component $\widetilde \Psi^0_{m+1}$ indeed belongs
to the space ${\mathcal E}^{(m+1)\kappa_0}_{(m+1)n}$. This
completes the proof of Proposition \ref{p:propagation} in the case
$\beta=0$.

\subsubsection{Induction for the construction of the parametrizations $(\Psi^{\beta}_{m+1})_{\beta \in \N^N}$.}
We assume here that the desired parametrizations
$\Psi_{j}^{\beta}$ for all $j\in
\{1,\ldots,m\}$ and all $\beta\in \N^N$ as well as the
parametrizations $\Psi_{m+1}^{\beta}$ for $|\beta|\leq r$ have already been constructed.
We
will construct the parametrization $\Psi^{\delta}_{m+1}$ for
$\delta\in \N^N$ of length $r+1$ with the desired properties.
We first reduce this to the construction for $\delta$ corresponding to purely transversal derivatives
in Lemma
\ref{l:course},
and start with the following computational lemma.

\begin{lem}\label{l:useful} With the above notation and assumptions,
the following holds.
\begin{enumerate}
\item[{\rm (i)}] For any $H\in \autM$, we have
\begin{multline}\label{e:firstorderderiv}
\frac{\partial \Psi_{m+1}^0}{\partial
t^{m+1}}(t^{[m+1]},\jetm{0}{1}\overline{H},j_0^{(m+1)\kappa_0}H)=\\
\left(\frac{\partial H}{\partial z}\circ
v^{m+1}\right)
\left(t^{[m+1]}\right)+ \left(\frac{\partial H}{\partial
w}\circ v^{m+1}\right)\left(t^{[m+1]}\right)\cdot
Q_z\left(t^{m+1},\bar{v}^{m}\left(t^{[m]}\right)\right),
\end{multline}
for all $t^{[m+1]}$ sufficiently close to the origin in
$\C^{(m+1)n}$.
\item[{\rm (ii)}] Let
  \[I_r=\left\{ \mu \in \N^N:|\mu|=r\right\}, \quad J_r=\left\{ \nu \in \N^d:|\nu|=r+1\right\}.\]
  For every multindex $\delta \in \N^N$ of length $r+1$,
there exist holomorphic $n\times 1$ matrices
$(h_{\mu,\delta})_{\mu\in I_r}$
and holomorphic functions $(\xi_{\nu,\delta})_{\nu \in J_r}$ defined near the origin in
$\C^{(m+1)n}$, such that for any $H\in \autM$,
\begin{equation}\label{e:nochoice}\begin{split}
\left((\partial^{\delta}H)\circ v^{m+1}\right)&\left(t^{[m+1]}\right)=\\
&\sum_{\mu \in I_r} 
\frac{\partial \Psi^{\mu}_{m+1}}{\partial
t^{m+1}}\left(t^{[m+1]},\jetm{0}{1}
\overline{H},j_0^{(m+1)\kappa_0+r}H\right)\cdot
h_{\mu,\delta}\left(t^{[m+1]}\right)
\\&+\sum_{\nu\in J_r
}\xi_{\nu,\delta}\left(t^{[m+1]}\right)
\left(\frac{\partial^{|\nu|}H}{\partial w^{\nu}}\circ v^m
\right)\left(t^{[m+1]}\right).\end{split}
\end{equation}
\end{enumerate}
\end{lem}

\begin{proof}
Part (i) of Lemma \ref{l:useful} follows by
differentiating the identity
\begin{align*}
   \left(H\circ v^{m+1}\right)\left(t^{[m+1]}\right)&=
   H\left( t^{m+1},Q \left( t^{m+1},\bar{v}^m \left( t^{[m]} \right) \right) \right) \\
   &=\Psi_{m+1}^0 \left( t^{[m+1]},\jetm{0}{1}\overline{H},\jetm{0}{(m+1)\kappa_0} H \right),
\end{align*}
which holds for all $H\in \autM$,
with respect to $t^{m+1}$.

To prove (ii) for every $\delta=(\delta_1,\ldots,\delta_N) \in
\N^N$ of length $r$, we set $q_\delta:=\sum_{i=1}^n\delta_i$, and
we prove \eqref{e:nochoice} by induction on $q_{\delta}$. If
$q_\delta=0$, \eqref{e:nochoice} is a trivial statement. Suppose
now that $q_\delta>0$ with $\delta \in \N^N$ of length $r+1$. Then
we may write ${\delta}=e+\omega$ with $e\in S$ ($S$ denotes the
set of multiindices of length $1$ having $0$ from the $n+1$-th
component to the $N$-th component, as defined in Proposition
\ref{p:propagation} (i)), $q_\omega<q_\delta$ and $|\omega|=r$.
Differentiating the identity
\begin{align*}\Psi^{\omega}_{m+1}(t^{[m+1]},\jetm{0}{1}\overline{H},j_0^{(m+1)\kappa_0+r}H)&=
\left((\partial^{\omega}H)\circ v^{m+1}\right)\left(t^{[m+1]}\right)\\
&= \left((\partial^{\omega}H)(t^{m+1},Q(t^{m+1},\bar{v}^m(t^{[m]})\right)\\
\end{align*} which holds for all $H\in \autM$ with
respect to $t^{m+1}$, we get (similarly to what was done for the
proof of (i)) that
\begin{multline}\label{e:sick}
 \frac{\partial
\Psi^{\omega}_{m+1}}{\partial
t^{m+1}}\left(t^{[m+1]},\jetm{0}{1}\overline{H},j_0^{(m+1)\kappa_0+r}H\right)=
\left(\left(\frac{\partial}{\partial z}(\partial^\omega H)\right)\circ
v^{m+1}\right)\left(t^{[m+1]}\right)\\+Q_z\left(t^{m+1},\bar{v}^{m}\left(t^{[m]}\right) \right)\cdot
\left(\left(\frac{\partial}{\partial w}(\partial^\omega H)\right)\circ
v^{m+1}\right)\left(t^{[m+1]}\right).\end{multline}
 Hence \eqref{e:sick} and the induction hypothesis shows that $(\partial^{\delta}H)\circ v^{m+1}$
 may be written in the form \eqref{e:nochoice}. The proof of the
 lemma is complete.
\end{proof}

We  let $O_{r+1}$ denote the set of multindices $\delta\in \N^N$
of length $r+1$ with $0$ from the first to the $n$-th component.
Note that the set of multiindices in $\N^d$ of length $r+1$ can be
identified with $O_{r+1}$ via the map $\N^d \ni \nu \mapsto
(0,\nu) \in O_{r+1}$. The purpose of Lemma \ref{l:course} is to
show that it is enough to find the parametrizations
$\Psi^{\delta}_{m+1}$ with $\delta \in O_{r+1}$ to get to the
parametrizations $\Psi_{m+1}^{\delta}$ for all $\delta \in \N^N$
of length $r+1$.

\begin{lem}\label{l:course}
Assume that one  has constructed the parametrizations
$(\Psi^\beta_{m+1})$ for all multiindices $\beta \in O_{r+1}$
satisfying the required properties of Proposition {\rm
\ref{p:propagation}}. Then one may construct all the
parametrizations $(\Psi^\delta_{m+1})$ for all multiindices
$\delta \in \N^N$ of length $r+1$ satisfying the required
properties of Proposition {\rm \ref{p:propagation}}.
\end{lem}

\begin{proof}  If $r>0$, for any multiindex $\delta\in
  \N^N\setminus O_{r+1}$, we set (writing $t=t^{[m+1]}$)
%\begin{equation}\label{e:k}
%\begin{aligned}
%\Psi^{\delta}_{m+1}(t^{[m+1]},\lambda^1,\Lambda^{(m+1)\kappa_0+r+1}):&=&\sum_{\mu
%\in I_r} \frac{\partial \Psi^{\mu}_{m+1}}{\partial
%t^{m+1}}(t^{[m+1]},\lambda^1,\Lambda^{(m+1)\kappa_0+r})\cdot
%h_{\mu,\delta}\left(t^{[m+1]}\right)\\
%&+& \sum_{\nu\in J_r }\xi_{\nu,\delta}\left(t^{[m+1]}\right)\cdot
%\Psi^{(0,\nu)}_{m+1}(t^{[m+1]},\lambda^1,\Lambda^{(m+1)\kappa_0+r+1}).
%\end{aligned}
%\end{equation}
\begin{multline}
\label{e:k}
\Psi^{\delta}_{m+1}(t,\lambda^1,\Lambda^{(m+1)\kappa_0+r+1}):=\\ \sum_{\mu
\in I_r} \frac{\partial \Psi^{\mu}_{m+1}}{\partial
t^{m+1}}(t,\lambda^1,\Lambda^{(m+1)\kappa_0+r})\cdot
h_{\mu,\delta}\left(\cdot\right)
+ \sum_{\nu\in J_r }\xi_{\nu,\delta}\left(t\right)\cdot
\Psi^{(0,\nu)}_{m+1}(t,\lambda^1,\Lambda^{(m+1)\kappa_0+r+1}).
\end{multline}
Then it follows from Lemma \ref{l:useful} (ii) that the desired
parametrizations are constructed for all $\delta$ of length $r+1$.
Moreover, it is clear that since each component of
$\Psi_{m+1}^\beta$ for $\beta\in O_{r+1}$ or for $|\beta|=r$
belongs to the space ${\mathcal E}_{(m+1)n}^{(m+1)\kappa_0+r+1}$
by assumption, the same fact also holds for each component of
$\Psi^{\delta}_{m+1}(t^{[m+1]},\lambda^1,\Lambda^{(m+1)\kappa_0+r+1})$
with $\delta \not \in O_{r+1}$.

In case $r=0$, note that $\N^d\setminus O_{r+1}$ coincides with
$S$ as defined in Proposition \ref{p:propagation} (i). If the
parametrizations $\Psi_{m+1}^\delta$ for $\delta \in O_{r+1}$ have
been constructed, then we may set
\begin{multline}\label{e:crap}
\left( \Psi_{m+1}^\delta\left(t^{[m+1]},\lambda^1,\Lambda^{(m+1)\kappa_0+1}\right)\right)_{\delta
\in S}:=\frac{\partial \Psi^0_{m+1}}{\partial
t^{m+1}}(t^{[m+1]},\lambda^1,\Lambda^{(m+1)\kappa_0})\\
-\left( \Psi_{m+1}^\delta \left(
t^{[m+1]},\lambda^1,\Lambda^{(m+1)\kappa_0+1}
\right)\right)_{\delta \in O_{r+1}}\cdot Q_z\left(
t^{m+1},\bar{v}^{m}\left(t^{[m]}\right)\right), \end{multline}
where we see $(\Psi_{m+1}^\delta)_{\delta \in S}$ as a $N\times n$
matrix and $(\Psi_{m+1}^\delta)_{\delta \in O_{r+1}}$
 as a $N\times d$ matrix. It follows from Lemma \ref{l:useful} (i) that the
parametrizations are constructed for all $\delta$ of length $1$.
It also follows from \eqref{e:crap} that our construction gives
$$(\widetilde \Psi_{m+1}^\delta
(0,\lambda^1,\Lambda^{(m+1)\kappa_0}))_{\delta \in
S}=\frac{\partial \widetilde \Psi_{m+1}^0}{\partial
t^{m+1}}(0,\lambda^1,\Lambda^{(m+1)\kappa_0}),$$ which proves the
required condition in \eqref{e:conditions++}. Finally, 
since each component of $\Psi_{m+1}^\mu$ for $\mu=0$ or for
$\mu\in O_{r+1}$ belongs to the space ${\mathcal
E}_{(m+1)n}^{(m+1)\kappa_0+1}$ by assumption, the same fact also
holds for each component of
$\Psi^{\delta}_{m+1}(t^{[m+1]},\lambda^1,\Lambda^{(m+1)\kappa_0+1})$
with $\delta \in S$ in view of \eqref{e:crap}. This completes the
proof of Lemma \ref{l:course}.
\end{proof}

By Lemma \ref{l:course}, it is now enough to construct the
parametrization $\Psi_{m+1}^{\delta}$ for $\delta \in \N^N$ with
$|\delta|=r+1$ and $\delta \in O_{r+1}$. In what follows, we fix
such a multiindex $\delta$, and we, therefore, write $\delta=(0,\mu)$
with $\mu\in \N^d$ and $|\mu|=r+1$.

We start by using \eqref{e:fundamental3} in which we set
$Z=v^{m+1}=v^{m+1}\left(t^{[m+1]}\right)$ and
$\zeta=\bar{v}^m=\bar{v}^{m}\left(t^{[m]}\right)$. We obtain that for every
$\alpha \in \N^n\setminus \{0\}$ and for every $H=(f,g)\in \autM$
the expression
\begin{multline}\label{e:follow0}
\left(\frac{\partial^{|\mu|}f}{\partial w^\mu}\circ
v^{m+1}\right)\cdot \\
\left(\bar{Q}_{\chi^{\alpha},z}(\bar{f}\circ
\bar{v}^m ,H\circ v^{m+1})+ Q_z(f\circ v^{m+1},\overline{H}\circ
\bar{v}^m)\cdot \bar{Q}_{\chi^{\alpha},w}(\bar{f}\circ \bar{v}^m
,H\circ v^{m+1})\right)
\end{multline}
is equal to \eqref{e:follow1}+\eqref{e:follow00} where
\begin{equation}\label{e:follow1}
{\mathcal T}_{\mu,\alpha}'\left(v^{m+1},\bar{v}^m,H\circ
v^{m+1},\overline{H}\circ \bar{v}^m, ((\partial^{\beta} H)\circ
v^{m+1})_{1\leq |\beta|\leq r},
%\jetm{v^{m+1}}{r} H,
((\partial^{\beta}\overline{H})\circ \bar{v}^m)_{1\leq |\beta|\leq
r+1}
%\jetm{\bar{v}^{m}}{r} \overline{H}
\right) \end{equation} \begin{equation}\label{e:follow00}
\frac{{\mathcal Q}_{\mu,\alpha}\left( v^{m+1},\bar{v}^m,
((\partial^{\beta}\overline{H})\circ \bar{v}^m)_{1\leq |\beta|\leq
|\alpha|+r+1}
%\jetm{\bar{v}^{m}}{|\alpha|+r+1} \overline{H}
\right) }{ \left(\D(v^{m+1},\bar{v}^m, ((\partial^\beta
\overline{H})\circ \bar{v}^m)_{|\beta|=1})
%\jetm{\bar{v}^{m}}{1} \overline{H})
\right)^{2|\alpha|+r}}.
\end{equation}
As in all previous analogous situations, we are only
interested in multiindices $\alpha \in \N^n\setminus \{0\}$ such
that $|\alpha|\leq \kappa_0$. For such $\alpha$'s, we define the
following $n\times d$ matrix with holomorphic coefficients near
the origin in $\C^{(m+1)n}$:
\begin{multline}\label{e:firstguy}
K_\alpha(T^1,T^2,\ldots,T^{m+1}):= \bar{Q}_{\chi^{\alpha},z}(T^m
,v^{m+1}(T^1,\ldots,T^{m+1}))+\\
Q_z(T^{m+1},\bar{v}^m(T^1,\ldots,T^m))\cdot
\bar{Q}_{\chi^{\alpha},w}(T^m ,v^{m+1}(T^1,\ldots,T^{m+1})).
\end{multline}
Here each $T^j\in \C^n$, $j=1,\ldots,m+1$. We also write
$K_\alpha:=(K_\alpha^1,\ldots,K_\alpha^d)$ and set
\begin{multline}\label{e:secondguy}
\phi_{\mu,\alpha}(t^{[m+1]},\lambda^1,\Lambda^{(m+1)\kappa_0+r+1}):=\\
{\mathcal T}_{\mu,\alpha}'\left(v^{m+1},\bar{v}^m,\Psi_{m+1}^0,
 \overline{\Psi_m^0},
 %(\Psi_{m+1}^{\beta})_{1\leq|\nu|\leq r},
 \Psi_{m+1}^{(r)}
 %(\overline{\Psi_m^{\beta}})_{1\leq |\beta|\leq r+1}
 \overline{\Psi_m^{(r+1)}}
 \right)+
\frac{{\mathcal Q}_{\mu,\alpha}\left(v^{m+1},\bar{v}^m,
%(\overline{\Psi_{m}^{\beta}})_{|\beta|\leq |\alpha|+r+1})
\overline{\Psi_{m}^{(r+1)}}\right) }{ \left(\D(v^{m+1},\bar{v}^m,
%(\overline{\Psi^\beta_m})_{|\beta|=1}))^{2|\alpha|+r}}
\overline{\Psi^{(1)}_m})\right)^{2|\alpha|+r}}
\end{multline}
By the induction assumption, \eqref{e:night} and Lemma
\ref{l:heavystuff},
$\phi_{\mu,\alpha}:=(\phi^1_{\mu,\alpha},\ldots,\phi^d_{\mu,\alpha})$
defines a $\Cd$-valued holomorphic map defined in a neighbourhood
of $\{0\} \times {\mathcal G}^1_0(\CN)\times {\mathcal
G}_0^{(m+1)\kappa_0+r+1}(\CN)$ in $\C^{(m+1)n}\times {\mathcal
G}_0^1(\CN)\times {\mathcal G}_0^{(m+1)\kappa_0+r+1}(\CN)$, whose
components belong to ${\mathcal E}_{(m+1)n}^{(m+1)\kappa_0+r+1}$.
Since $M$ is assumed to be in $\mcl$, we may again, by Proposition
\ref{p:classnormal}, choose $n$ integers $i_1,\ldots,i_n\in
\{1,\ldots,d\}$ and $n$ multiindices
$\alpha^{(1)},\ldots,\alpha^{(n)}\in \N^n\setminus \{0\}$ of
length $\leq \kappa_0$ so that $z\mapsto
(\bar{q}^{i_1}_{\alpha^{(1)}}(z,0),\ldots,\bar{q}^{i_n}_{\alpha^{(n)}}(z,0))$
is of generic rank $n$. In what follows, we write
$Y^k=(Y^k_1,\ldots,Y^k_n)\in \Cn$, $k\in \{1,\ldots, m\}$ and for
$j=1,\ldots,n$ we denote by $\widetilde \Psi^{0,j}_{k}$ (resp.\
$\overline{\widetilde \Psi^{0,j}_{k}}$) the $j$-th component of
$\widetilde \Psi^{0}_{k}$ (resp.\ $\overline{\widetilde
\Psi^{0}_{k}}$). Consider the following linear singular analytic
system in the unknowns $(X,Y^1,\ldots,Y^{m})\in \C^{n}\times \Cn
\ldots \times \Cn =\C^{(m+1)n}$ given by
\begin{equation}\label{e:system4}
\begin{cases}
\begin{aligned}
X\cdot
K^{i_\nu}_{\alpha^{(\nu)}}&(T^1,\ldots,T^{m+1})=\phi^{i_{\nu}}_{\mu,\alpha^{(\nu)}}
(t^{[m+1]},\lambda^1,\Lambda^{(m+1)\kappa_0+r+1}), \quad 1\leq \nu \leq n \\
&Y_j^k\cdot T_j^k=\overline
{\widetilde \Psi_k^{0,j}}
(t^{[k]},\lambda^1,\Lambda^{k\kappa_0}),
\quad 1\leq j\leq
n,\ 1\leq k\leq m,\ k\ {\rm odd} \\
& Y_j^k\cdot
T_j^k=\widetilde
\Psi_k^{0,j}(t^{[k]},\lambda^1,\Lambda^{k\kappa_0}),\quad 1\leq j\leq
n,\ 1\leq k\leq m,\ k\ {\rm even}
\end{aligned}
\end{cases}
\end{equation}
in which for every $k\in \{1,\ldots,m+1\}$ we have set
\begin{equation}\label{e:system5}
T^k:=\begin{cases} \widetilde
\Psi_k^0(t^{[k]},\lambda^1,\Lambda^{k\kappa_0})\ {\rm if}\  k\
{\rm is}\ {\rm even}\cr \overline{\widetilde
\Psi_k^0}(t^{[k]},\lambda^1,\Lambda^{k\kappa_0})\ {\rm if}\  k\
{\rm is}\ {\rm odd}
\end{cases}
 \end{equation}
In view of \eqref{e:follow0}, \eqref{e:follow1},
 \eqref{e:add},
\eqref{e:firstguy},\eqref{e:secondguy} and the induction
assumption, we know that the vector
$((\partial_{w}^{\mu}f)\circ v^{m+1} ,1,1\ldots,1)$ is a solution
of the system \eqref{e:system4} with
$\lambda^1=\jetm{0}{1}\overline{H}$ and
$\Lambda^{(m+1)\kappa_0+r+1}=j_0^{(m+1)\kappa_0+r+1}H$ 
for any $H\in \autM$. To apply
Proposition \ref{t:lineqn2} to get a parametrization of the
solutions of the system \eqref{e:system4}-\eqref{e:system5}, we
need to check that the $(m+1)n\times (m+1)n$ matrix associated to
this linear system has maximal generic rank and that for all
$(\lambda^1,\Lambda^{(m+1)\kappa_0+r+1})\in {\mathcal
G}_0^1(\CN)\times {\mathcal G}_0^{(m+1)\kappa_0+r+1}(\CN)$, the
mapping ${\mathcal I}\colon (\C^{(m+1)n},0)\to (\C^{(m+1)n},0)$
given by
\begin{equation}\label{e:endsoon?}
 t^{[m+1]}\mapsto \left(\overline{\widetilde
\Psi_1^0}(t^{[1]},\cdot),\widetilde
\Psi_2^0(t^{[2]},\cdot),\ldots,\overline{\widetilde
\Psi_m^{0}}(t^{[m]},\cdot),\widetilde\Psi_{m+1}^0(t^{[m+1]},\cdot)\right)
\end{equation}
is a local biholomorphism at the origin. Let us check the first
condition.  The  $(m+1)n\times (m+1)n$ matrix associated to the
 linear singular system \eqref{e:system4} is given by the following
\begin{equation}\label{e:matrixmore}
\left(\begin{array}{cccccccc}
K(T^1,\ldots,T^{m+1})&0&\ldots&\ldots&\ldots&\ldots&\ldots&0\\
0&T^1_1&0&\ldots&\ldots&\ldots&\ldots&0\\
\vdots&0&\ddots&0&\ldots&\ldots&\ldots&0\\
\vdots&\vdots&0&T^1_n&0&\ldots&\ldots&0\\
\vdots&\vdots&\ldots&0&T^2_1&0&\ldots&0\\
\vdots&\vdots&\ldots&\ldots&0&\ddots&0&0\\
\vdots&\vdots&\ldots&\ldots&\ldots&0&\ddots&0\\
0&0&\ldots&\ldots&\ldots&\ldots&0&T^m_n
\end{array} \right)
\end{equation}
where $T^k=(T_1^k,\ldots,T_n^k)$ for $k=1,\ldots,m$ and $K$ is the
$n\times n$ matrix whose columns are given by the vectors
$K^{i_\nu}_{\alpha^{(\nu)}}$, $\nu=1,\ldots,n$. Hence it is enough
to show that $K$ has generic rank $n$. In fact, we claim that $K$
restricted to the subspace $T^1=\ldots=T^m=0$ has already generic
rank $n$. Indeed, in view of \eqref{e:firstguy} and the normality
of the chosen coordinates $(z,w)$, $K(0,\ldots,0,T^{m+1})$ is the
$n\times n$ matrix whose columns are the vectors
$\alpha^{(\nu)}!\, \displaystyle \frac{\partial
\bar{q}^{i_\nu}_{\alpha^{(\nu)}}}{\partial z}(T^{m+1},0)$ for
$\nu=1,\ldots,n$, which proves the above claim and the fact that
the matrix \eqref{e:matrixmore} has generic rank $(m+1)n$. We now
prove that the local holomorphic map given by \eqref{e:endsoon?}
is invertible at the origin.  For this, we note that its Jacobian
matrix at the origin is a $(m+1)n$ lower triangular matrix by
blocks with $(m+1)$ blocks of size $n\times n$ in the diagonal.
Moreover, in view of \eqref{e:conditions++}, we obtain that this
matrix is of the following form
\begin{equation}\label{e:jacobi}
\left(\begin{array}{ccccc} \widetilde \lambda^1&0&0&\ldots&0\\
*&\widetilde \Lambda^1&0&\ldots&0\\
*&*& \widetilde \lambda^1&0&0\\
*&*&*&\ddots&0\\
*&*&*&*&\widetilde \Lambda^1
\end{array} \right)
\end{equation}
and hence is invertible. This shows that the map given by
\eqref{e:endsoon?} is a local biholomorphism. We may now apply
Proposition \ref{t:lineqn2} to conclude that there exists a
holomorphic map
\begin{equation}\label{e:gammamuend}
\Gamma_{m+1}^\mu
%=\Gamma_{m+1}^\mu
%(t^{[m+1]},\lambda^1,\Lambda^{(m+1)\kappa_0+r+1})
\colon
\C^{(m+1)n}\times {\mathcal G}^1_0(\CN)\times {\mathcal
G}_0^{(m+1)\kappa_0+r+1}(\CN)\to \C^n\times \C^{mn}
\end{equation} defined in a neighbourhood of $\{0\}\times
{\mathcal G}^1_0(\CN)\times {\mathcal
G}_0^{(m+1)\kappa_0+r+1}(\CN)$
 such that if $$(X,Y)\colon (\C^{(m+1)n},0)\to (\C^{(m+1)n},0)$$ 
 solves
 \eqref{e:system4}-\eqref{e:system5} for some $(\lambda^1,\Lambda^{(m+1)\kappa_0+r+1})\in
 {\mathcal G}^1_0(\CN)\times {\mathcal G}_0^{(m+1)\kappa_0+r+r}(\CN)$, then
$(X,Y)=(X\left(t^{[m+1]}\right),Y(t^{[m+1]})=\Gamma_{m+1}^\mu
(t^{[m+1]},\lambda^1,\Lambda^{(m+1)\kappa_0+r+1})$. We write
$$\Gamma_{m+1}^\mu:=(\Gamma_{m+1}^{\mu,1},\Gamma_{m+1}^{\mu,2},
\ldots, \Gamma_{m+1}^{\mu,m+1})\in \Cn\times \Cn \times \ldots
\times \Cn$$ and set
\begin{equation}\label{e:tired3}
\widetilde
\Psi_{m+1}^\delta(t^{[m+1]},\lambda^1,\Lambda^{(m+1)\kappa_0+r+1}):=\Gamma_{m+1}^{\mu,1}
(t^{[m+1]},\lambda^1,\Lambda^{(m+1)\kappa_0+r+1}),
\end{equation}
\begin{multline}
\label{e:tired4}
\widehat
\Psi_{m+1}^\delta(t^{[m+1]},\lambda^1,\Lambda^{(m+1)\kappa_0+r+1}):= \\
\widetilde
\Psi_{m+1}^{\delta}\cdot Q_z(\widetilde
\Psi_{m+1}^0,\overline{\Psi_m^0})+
W_\mu
\left(v^{m+1},\bar{v}^m,\Psi_{m+1}^0,\overline{\Psi^0_{m}},
%(\Psi_{m+1}^{\beta})_{1\leq |\beta|\leq r},
\Psi_{m+1}^{(r)},
%(\overline{\Psi_m^{\beta}})_{1\leq |\beta|\leq r+1}
\overline{\Psi_m^{(r+1)}} \right),
\end{multline}
where $W_\mu$ is the $\Cd$-valued polynomial map given by Lemma
\ref{l:gderivative}. In view of the above construction and Lemma
\ref{l:gderivative}, the map $\Psi_{m+1}^\delta$ is the desired
parametrization for the $\delta$-th derivative of all elements of
$\autM$. It remains to check that its components belong to the
space ${\mathcal E}_{(m+1)n}^{(m+1)\kappa_0+r+1}$. In view of
Lemma \ref{l:heavystuff} (i), \eqref{e:tired3}, \eqref{e:tired4}
and the induction assumption, it is enough to check that the
components of the map $\Gamma_{m+1}^\mu$ given by
\eqref{e:gammamuend} belong to the above mentioned space. For this
we write the Taylor expansion
\begin{equation}\label{e:bernie}
\Gamma_{m+1}^\mu(t^{[m+1]},\lambda^1,\Lambda^{(m+1)\kappa_0+r+1}):=\sum_{\gamma
\in
\N^{(m+1)n}}\Gamma_{m+1,\gamma}^\mu(\lambda^1,\Lambda^{(m+1)\kappa_0+r+1})\left(t^{[m+1]}\right)^{\gamma}.
\end{equation}
From Proposition \ref{t:lineqn2} we know that
$\Gamma_{m+1,\gamma}^\mu(\lambda^1,\Lambda^{(m+1)\kappa_0+r+1})$
may be written in the following form
\begin{equation}\label{e:form1}
  \frac{B_{\gamma}
  \left[\left((\partial_{t^{[m+1]}}^{\nu}\phi_{\mu,\alpha})
  (0,\lambda^1,\Lambda^{(m+1)\kappa_0+r+1})
  ,
  (\partial_{t^{[k]}}^{\nu}\widetilde
  \Psi_{k}^0)
  %(0,\lambda^1,\Lambda^{k\kappa_0})
  ,
  (\partial_{t^{[k]}}^{\nu}\overline{\widetilde
  \Psi_k^0})
  %(0,\lambda^1,\Lambda^{k\kappa_0})
  \right)_{\substack{1\leq |\alpha|\leq
  \kappa_0 \\ |\nu|\leq d_\gamma;1\leq k\leq m+1 }}\right]}{\left[{\det}\,
  \left(\displaystyle \frac{\partial {\mathcal I}}{\partial
  t^{[m+1]}}(0,\lambda^1,\Lambda^{(m+1)\kappa_0+r+1})\right)\right]^{c_\gamma}},
\end{equation}
 where $B_\gamma$ is a $\C^{(m+1)n}$-valued polynomial map of its arguments,
 $c_\gamma$ and $d_\gamma$ are positive integers, and
 the last two arguments of $B_\gamma$ are evaluated at
 $(0,\lambda^1,\Lambda^{k\kappa_0})$. Since each component
of the $\phi_{\mu,\alpha}$, $\widetilde \Psi_k^0$ and
$\overline{\widetilde \Psi_k^0}$ for $1\leq |\alpha|\leq
\kappa_0$, $1\leq k\leq m+1$ belongs to the space ${\mathcal
E}_{(m+1)n}^{(m+1)\kappa_0+r+1}$ by the above construction and the
induction assumption and since the Jacobian matrix of ${\mathcal I
}$ is given by \eqref{e:jacobi}, it is clear that each component
of $\Gamma_{m+1}^\mu$ belongs to the space ${\mathcal
E}_{(m+1)n}^{(m+1)\kappa_0+r+1}$. This completes the induction
part of the proof of Proposition \ref{p:propagation} and hence the
proof of that proposition.

%%%%%%%%%%%%%%%%%%%%%%%%%%%%%%%%%%%%%%%%%%%%%%%%%%%%%%%%%%%%%%%%%%%%%%%%

\section{Proofs of Theorems \ref{t:further1},
\ref{t:further2}, \ref{t:main1}, \ref{t:main3}, \ref{t:main2},
\ref{t:main1hsf}, \ref{t:main3hsf} }\label{s:conclusion} In this
section, we keep the notation introduced in \S \ref{s:jetpar}.

\begin{proof}[Proof of Theorem {\rm \ref{t:further1}}]
By the minimality criterion of Baouendi, Ebenfelt and Rothschild
(see e.g.\ \cite{BER4}), there exists $1\leq k_1\leq d+1$ (where
$d$ is the codimension of $M$ in $\CN$) such that the generic rank
of the Segre set map $v^{k_1}$ is $N$. Applying firstly
Proposition \ref{p:propagation} with $k=2k_1$, we obtain that
there exists a map $\Psi^0_{2k_1}\colon \C^{2nk_1}\times {\mathcal
G}_0^1(\CN)\times {\mathcal G}^{2k_1\kappa_0}_0(\CN)\to \CN$
holomorphic in a neighbourhood of $\{0\}\times {\mathcal
G}_0^1(\CN)\times {\mathcal G}_0^{2k_1\kappa_0}(\CN)$ whose
components belong to the space ${\mathcal
E}_{2k_1n}^{2k_1\kappa_0}$ such that for every $H\in \autM$ and
for every $t^{[2k_1]}$ sufficiently close to the origin in
$\C^{2k_1n}$ the following identity holds
$$(H\circ
v^{2k_1})(t^{[2k_1]})=\Psi^0_{2k_1}(t^{[2k_1]},j^1_0\overline{H},j_0^{2k_1\kappa_0}H).$$
Here we recall that $\kappa_0=\kappa (M,0)$. By mimicing the
proofs of \cite[\S 4.1]{BER4}, one may construct a map
$\Psi^0\colon \CN\times {\mathcal G}_0^1(\CN)\times {\mathcal
G}^{2k_1\kappa_0}_0(\CN)\to \CN$ holomorphic in a neighbourhood of
$\{0\}\times  {\mathcal G}_0^1(\CN)\times {\mathcal
G}^{2k_1\kappa_0}_0(\CN)$ whose components belong to the space
${\mathcal E}_{N}^{2k_1\kappa_0}$ such that for every $H\in \autM$
and for every $Z$ sufficiently close to the origin in $\C^{N}$ the
following identity holds
$$H(Z)=\Psi^0(Z,j^1_0\overline{H},j_0^{2k_1\kappa_0}H).$$
Now the map
$\psi(Z,\Lambda^{2k_1\kappa_0}):=\Psi^0(Z,\overline{\Lambda^1},\Lambda^{2k_1\kappa_0})$
satisfies all the conclusions of Theorem \ref{t:further1} except
that $\psi$ is a function of $\Lambda^{2k_1\kappa_0}\in
G_0^{2(d+1)\kappa_0}(\CN)$, which is twice more than claimed. To
complete the proof of the theorem, we apply again Proposition
\ref{p:propagation} with $k=k_1$, which may be assumed to be even
(the other case being analogous) and obtain the existence of a map
$\Psi^0_{k_1}\colon \C^{nk_1}\times {\mathcal G}_0^1(\CN)\times
{\mathcal G}^{k_1\kappa_0}_0(\CN)\to \CN$ holomorphic in a
neighbourhood of $\{0\}\times {\mathcal G}_0^1(\Cn)\times
{\mathcal G}_0^{k_1\kappa_0}(\CN)$ whose components belong to the
space ${\mathcal E}_{k_1n}^{k_1\kappa_0}$ such that for every
$H\in \autM$ and for every $t^{[k_1]}$ sufficiently close to the
origin in $\C^{k_1n}$ the following identity holds
\begin{equation}\label{e:holds3}
(H\circ
v^{k_1})(t^{[k_1]})=\Psi^0_{k_1}(t^{[k_1]},j^1_0\overline{H},j_0^{k_1\kappa_0}H).
\end{equation}
By mimicing the arguments of \cite[\S 4.2]{BER4}, one may
construct by using \eqref{e:holds3} and the above map $\Psi^0$,
another map $\Phi\colon \CN\times {\mathcal G}_0^1(\CN)\times
{\mathcal G}^{k_1\kappa_0}_0(\CN)\to \CN$ holomorphic in a
neighbourhood of $\{0\}\times  {\mathcal G}_0^1(\CN)\times
{\mathcal G}^{k_1\kappa_0}_0(\CN)$ whose components belong to the
space ${\mathcal E}_{N}^{k_1\kappa_0}$ such that for every $H\in
\autM$ and for every $Z$ sufficiently close to the origin in
$\C^{N}$ the following identity holds
$$H(Z)=\Phi(Z,j^1_0\overline{H},j_0^{k_1\kappa_0}H);$$
the details are left to the reader. Then setting
$\Psi(Z,\Lambda^{k_1\kappa_0}):=\Phi
(Z,\overline{\Lambda^1},\Lambda^{k_1\kappa_0})$ yields the
required parametrization. The proof of Theorem \ref{t:further1} is
complete.
\end{proof}

\begin{proof}[Proof of Theorem {\rm \ref{t:further2}}]
Theorem \ref{t:further2} follows directly from Theorem
\ref{t:further1} and from repeating the arguments of \cite[\S
4.3]{BER4}. We leave the details to the reader.
\end{proof}

\begin{proof}[Proofs of Theorem {\rm \ref{t:main1}} and Theorem {\rm \ref{t:main3}}] Firstly, note that to prove Theorem~\ref{t:main1}, it is enough to prove the parametrization property for the stability group of $(M,p)$ for every $p\in M$. Then Theorem \ref{t:main1} (resp.\ Theorem \ref{t:main3}) follows
immediately from Theorem \ref{t:further1} (resp.\ Theorem
\ref{t:further2}), Proposition \ref{t:prop1} and the
upper-semicontinuity of the map $M\ni q\mapsto \kappa_M(q)$ (see
\S \ref{ss:nondegeneracy}).
\end{proof}

\begin{proof}[Proof of Theorem {\rm \ref{t:main2}}]
Firstly, when $M$ is generic, Theorem \ref{t:main2} follows from
the fact that compact real-analytic CR submanifolds of $\CN$ do
not contain any complex-analytic subvariety of positive dimension
(\cite{DF2}) and hence are essentially finite at all points, from
the upper-semicontinuity of the map $p\mapsto \ell_p$ in Theorem
\ref{t:main3} and from the reflection principle proved in
\cite{BJT}. When $M$ is not generic, one proceeds as follows. By
\cite{BERbook}, there exists $r\in \{1,\ldots, N-1\}$ such that
for every point $p\in M$ there exists a neighbourhood $U_p$ of $p$
in $\CN$ and a neighbourhood $W_p$ of $0$ in $\CN$ such $M\cap
U_p$ is biholomorphic to $W_p\cap (\widehat M_p\times \{0\})$
where $\widehat M_p$ is a real-analytic generic submanifold of
$\C^r$ containing the origin, minimal at each of its points and
not containing any complex-analytic subvariety of positive
dimension. By using Theorem \ref{t:main3}, for every $p\in M$ we
may assume, shrinking $\widehat M _p$ near the origin if
necessary, that there exists a positive integer $k_p$ so that
germs at any point $q\in \widehat M_p$ of smooth CR
diffeomorphisms mapping $\widehat M_p$ into another real-analytic
generic submanifold of $\C^r$ of the same CR dimension as that of
$M$ are uniquely determined by their $k_p$-jets at $q$. Then the
theorem easily follows by using a finite covering of $M$ by the
the neighbourhoods $U_p$ for $p\in M$. The proof of Theorem
\ref{t:main2} is therefore complete.
\end{proof}

\begin{proof}[Proofs of Theorems {\rm \ref{t:main1hsf}} and {\rm \ref{t:main3hsf}}]
Since real-analytic hypersurfaces which contain no complex-analytic
subvariety of positive dimension are automatically minimal (see
e.\ g.\ \cite{BERbook}), Theorem \ref{t:main1hsf} (resp.\ Theorem
\ref{t:main3hsf}) follows immediately from Theorem \ref{t:main1}
(resp.\ from Theorem \ref{t:main3}).
\end{proof}

\bibliographystyle{amsplain}
\bibliography{bibliography_oct04}
\end{document}